\documentclass[11pt]{article}

\usepackage{amsmath}
\usepackage{amssymb}
\usepackage{latexsym}
\usepackage{enumitem}
\usepackage{mathrsfs}
\usepackage{comment}
\usepackage{color}

\newtheorem{assumption}{Assumption}

\def\qed{ \ \vrule width.2cm height.2cm depth0cm\smallskip}
\newenvironment{proof}{\noindent {\bf Proof.\/}}{$\qed$\vskip 0.1in}

\newcommand{\la}{\langle}
\newcommand{\ra}{\rangle}

\newcommand{\ol}{\overline}
\newcommand{\ul}{\underline}
\newcommand{\eps}{\varepsilon}

\newcommand{\ba}{\begin{array}}
\newcommand{\ea}{\end{array}}
\newcommand{\be}{\begin{equation}}
\newcommand{\ee}{\end{equation}}
\newcommand{\bea}{\begin{eqnarray}}
\newcommand{\eea}{\end{eqnarray}}
\newcommand{\beaa}{\begin{eqnarray*}}
\newcommand{\eeaa}{\end{eqnarray*}}

\def\dbD{\mathbb{D}}
\def\dbE{\mathbb{E}}
\def\dbF{\mathbb{F}}

\def\dbL{\mathbb{L}}

\def\dbP{\mathbb{P}}
\def\dbR{\mathbb{R}}

\def\dbT{\mathbb{T}}

\def\sL{\mathscr{L}}

\def\a{\alpha}

\def\g{\gamma}
\def\d{\delta}
\def\e{\varepsilon}

\def\k{\kappa}
\def\l{\lambda}
\def\m{\mu}

\def\si{\sigma}

\def\f{\varphi}
\def\th{\theta}
\def\o{\omega}

\def\G{\Gamma}
\def\D{\Delta}
\def\Th{\Theta}

\def\O{\Omega}
\def\cA{{\cal A}}
\def\cB{{\cal B}}

\def\cE{{\cal E}}
\def\cF{{\cal F}}

\def\cJ{{\cal J}}

\def\cL{{\cal L}}

\def\cP{{\cal P}}

\def\cU{{\cal U}}

\def\cW{{\cal W}}

\def\no{\noindent}

\def\ms{\medskip}
\def\bs{\bigskip}
\def\q{\quad}
\def\qq{\qquad}

\def\pa{\partial}
\def\cd{\cdot}
\def\cds{\cdots}

\def\qed{ \hfill \vrule width.25cm height.25cm depth0cm\smallskip}

\newcommand{\basa}{\begin{assumption}}
\newcommand{\easa}{\end{assumption}}

\newcommand{\bas}{\begin{assum}}
\newcommand{\eas}{\end{assum}}

\def\limsup{\mathop{\overline{\rm lim}}}
\def\liminf{\mathop{\underline{\rm lim}}}

\def\pa{\partial}

 \def\cd{\cdot}
\def\cds{\cdots}

\def\supp{\hbox{\rm supp$\,$}}

\def\dis{\displaystyle}

\def\wh{\widehat}

\def\cad{c\`{a}dl\`{a}g}

\def\1{{\bf 1}}

\def\:{\!:\!}
\def\reff#1{{\rm(\ref{#1})}}
\def \proof{{\noindent \bf Proof\quad}}

at 9pt

\begin{document}

\newtheorem{thm}{Theorem}[section]
\newtheorem{lem}[thm]{Lemma}
\newtheorem{cor}[thm]{Corollary}
\newtheorem{prop}[thm]{Proposition}
\newtheorem{rem}[thm]{Remark}
\newtheorem{eg}[thm]{Example}
\newtheorem{defn}[thm]{Definition}
\newtheorem{assum}[thm]{Assumption}

\renewcommand {\theequation}{\arabic{section}.\arabic{equation}}
\def\thesection{\arabic{section}}

\title{ Viscosity Solutions to Parabolic Master Equations and McKean-Vlasov SDEs with Closed-loop Controls}
\author{Cong Wu\thanks{\noindent  Wells Fargo Securities, San Francisco, CA 94105. Email: wucong085@gmail.com. } ~ and ~ Jianfeng Zhang\thanks{\noindent  Department of Mathematics,
University of Southern California, Los Angeles, CA 90089. E-mail:
jianfenz@usc.edu. This author is supported in part by NSF grant  \#1413717. }  
}
\date{\today}
\maketitle

\begin{abstract}
The master equation is a type of PDE whose state variable involves the distribution of certain underlying state process.  It is a powerful tool for studying the limit behavior of large interacting systems, including mean field games and systemic risk. It also appears naturally in stochastic control problems with partial information and in time inconsistent problems. In this paper we propose a novel notion of viscosity solution for parabolic master equations, arising mainly from control problems, and establish its wellposedness. Our main innovation is to restrict the involved measures to certain set of semimartingale measures which satisfy the desired compactness. As an important example, we study the HJB master equation associated with the control problems for McKean-Vlasov SDEs. Due to practical considerations, we consider closed-loop controls. It turns out that the regularity of the value function becomes much more involved in this framework than the counterpart in the standard control problems. Finally, we build the whole theory in the path dependent setting, which is often seen in applications. The main result in this part is an extension of Dupire \cite{Dupire}'s functional It\^{o} formula. This It\^{o} formula requires a special structure of the derivatives with respect to the measures, which was originally due to Lions \cite{Lions4} in the state dependent case. We provided an elementary proof for this well known result in the short note \cite{WZ}, and the same arguments work in the path dependent setting here. 
\end{abstract}

\no{\bf Keywords.}  Master equation, McKean-Vlasov SDEs, viscosity solution, functional It\^{o} formula, path dependent PDEs, Wasserstein spaces, dynamic programming principle

\bs

\no{\it 2000 AMS Mathematics subject classification:}  35K55,  49L25, 60H30,	35R15,  49L20,  93E20

\vfill\eject
\section{Introduction}
\label{sect-Introduction}
\setcounter{equation}{0} 
Initiated independently by Caines, Huang, \& Malhame \cite{CHM} and Lasry \& Lions \cite{LL}, mean field games and the closely related mean field  control problems have received very strong attention in the past decade. Such problems consider the limit behavior of large systems where the agents interact with each other in certain symmetric way, with the systemic risk as a notable application. There have been numerous publications on the subject, see e.g. Cardaliaguet \cite{Cardaliaguet},  Bensoussan, Frehse, \& Yam \cite{BFY}, Carmona \& Delarue \cite{CD1, CD2}, and the references therein. The master equation is a powerful and inevitable tool in this framework, which plays the role of the PDE in the standard literature of controls and games. The main feature of the master equation is that its state variable contains  probability measures, typically the distribution of certain underlying state process, so it can be viewed as a PDE on the Wasserstein space. By nature this is an infinite dimensional problem.  The master equation is also a convenient tool for (standard) control problems with partial information, see e.g. Bandini,  Cosso,  Fuhrman, \& Pham \cite{BCFP1, BCFP2} and Saporito \& Zhang \cite{SZ}, and for some time inconsistent problems as we will see in this paper.

Our main goal of this paper is to propose an intrinsic notion of viscosity solutions for  parabolic master equations which  mainly arise from control problems or zero-sum game problems in the McKean-Vlasov setting. There have been serious efforts on classical solutions for master equations in various settings,   see e.g. Buckdahn, Li,  Peng, \& Rainer \cite{BLPR}, Cardaliaguet,  Delarue, Lasry, \& Lions \cite{CDLL},  Chassagneux, Crisan, \& Delarue \cite{CCD}, Saporito \& Zhang \cite{SZ}, and Bensoussan,  Graber, \& Yam \cite{BGY}. However, due to its infinite dimensionality, all these works require very strong technical conditions. So there is a cry for an appropriate notion of weak solutions. {\color{black}We remark that a classical solution requires the candidate solution (typically the value function of certain control/game problem) to be  in $C^{1,2}$ (in appropriate sense), while a viscosity solution theory will allow us to reduce the regularity requirement  to $C^0$.  It is in general very challenging to establish the differentiability of the value function (especially that with respect to the measures), so such a relaxation of regularity requirement is desirable in many applications. }

There have already been some works on viscosity solutions. A natural approach is to use smooth test functions on the  Wasserstein space, see e.g. Carmona \& Delarue \cite{CD}. However, the involved space lacks the local compactness, which is crucial for the viscosity theory, and thus the comparison principle does not seem possible  in this approach. In an alternative approach Pham \& Wei \cite{PW}  lift the functions on the Wasserstein space to those on the Hilbert space of random variables and then apply the existing viscosity theory on Hilbert spaces, see e.g. Lions \cite{Lions1, Lions2, Lions3} and  Fabbri, Gozzi, \& Swiech \cite{FGS}. Along this approach one could obtain both existence and uniqueness. However, this notion is not intrinsic, in particular, it is not clear to us that a classical solution (with smoothness in the Wasserstein space of probability measures instead of the Hilbert space of random variables) would be  a viscosity solution in their sense. Moreover,  the viscosity theory on Hilbert spaces is not available in the path dependent case (see Ren \& Rosestolato \cite{RR} for some recent progress along this direction though), and thus it will be difficult to extend their results to the path dependent case which is important in applications and is another major goal of this paper.  We remark that we are in the stochastic setting and thus the master equation is of second order (in certain sense, see Remark \ref{rem-C12hat}). There are several works for first order master equations corresponding to the deterministic setting, see e.g.    Gangbo \& Swiech \cite{GS1, GS2} and  Bensoussan \& Yam \cite{BY}. 

We shall propose a new notion of viscosity solutions, motivated from our previous works  Ekren, Keller, Touzi, \& Zhang \cite{EKTZ} and Ekren, Touzi, \& Zhang \cite{ETZ1, ETZ2} for viscosity solutions of path dependent PDEs. Our main innovation is to modify the set of test functions so as to ensure certain desired compactness. To be precise, let $V(t,\mu)$ be a candidate solution, where $\mu$ is a probability measure, and $\f$ be a smooth (in certain sense) test function at $(t, \mu)$, we shall require  $[\f - V]$ achieves maximum/minimum at $(t,\mu)$ only over the set $[t, t+\d] \times \cP_L(t,\mu)$, where $\cP_L(t, \mu)$ is a compact set of semimartingale measures with drift and diffusion characteristics bounded by a constant  $L$. We note that, if  we replace the above $\cP_L(t, \mu)$ with the $\d$-neighborhood of $\mu$ under  the Wasserstein distance, as in \cite{CD}, then the latter set is not compact under the Wasserstein distance and we will encounter serious difficulties for establishing the comparison principle. We should also note that, if the underlying state space (on which the probability measures are defined) is a torus  $\dbT^d$ instead of $\dbR^d$, then in the state dependent case the $\d$-neighborhood of $\mu$ under the Wasserstein distance is compact and thus the theory is quite hopeful. However, for the applications in our mind it is more natural to consider $\dbR^d$ as the underlying state space, and in the mean time we are interested in the path dependent case for which the $\d$-neighborhood wouldn't work for the torus either.

Our choice of $\cP_L(t,\mu)$ is large enough so that, in many applications we are interested in,  the value function will be a viscosity solution to the corresponding master equation. On the other hand, the  compactness of $\cP_L(t,\mu)$ enables us to establish the basic properties of viscosity solutions following rather standard  arguments:  consistency with classical solutions, equivalence to the alternative definition through semi-jets, stability, and partial comparison principle. The comparison principle is of course the main challenge. We nevertheless establish some partial results  in the general case and prove the full comparison principle completely in some special cases. To our best knowledge this is the first uniqueness result in the literature for an intrinsic notion of viscosity solutions for second order master equations.

As far as we know, all works on master equations in the existing literature consider only the state dependent case, where the measures are defined on the finite dimensional space $\dbR^d$ (or the torus $\dbT^d$). However, in many applications the problem can be path dependent, for example,  lookback options, variance swap, rough volatility, delayed SDEs, to mention a few. In particular,  Saporito \& Zhang \cite{SZ} studied control problems with information delay, which naturally induces a path dependent master equation. The second goal of this paper is to establish the whole theory in the path dependent setting, namely the involved probability measure $\mu$ is the distribution of the stopped underlying process $X_{\cd\wedge t}$, rather than the distribution of the current state $X_t$. The main result in this regard is a functional It\^{o} formula in the  McKean-Vlasov setting, extending the well known result of Dupire \cite{Dupire} in the standard setting. To establish this, we require a special structure of  the path derivative with respect to the measure, see \reff{pamuf} below. In the state dependent case, such structure was established by Lions \cite{Lions4}, see also Cardaliaguet \cite{Cardaliaguet} and Gangbo \& Tudorascu \cite{GT}, by using quite advanced tools. We provided an elementary proof for this well known result, which was reported separately in the short note \cite{WZ}, and the same arguments work well in our path dependent framework here. We emphasize that, while this paper is in the path dependent setting, our results on viscosity solutions of master equations are new even in the state dependent case. 

Our third goal is to study McKean-Vlasov SDEs with closed-loop controls, whose value function is a viscosity solution to the HJB type master equation. We note that in many applications  closed-loop controls (i.e. the control depends on the state process) are more appropriate than open-loop controls (i.e. the control depends on the noise),  especially when games are considered, see e.g. Zhang \cite{Zhang} Section 9.1 for detailed discussions. For McKean-Vlasov SDEs, the two types of controls have very subtle difference even for control problems (and more subtle for games), and under closed-loop controls, the regularity of the value function becomes rather technical. By choosing the  admissible controls carefully and by using some sophisticated approximations, we manage to prove the desired regularity and then verify the viscosity solution  property.  Again, while we are in the path dependent setting, our result is new even in the state dependent case, and we believe our approximations will be quite useful for more thorough analysis on functions of probability measure. 

Finally, we emphasize that our master equation is parabolic, which mainly corresponds to control problems or zero-sum game problems in the McKean-Vlasov setting, and the solution takes the form $V(t, \mu)$.  {\color{black} The master equation induced by mean field games involves functions in the form $V(t,x,\mu)$, and in the path dependent setting this becomes $V(t,\o, \mu)$. The two types of equations have some fundamental differences.   On one hand, our master equation could be nonlinear in $\pa_\mu V$, the derivative of $V$ with respect to the probability measure $\mu$, while mean field game master equation is typically linear in $\pa_\mu V$ (but could be nonlinear in $\pa_x V$). On the other hand, mean field game master equation is non-local in $\pa_x V$, which destroys certain crucial monotonicity property and thus the comparison principle does not hold (even for classical solutions). In fact, due to these differences, in many works master equations refer only to the equations arising from mean field games, while those from mean field control problems are called HJB equations in Wasserstein space. We nevertheless call both master equations, since they share many properties and require similar technical tools. So this paper studies mean field control  master equations, and we refer to the recent work Mou \& Zhang \cite{MZ} for weak solutions (instead of viscosity solutions) to mean field game master equations. }

The rest of the paper is organized as follows. In Section \ref{sect-Ito} we establish the functional It\^{o} calculus in the Wasserstein space. In Section \ref{sect-classical} we introduce parabolic master equations and present several examples, which in particular show some applications of master equations. In Section \ref{sect-Viscosity} we introduce our notion of viscosity solutions and establish its wellposedness. In Section \ref{sect-HJB} we study the McKean-Vlasov SDE with closed-loop controls and show its value function is a viscosity solution to the HJB master equation.   

\section{Functional It\^o calculus in the Wasserstein space}
\label{sect-Ito}
\setcounter{equation}{0}

\subsection{A brief overview in the state dependent setting}
\label{sect-state}
We first recall the Wasserstein metric on the space of probability measures. Let $(\O, \cF)$ be an arbitrary measurable space equipped with a metric $\|\cd\|$. For any probability measures $\mu, \nu$ on $\cF$,  let $\cP(\mu, \nu)$ denote the space of probability measures $\dbP_{\mu,\nu}$ on the product space $(\O\times \O, \cF\times \cF)$ with marginal measures $\mu$ and $\nu$. Then the $2$-Wasserstein distance  of $\mu$ and $\nu$ is defined as (assuming $(\O, \cF)$ is rich enough):
\bea
\label{cW21}
\cW_2(\mu,\nu):=\inf_{\dbP_{\mu, \nu}\in\cP(\mu,\nu)}\Big(\int_{\O\times \O} \|\o_1-\o_2\|^2 d\dbP_{\mu,\nu}(\o_1, \o_2)\Big)^{\frac{1}{2}}.
\eea

In the state dependent setting, one may set the measurable space as $\big(\dbR^d, \cB(\dbR^d)\big)$ (or the torus $\big(\dbT^d, \cB(\dbT^d)\big)$ as in some works). Let $\cP_2(\dbR^d)$ denote the set of square integrable measures on $\big(\dbR^d, \cB(\dbR^d)\big)$, equipped with the metric $\cW_2$. For an arbitrary probability space $(\O, \cF, \dbP)$, let $\dbL^2(\cF, \dbP)$ denote the Hilbert space of $\dbP$-square integrable $\cF$-measurable $\dbR^d$-valued random variables. Given a function $f: \cP_2(\dbR^d)\to \dbR$,  we may lift $f$ to a function on $\dbL^2(\cF, \dbP)$: $F(\xi) := f(\dbP_\xi)$, where $\dbP_\xi$ is the $\dbP$-distribution of $\xi\in \dbL^2(\cF,\dbP)$. Assume $F$ is continuously Fr\'{e}chet differentiable, Lions \cite{Lions4} showed that the Fr\'{e}chet derivative $D F$ takes the following form: for some deterministic function $h: \cP_2(\dbR^d) \times  \dbR^d\to \dbR^d$,
\bea
\label{stateFrechet}
D F (\xi) = h(\dbP_\xi, \xi),
\eea
see also Cardaliaguet \cite{Cardaliaguet},  Gangbo \& Tudorascu \cite{GT}, and Wu \& Zhang \cite{WZ}.  Thus naturally we may define $\pa_\mu f := h$.  Note that $\pa_\mu f$ is essentially equivalent to the Wasserstein gradient in the optimal transportation theory, see e.g. Carmona \& Delarue \cite{CD1}. Assume further that $\pa_\mu f$ is continuously differentiable with respect to the second variable $x$, then we have the following It\^{o} formula, due to Buckdahn, Li, Peng, \& Rainer \cite{BLPR} and Chassagneux, Crisan, \&  Delarue \cite{CCD},
\bea
\label{stateIto}
f(\dbP_{X_t}) = f(\dbP_{X_0}) + \dbE^\dbP\Big[\int_0^t \pa_\mu f(\dbP_{X_s}, X_s) \cd dX_s + {1\over 2} \int_0^t \pa_x \pa_\mu f (\dbP_{X_s}, X_s) : d\la X\ra_s\Big],
\eea
for any $\dbP$-semimartingale $X$ satisfying certain technical conditions, where $\cd$ and $:$ denote inner product and trace, respectively.  

Our goal of this section is to extend both \reff{stateFrechet} and \reff{stateIto} to the path dependent setting. {\color{black} We remark that path dependence appears naturally in many applications. For example, in option pricing theory, many exotic options like lookback options and Asian options are path dependent, then their prices would satisfy certain path dependent PDEs.  Another interesting example is the rough volatility model, where the state process is non-Markovian and  a path dependent PDE is induced naturally even in state dependent models, see Viens \& Zhang \cite{VZ}.  All these models will naturally lead to path dependent master equations when extended to the mean field framework. A more interesting example is the stochastic optimization  in standard framework but with constant controls, where a state dependent model will naturally induce a path dependent master equation, see Theorem \ref{thm-deterministic2}   below.   }

Throughout the paper, for an arbitrary process $X$, we introduce the notation:
\bea
\label{Xst}
X_{s,t} := X_t - X_s,\q 0\le s\le t\le T.
\eea

\subsection{The canonical setup in the path dependent setting}
\label{sect-canonical}
Throughout this paper, we shall fix the canonical space  $\O:=C([0,T],\dbR^d)$, equipped with the uniform norm $\|\o\|:=\sup_{t\in[0,T]}|\o_t|$. Let $X$ denote the canonical process, namely $X_t(\o):=\o_t$, $\dbF:= \{\cF_t\}_{0\le t\le T} := \dbF^X$ the natural filtration generated by $X$,  $\cP_2$ the set of probability measures $\mu$ on $(\O,\cF_T)$ such that $\dbE^\mu[\|X\|^2] < \infty$, equipped with the Wasserstein distance $\cW_2$ defined by \reff{cW21}. Note that $(\O, \|\cd\|)$ and $(\cP_2, \cW_2)$ are Polish spaces, namely they are complete and separable. We may also use the notation $\dbP$ to denote probability measures. Quite often we shall use $\mu$ when viewing it as a variable of functions,  and use $\dbP$ when considering the distribution of some random variables or processes.  Moreover, given a random variable or a stochastic process $\xi$ under certain probability measure $\dbP$,   we also use $\dbP_\xi := \dbP\circ \xi^{-1}$ to denote its distribution under $\dbP$. When the measure $\dbP$ is clear from the context, we may also use the notation $\cL_\xi := \dbP_\xi$.
 
The state space of our master equation is $\Th:=[0,T]\times \cP_2$.  For each $(t,\mu)\in\Th$, let $\mu_{[0,t]}\in \cP_2$ be the distribution of the stopped process $X_{t\wedge \cd}$ under $\mu$. Since $\cF_T^{X_{t\wedge \cd}}=\cF_t$, $\mu_{[0,t]}$ is completely determined by the restriction of $\mu$ on $\cF_t$. For $(t,\mu),(t', \mu')\in\Th$, by abusing the notation $\cW_2$ we define the $2$-Wasserstein pseudometric on $\Th$ as
\bea
\label{cW22}
\cW_2((t,\mu),(t',\mu')):=\Big(|t-t'|+\cW^2_2\big(\mu_{[0,t]},\mu'_{[0,t']}\big)\Big)^{\frac{1}{2}}.
\eea
If a function $f:\Th\to\dbR$ is Borel measurable, with respect to the topology induced by $\cW_2$, then it must be $\dbF$-adapted in the sense that $f(t,\mu)=f(t,\mu_{[0,t]})$ for any $(t,\mu)\in \Th$. In particular, if $f$ is continuous, then it is $\dbF$-adapted. Moreover, for $(t,\mu)\in \Th$, let $\mu_t := \mu \circ X_t^{-1}$ denote the distribution of the random variable $X_t$. We say $f$ is state dependent if $f(t, \mu)$ depends only on $\mu_t$, and in this case we may abuse the notation and denote $f(t,\mu_t) = f(t, \mu)$.

{\color{black}In order to  establish the functional It\^{o} formula on $\Th$, as in Dupire \cite{Dupire} we extend the canonical space to the {\cad} space $\wh \O:= \dbD([0, T), \dbR^d)$ (we use $\wh\cd$ to denote the extensions to the {\cad} space), equipped with the Skorohod distance:
\bea
\label{Skorohod}
d_{SK}(\wh \o, \wh \o') := \inf_\l \sup_{0\le t\le T} \big[|t-\l(t)| + |\wh\o_t - \wh \o'_{\l(t)}|\big]
\eea
where $\l: [0, T]\to [0, T]$ is continuous, strictly increasing, with $\l(0) = 0$ and $\l(T) = T$.  Extend the notations $\wh X$, $\wh \dbF$, $\wh \cP_2$, $\wh \Th$, as well as the $2$-Wasserstein pseudometric on $\wh \Th$  in an obvious way, in particular, in \reff{cW21} the metric $\|\o^1 - \o^2\|$ should be replaced with $d_{SK}(\wh \o^1, \wh \o^2)$.  Then $(\wh\O, d_{SK})$ and $(\wh\cP_2, \cW_2)$ are also Polish spaces.
}

\subsection{Pathwise derivatives in the Wasserstein  space}

Let $f: \wh \Th \to \dbR$  be continuous (and thus $\wh\dbF$-adapted). We define its time derivative as:
\bea
\label{patf}
\pa_t f(t, \wh \mu) := \lim_{\d\downarrow 0} {f(t+\d, \wh\mu_{[0, t]}) - f(t,\wh\mu)\over \d},
\eea
provided the limit in the right side above exists. 

\begin{rem}
\label{rem-patf}
{\rm
The $\pa_t f$ in \reff{patf} is actually the right time derivative.  Due to the adaptedness requirement, similar to the pathwise analysis in Dupire \cite{Dupire}, the left time derivative  is not convenient to define. Nevertheless, for the theory which we will develop in the paper, in particular for the functional It\^{o} formula, the  right time derivative is sufficient. 
\qed}
\end{rem}
   
The spatial derivative is much more involved.  Consider an arbitrary atomless Polish probability space $(\tilde \O, \tilde \cF, \tilde \dbP)$. Let $\dbL^2(\tilde \O; \dbR^d)$ and  $\dbL^2(\tilde \O; \wh\O)$ denote the sets of $\tilde \dbP$-square integrable $\tilde \cF$-measurable mappings $\xi: \tilde \O\to\dbR^d$ and  $\tilde X: \tilde \O \to \wh \O$, respectively.  We first lift $f$ to a function $F: [0, T]\times  \dbL^2(\tilde \O; \wh \O) \to \dbR$:
\bea
\label{lift}
F(t, \tilde X) := f(t, \tilde\dbP_{\tilde X}) =  f(t, \tilde\dbP_{\tilde X_{t\wedge \cd}}).
\eea
We say $F$ is Fr\'{e}chet differentiable at $(t, \tilde X)$ with derivative $D F(t, \tilde X)\in \dbL^2(\tilde \O; \dbR^d)$ if 
\bea
\label{Frechet}
F(t, \tilde X + \xi \1_{[t, T]}) - F(t, \tilde X)  = \dbE^{\tilde \dbP}\big[ D F(t, \tilde X) \cd \xi\big] + o(\|\xi\|_2)\q\mbox{for all}~ \xi \in   \dbL^2(\tilde \O; \dbR^d),
\eea
where $\|\xi\|_2^2 := \dbE^{\tilde \dbP}[|\xi|^2]$.
In particular, this implies that $D F(t, \tilde X)$ is the G\^{a}teux derivative:
\bea
\label{Gateux}
\dbE^{\tilde \dbP}\big[ DF(t, \tilde X)\cd  \xi\big] = \lim_{\e\to 0} {F(t, \tilde X +  \e \xi \1_{[t, T]}) - F(t, \tilde X)  \over \e},\q \mbox{for all}~ \xi \in   \dbL^2(\tilde \O; \dbR^d).
\eea
We emphasize that the above derivative involves only the perturbation of  $\tilde X$ on $[t, T]$, but not on $[0, t)$. Moreover, since $f$ is $\wh\dbF$-adapted, so $D F(t, \tilde X)$ actually involves only the perturbation of $\tilde X$ at $t$.  Our main result in this subsection is:
\begin{thm}
\label{thm-pamu}  Let $f: \wh \Th \to \dbR$ be continuous.  Assume the lifted function $F$ defined by \reff{lift} is Fr\'{e}chet differentiable and $D F$ is continuous in the sense that
\bea
\label{DFcont}
 \lim_{n\to\infty} \dbE^{\tilde \dbP}\Big[|DF(t, \tilde X^n) - D F(t, \tilde X)|^2\Big]=0,\q\mbox{whenever} \q\lim_{n\to\infty} \dbE^{\tilde \dbP} \big[ d^2_{SK}(\tilde X^n, \tilde X)\big]=0.
 \eea
  Then there exists an $\wh\cF_t$-measurable function $\psi:  \wh \O \to \dbR^d$  such that 
  \bea
  \label{DFX}
  D F(t, \tilde X) = \psi(\tilde X_{t\wedge \cd}),\q \tilde \dbP\mbox{-a.s.}
  \eea
  Moreover, $\psi$ is determined by $f$ and $\tilde \dbP_{\tilde X}$, and is unique $\tilde \dbP_{\tilde X}$-a.s.  
\end{thm}
\proof The uniqueness of $\psi$ follows from \reff{DFX} and the uniqueness of the Fr\'{e}chet derivative $D F(t, \tilde X)$. Moreover, by the $\wh\dbF$-adaptedness of $f$, clearly $D F(t, \tilde X)$ is determined by $\tilde X_{t\wedge \cd}$, and thus so is $\psi$. We prove the rest of the theorem  in two steps.

{\it Step 1.} We first construct $\psi$ in the case that $\tilde X$ is discrete: there exist $\wh\o_i \in \wh \O$, $i\ge 1$, such that  $\sum_{i\ge 1} p_i = 1$, where $p_i := \tilde \dbP(\tilde X = \wh\o_i)>0$.  For any $x \in \dbR^d\backslash \{0\}$,  $E \subset E_i := \{\tilde X = \wh\o_i\}$, and $\e>0$,  denote $\wh\o^\e_i:= \wh\o_i + \e x \1_{[t, T]}$ and $\tilde X^\e:=  \tilde X +  \e  x\1_E  \1_{[t, T]}$. {\color{black}Note that, 
\beaa
\tilde X^\e(\tilde \o) = \sum_{j\neq i} \wh\o_j \1_{E_j}(\tilde \o) + \wh \o_i \1_{E_i \backslash E} (\tilde \o) + \wh \o^\e_i \1_E(\tilde \o),\q \tilde \o\in \tilde \O.
\eeaa
Then, denoting by $\d_\cd$ the Dirac-measure,
\beaa
\cL_{ \tilde X^\e} = \sum_{j\neq i} p_j \d_{\{\wh\o_j\}} + [p_i - \tilde \dbP(E)] \d_{\{\wh\o_i\}} + \tilde \dbP(E) \d_{\{\wh\o^\e_i\}}, 
\eeaa
} and thus
\beaa
&&\dbE^{\tilde \dbP}\big[D F(t, \tilde X) \cd x \1_E\big] = \lim_{\e\to 0}  {F(t, \tilde X +  \e  x\1_E  \1_{[t, T]}) - F(t, \tilde X)  \over \e}\\
&&= \lim_{\e\to 0} {f\big(t,  \sum_{j\neq i} p_j \d_{\{\wh\o_j\}} + [p_i - \tilde \dbP(E)] \d_{\{\wh\o_i\}} + \tilde \dbP(E) \d_{\{\wh\o^\e_i\}} \big) - f\big(t,  \sum_{j\ge 1} p_j \d_{\{\wh\o_j\}}  \big) \over \e}.
\eeaa
This implies that $\dbE^{\tilde \dbP}\big[D F(t, \tilde X) \cd x \1_E\big] =\dbE^{\tilde \dbP}\big[D F(t, \tilde X) \cd x \1_{E'}\big]$ for any $E, E'\subset E_i$ such that $\tilde \dbP(E)=\tilde \dbP(E')$. By  Wu \& Zhang \cite{WZ} Lemma 2,  we see that $D F(t, \tilde X) \cd x$ is a constant on $E_i$: by setting $E=E_i$,
\bea
\label{DFx}
D F(t, \tilde X) \cd x = \lim_{\e\to 0}  {f\big(t,  \sum_{j\neq i} p_j \d_{\{\wh\o_j\}} +  p_i \d_{\{\wh\o_i+  \e x \1_{[t, T]}\}} \big) - f\big(t,  \sum_{j\ge 1} p_j \d_{\{\wh\o_j\}}  \big) \over \e p_i}.
\eea
 Since $x\in \dbR^d$ is arbitrary, $D F(t, \tilde X) = y_i\in \dbR^d$, $\tilde \dbP$-a.s. on $E_i$. Clearly there exists a  Borel-measurable function $\psi: \wh \O\to \dbR^d$ such that $\psi( \wh\o_i) = y_i$, $i\ge 1$, and thus $D F(t, \tilde X) = \psi(\tilde X)$, $\tilde \dbP$-a.s.  Note that $\psi$ is unique in $\tilde \dbP_{\tilde X}$-a.s. sense, and is determined by $f$ and $\tilde \dbP_{\tilde X}$.
  
{\it Step 2.} We now consider the  general distribution of $\tilde X$. For each $n\ge 1$, since $(\wh \O, d_{SK})$ is separable, there exists a partition $\{O^n_i, i\ge 1\}\subset \wh\O$ such that $d_{SK}(\wh\o, \wh\o^n_i) \le 2^{-n}$ for all $\wh\o \in O^n_i$, where $\wh\o^n_i\in O^n_i$ is fixed. Denote $\tilde X^n := \sum_{i\ge 1} \wh\o^n_i \1_{O^n_i}(\tilde X)$.  We remark that $\tilde X^n$ may not be $\dbF^{\tilde X}$-adapted, but such adaptedness is not needed here.  Since $\tilde X^n$ is discrete, by {\it Step 1}  we have $DF(t, \tilde X^n) = \psi_n(\tilde X^n)= \tilde \psi_n(\tilde X)$, where $\psi_n$ is defined by {\it Step 1} corresponding to $\tilde X^n$, and 
$\tilde \psi_n(\wh\o) := \sum_{i\ge 1}  \psi_n(\wh\o^n_i) \1_{O^n_i}(\wh\o)$, $\wh\o \in \wh\O$.
Clearly $\dbE^{\tilde \dbP}[ d^2_{SK}(\tilde X^n, \tilde X)] \le 2^{-2n}$, then by \reff{DFcont} we have 
\bea
\label{DFcont2}
\lim_{n\to\infty} \dbE^{\tilde \dbP}[|\tilde \psi_n( \tilde X) - DF(t, \tilde X)|^2] = 0.
\eea
  Thus there exists a subsequence $\{n_k\}_{k\ge 1}$ such that $\tilde \psi_{n_k} ( \tilde X) \to D F(t, \tilde X)$, $\tilde\dbP$-a.s.  Define
\bea
\label{Kpsin}
  \psi(\wh\o) := \limsup_{k\to \infty} \tilde \psi_{n_k}(\wh\o),\q K := \big\{\wh\o\in \wh\O: \limsup_{k\to\infty} \tilde \psi_{n_k}(\wh\o) = \liminf_{k\to\infty} \tilde \psi_{n_k}(\wh\o)\big\}.
\eea 
Then $\tilde\dbP(\tilde X\in K)=1$ and  $D F(t, \tilde X) = \psi(\tilde X)$, $\tilde \dbP$-a.s.

Moreover, let $\tilde X' \in \dbL^2(\tilde \O; \wh \O)$ be another process such that $\tilde \dbP_{\tilde X'}=\tilde \dbP_{\tilde X}$, and define $\tilde X^{'n}$ similarly by using  the same $\{O^n_i, \wh\o^n_i, i\ge 1\}$. Then $D F(t, \tilde X^{'n}) = \tilde \psi_n(\tilde X')$ for the same function $\tilde \psi_n$. Note that $\tilde \dbP(\tilde X'\in K) = \tilde \dbP(\tilde X\in K) =1$, then $\lim_{k\to\infty} \tilde \psi_{n_k}(\tilde X') = \psi(\tilde X')$, $\tilde \dbP$-a.s. On the other hand, $D F(t, \tilde X^{'n_k}) \to DF (t, \tilde X')$ in $\dbL^2$. So $D F(t, \tilde X') = \psi (\tilde X')$, $\tilde \dbP$-a.s., and thus $\psi$ does not depend on the choice of $\tilde X$.
\qed

Given the above theorem, particularly the fact that $\psi$ is determined by $\tilde \dbP_{\tilde X}$, we may introduce a function $\pa_\mu f: \wh\Th \times \wh\O\to \dbR^d$ such that $\pa_\mu f(t, \tilde \dbP_{\tilde X}, \wh\o) = \psi(\wh\o)$.  In particular, this implies: for any  $\wh\cF_t$-measurable $\mu$-square integrable  random variable $\xi: \wh \O\to \dbR^d$,
\bea
\label{pamuf}
\dbE^{\wh\mu}\big[ \pa_\mu f(t, \wh\mu, \wh X)\cd  \xi\big] = \lim_{\e\to 0} {f(t, \wh\mu \circ (\wh X +  \e \xi \1_{[t, T]})^{-1}) - f(t, \wh\mu)  \over \e}.
\eea

\begin{cor}
\label{cor-pamu}  Let  all the conditions in Theorem \ref{thm-pamu} hold true. Assume further that the continuity of $D F$ in \reff{DFcont} is uniform. Then  there exists a jointly Borel-measurable function $\pa_\mu f: \wh\Th\times \wh \O \to \dbR^d$  such that 
  \bea
  \label{DF}
  D F(t, \tilde X) = \pa_\mu f(t, \tilde \dbP_{\tilde X_{t\wedge \cd}}, \tilde X_{t\wedge \cd}),\q \tilde \dbP\mbox{-a.s.}
  \eea

 {\color{black} Moreover, if $\pa_\mu f(t, \cd)$ is jointly continuous in $\wh\cP_2 \times \wh\O$ for all $t$, then $\pa_\mu f$ is unique. }
\end{cor}
\proof In Theorem \ref{thm-pamu} {\it Step 1}, noting that $f$ is Borel measurable, then by \reff{DFx} one can easily see that $\pa_\mu f(t, \sum_{j\ge 1} p_j \d_{\{\wh\o_i\}}, \wh\o_i) := \psi(\wh\o_i)$ is jointly measurable.  Now consider the notations in Theorem \ref{thm-pamu} {\it Step 2}, and denote $\tilde \psi_n(t, \wh \mu, \wh\o) :=  \tilde\psi_n(\wh\o)$ which is jointly measurable in $(t, \wh\mu, \wh\o)$. By the uniform continuity of $D F$, one can choose a common subsequence $\{n_k, k\ge 1\}$ such that $\tilde \psi_{n_k} ( \tilde X) \to D F(t, \tilde X)$, $\tilde\dbP$-a.s. for all $\tilde X$. Denote 
$
\pa_\mu f(t, \wh \mu, \wh\o) := \limsup_{k\to\infty} \tilde \psi_{n_k}(t, \wh \mu, \wh\o).
$
Then $\pa_\mu f$ is jointly measurable and \reff{DF} holds true.

{\color{black}We now assume  $\pa_\mu f(t, \cd)$ is jointly continuous in $\wh\cP_2 \times \wh\O$ for all $t$. Notice again that $\pa_\mu f(t, \wh \mu,\cd)$ is unique, $\wh\mu$-a.s. Then, when $\supp(\wh\mu) = \wh\O$, by the continuity of $\pa_\mu f(t, \wh \mu,\cd)$ we see that $ \pa_\mu f(t, \wh \mu,\cd)$ is pointwise unique.   Finally, for any $\wh\mu\in \wh\cP_2$, there exist $\wh \mu_n \in \wh\cP_2$ such that $\supp(\wh\mu_n) = \wh\O$ for each $n$ and $\lim_{n\to\infty}\wh\cW_2(\wh\mu_n, \wh\mu)=0$. Then $\pa_\mu f(t, \wh\mu_n, \cd)$ is unique and $\dis\lim_{n\to\infty} \pa_\mu f(t, \wh \mu_n, \wh \o) = \pa_\mu f(t, \wh\mu, \wh\o)$. This clearly implies the uniqueness of $\pa_\mu f(t, \wh\mu, \wh\o)$.
}\qed

Now given $\pa_\mu f: \wh\Th \times \wh\O\to \dbR^d$, assume $\pa_\mu f(t,  \cd)$ is continuous and thus is unique.  In the spirit of Dupire \cite{Dupire} we may define further the derivative  function $\pa_{\o } \pa_\mu f: \wh\Th \times \wh\O\to \dbR^{d\times d}$ determined by:
\bea
\label{paxpamuf}
\pa_{\o} \pa_\mu f(t, \wh\mu, \wh\o) ~ x := \lim_{\e\to 0} { \pa_\mu f(t, \wh\mu, \wh \o + \e x\1_{[t, T]} ) - \pa_\mu f(t, \wh\mu, \wh \o)\over \e}, \q\mbox{for all}~ x\in \dbR^d.
\eea

\begin{eg}
\label{eg-derivative}
Let $d=1$ and $f(t,\wh\mu):=\dbE^{\wh\mu}\Big[ \wh X_t\int_0^t\wh X_sds\Big] -  \dbE^{\wh\mu}[ \wh X_t^2]\dbE^{\wh\mu}\Big[\int_0^t\wh X_sds\Big]$.
  Then 
\beaa
&\dis  \pa_t f(t, \wh \mu) = \dbE^{\wh\mu}[ \wh X_t^2] - \dbE^{\wh\mu}[\wh X_t^2]\dbE^{\wh\mu}[\wh X_t],&\\
&\dis \pa_\mu f(t, \wh\mu, \wh\o) = \int_0^t\wh \o_sds - 2 \wh \o_t \dbE^{\wh\mu}\Big[\int_0^t\wh X_sds\Big],\q \pa_\o \pa_\mu f(t, \wh\mu, \wh\o) =  - 2  \dbE^{\wh\mu}\Big[\int_0^t\wh X_sds\Big].&
\eeaa
\end{eg}
\proof First, note that
\beaa
f(t+\d, \wh \mu_{[0, t]}) &=& \dbE^{\wh\mu}\Big[ \wh X_t \int_0^{t+\d}\wh X_{t\wedge s}ds\Big] -  \dbE^{\wh\mu}[ \wh X_t^2]\dbE^{\wh\mu}\Big[\int_0^{t+\d}\wh X_{t\wedge s}ds\Big] \\
& =&  f(t, \wh \mu) + \d \dbE^{\wh\mu}[ \wh X_t^2] - \d\dbE^{\wh\mu}[\wh X_t^2]\dbE^{\wh\mu}[\wh X_t].
\eeaa
Then by \reff{patf} one can easily see that $\pa_t f(t, \wh \mu) = \dbE^{\wh\mu}[ \wh X_t^2] - \dbE^{\wh\mu}[\wh X_t^2]\dbE^{\wh\mu}[\wh X_t]$.

Next, for any appropriate $\tilde \dbP$ and $\tilde X$ on $(\tilde \O, \tilde \cF)$, we have
\beaa
F(t, \tilde X) = \dbE^{\tilde \dbP}\Big[ \tilde X_t\int_0^t\tilde X_sds\Big] -   \dbE^{\tilde \dbP}[ \tilde X_t^2]  \dbE^{\tilde \dbP}\Big[\int_0^t\tilde X_sds\Big] 
\eeaa
Then, 
\beaa
F(t, \tilde X + \xi \1_{[t, T]}) &=&  \dbE^{\tilde \dbP}\Big[ [\tilde X_t+\xi]\int_0^t\tilde X_sds\Big] -   \dbE^{\tilde \dbP}[ [\tilde X_t+\xi]^2]  \dbE^{\tilde \dbP}\Big[\int_0^t\tilde X_sds\Big] \\
&=& F(t, \tilde X) + \dbE^{\tilde \dbP}\Big[  \xi \int_0^t\tilde X_sds\Big] -   \dbE^{\tilde \dbP}[ 2\xi\tilde X_t+\xi^2]  \dbE^{\tilde \dbP}\Big[\int_0^t\tilde X_sds\Big] .
\eeaa
This implies
\beaa
D F(t, \tilde X) = \int_0^t\tilde X_sds - 2 \tilde X_t ~\!\dbE^{\tilde \dbP}\Big[\int_0^t\tilde X_sds\Big],
\eeaa
and thus
\beaa
\pa_\mu f(t, \wh \mu, \wh\o) =  \int_0^t\wh \o_sds - 2 \wh \o_t ~\! \dbE^{\wh\mu}\Big[\int_0^t\wh X_sds\Big].
\eeaa

Finally, by \reff{paxpamuf} it is straightforward to derive: $\pa_\o \pa_\mu f(t, \wh \mu, \wh\o) =  - 2 \dbE^{\wh\mu}\big[\int_0^t\wh X_sds\big]$.
\qed

\begin{defn}
\label{defn-C12hat}
Let $C^{1,1,1}(\wh \Th)$ be the set of  continuous mappings $f: \wh \Th \to \dbR$ such that there  exist  continuous functions $\pa_t f: \wh \Th\to \dbR$, $\pa_\mu f: \wh \Th \times \wh\O \to \dbR^d$, and $\pa_\o \pa_\mu f: \wh \Th \times \wh\O \to \dbR^{d\times d}$. 

Moreover, let $C^{1,1,1}_b(\wh \Th)\subset C^{1,1,1}(\wh \Th)$ denote the subset such that   $\pa_t f$ is bounded,  and $\pa_\mu f,  \pa_\o \pa_\mu f$ have linear growth in $\wh\o$: 
\bea
\label{poly}
|\pa_\mu f(t,\wh\mu, \wh\o)|+|\pa_\o \pa_\mu f(t, \wh\mu,\wh \o)| \le C\big[1+\|\wh\o\|\big],\q \mbox{for all}~(t,\wh\mu, \wh\o) \in \wh \Th \times \wh\O.
\eea
\end{defn}

\begin{rem}
\label{rem-C12hat}
{\rm
Our master equation \reff{master} below will involve the derivatives $\pa_t f, \pa_\mu f, \pa_\o\pa_\mu f$, but does not involve $\pa_\mu \pa_\mu f$ which can be defined in a natural way. The existence of $\pa_\o\pa_\mu f$ is of course a stronger requirement than that of $\pa_\mu f$, but roughly speaking it is  weaker than the existence of $\pa_\mu\pa_\mu f$. In the literature people call master equations involving $\pa_\mu \pa_\mu f$ second order, so our master equation is somewhat between first order and second order. 
\qed}
\end{rem}

\subsection{The functional It\^{o} formula}
\label{sect-FIto}
For any $L>0$, denote by $\wh \cP_L$ be the subset of $\mu\in \wh \cP_2$ such that  $\mu$ is a semimartingale measure with both the drift and diffusion characteristics bounded by $L$. {\color{black}To be precise, $\mu = \tilde \dbP \circ \tilde X^{-1}$, where $(\tilde \O, \tilde \dbF, \tilde \dbP)$ is a filtered probability space, $\tilde X_t = \tilde X_0 + \int_0^t \tilde b_s ds + \int_0^t \tilde \si_s d\tilde B_s$, $\tilde X_0 \in \dbL^2(\tilde \cF_0, \tilde \dbP; \dbR^d)$, $\tilde b: [0, T]\times \tilde \O \to \dbR^d$ and $\tilde \si: [0, T]\times \tilde \O\to \dbR^{d\times d}$ are $\tilde\dbF$-progressively measurable with $|\tilde b|, {1\over 2} |\tilde \si|^2\le L$, and $\tilde B$ is a $d$-dimensional $(\tilde \dbF, \tilde \dbP)$-Brownian motion.}  Note that, in particular, $\wh X$ is continuous in $t$, $\mu$-a.s., namely supp$(\mu) \subset \O \subset \wh\O$.  So $\mu$ can actually be viewed as a measure on $\O$ and thus we use the notation $\mu$ instead of $\wh \mu$ here.

\begin{thm}
\label{thm-Ito}
Let $f\in C^{1,1,1}_b(\wh \Th)$ and $\mu \in \wh \cP_L$ for some $L>0$.  Then
\bea
\label{functionalIto}
&&f(t,\mu)=f(0,\mu)+\int_0^t\pa_t f(s,\mu) ds\\
&&\qq +\dbE^\mu\Big[\int_0^t \pa_\mu f(s,\mu, \wh X)\cd d \wh X_s +\frac{1}{2} \int_0^t \pa_\o\pa_\mu f(s,\mu,\wh X): d\la \wh X\ra_s \Big].\nonumber
\eea
\end{thm}
\proof For notational simplicity, assume $d=1$ and $t=T$.  The general case can be proved without any additional difficulty.   Fix $\mu \in \wh\cP_L$ and let $(\tilde \O, \tilde \dbP, \tilde X)$ be the desired setting so that $\mu = \tilde \dbP \circ \tilde X^{-1}$.  Fix $n\ge 1$ and let $\pi: 0=t_0<t_1<\cds<t_n=T$ be a  uniform partition of $[0,T]$. Recall \reff{Xst} and denote 
\beaa
&\tilde X^n:=\sum_{i=0}^{n-1}\tilde X_{t_i}1_{[t_i,t_{i+1})}+\tilde X_T1_{\{T\}},\q \mu^n:= \tilde \dbP \circ (\tilde X^n)^{-1};&\\
&\tilde X^{n,\th} := \tilde X^n_{t_i\wedge \cd} + \th \tilde X_{t_i, t_{i+1}}\1_{[t_{i+1}, T]}, \q \mu^{n,\th} := \tilde \dbP \circ (\tilde X^{n,\th})^{-1},\q \th \in [0,1].&
\eeaa 
 Note that $\tilde X^n_{ t_{i+1}\wedge \cd} = \tilde X^n_{ t_i\wedge \cd} + \tilde X_{t_i, t_{i+1}} \1_{[t_{i+1}, T]}$.  Then
\bea
\label{I1234}
&& f(T, \mu^n) - f(0, \mu^n)= \sum_{i=0}^{n-1}\Big[f(t_{i+1}, \mu^n_{[0, t_{i+1}]})-f(t_i,\mu^n_{[0, t_{i}]})\Big]\nonumber\\
&&=\sum_{i=0}^{n-1}\Big[[f(t_{i+1}, \mu^n_{[0, t_i]})-f(t_i,\mu^n_{[0, t_i]})] + [f(t_{i+1}, \mu^n_{[0, t_{i+1}]})-f(t_{i+1},\mu^n_{[0, t_{i}]})]\Big]\nonumber\\
&&=\sum_{i=0}^{n-1}\Big[\int_{t_i}^{t_{i+1}}\pa_t f(t,\mu^n_{[0, t_i]})dt +  \int_0^1 \dbE^{\tilde \dbP} \big[ \pa_\mu f\big(t_{i+1}, \mu^{n,\th}, \tilde X^{n,\th})  \tilde X_{t_i, t_{i+1}} \big] d\th\Big] \nonumber\\
&&=\sum_{i=0}^{n-1}\int_{t_i}^{t_{i+1}}\pa_t f(t,\mu^n_{[0, t_i]})dt + \sum_{i=0}^{n-1}\int_0^1 \dbE^{\tilde \dbP} \big[ \pa_\mu f\big(t_{i+1}, \mu^{n,\th}, \tilde X^n_{t_i\wedge \cd})  \tilde X_{t_i,t_{i+1}} \big] d\th\nonumber \\
&&\qq+ \sum_{i=0}^{n-1} \int_0^1 \int_0^1 \dbE^{\tilde \dbP} \big[ \pa_\o \pa_\mu f\big(t_{i+1}, \mu^{n,\th}, \tilde X^{n, \tilde \th \th}) \th |\tilde X_{t_i,t_{i+1}}|^2 \big] d\tilde \th d\th \nonumber\\
&&=: I^n_1 + I^n_2 + I^n_3,
\eea
where $I^n_i$, $i=1,2,3$, are defined in an obvious way.

We now send $n\to \infty$. Since $\tilde X$ is continuous, $\tilde \dbP$-a.s.,  then, for any $t\in [0, T]$ and $\th \in [0, 1]$, 
\bea
\label{dSKconv}
d_{SK}(\tilde X^n, \tilde X) +  d_{SK}(\tilde X^n_{t_i\wedge \cd}, \tilde X_{t\wedge \cd}) + d_{SK}(\tilde X^{n,\th}_{t_{i+1}\wedge \cd}, \tilde X_{t\wedge \cd}) \to 0, \q \tilde \dbP\mbox{-a.s.}
\eea
where we always choose $i$ such that $t_i \le t< t_{i+1}$. Since $\|\tilde X^n\|\le \|\tilde X\|, \|\tilde X^{n,\th}\|\le \|\tilde X\|$, by the dominated convergence theorem we have
\beaa
\cW_2(\mu^n_{[0, t_i]}, \mu_{[0, t]}) + \cW_2(\mu^{n,\th}_{[0, t_{i+1}]}, \mu_{[0, t]}) \to 0.
\eeaa
Then, by the desired regularity of $f$, together with the boundedness of $\pa_t f$, \reff{poly}, and the fact that  the $\tilde b$ and $\tilde \si$ associated with $\tilde X$ are bounded, we can easily have
\beaa
&\dis\lim_{n\to\infty} \Big[f(T, \mu^n) - f(0, \mu^n)\Big] = f(T, \mu) - f(0, \mu),\q \lim_{n\to\infty} I^n_1 = \int_0^T \pa_t f(t, \mu) dt;&\\
&\dis\lim_{n\to\infty} \dbE^{\tilde \dbP} \Big[ \Big|\int_0^1\pa_\mu f\big(t_{i+1}, \mu^{n,\th}, \tilde X^n_{t_i\wedge \cd}) d\th - \pa_\mu f(t, \mu, \tilde X)\Big|^2\Big]=0; &\\
&\dis\lim_{n\to\infty} \dbE^{\tilde\dbP} \Big[ \Big|\int_0^1 \int_0^1 \pa_\o \pa_\mu f\big(t_{i+1}, \mu^{n,\th}, \tilde X^{n, \tilde \th \th}) \th d\tilde \th d\th -{1\over 2}\pa_\o \pa_\mu f(t, \mu, \tilde X)\Big|^2\Big]=0.&
\eeaa
Plug all these into \reff{I1234}, and recall that $\tilde \dbP\circ \tilde X^{-1} = \mu$, we can easily obtain \reff{functionalIto}.
\qed

We remark that it is possible to relax the technical conditions required for the functional It\^{o} formula \reff{functionalIto}, in particular we can allow  $\wh \mu \in \wh \cP_2$ to be semimartingale measures with supp$(\wh\mu)$ not within $\O$.  We also remark that, since $\la \wh X\ra$ is symmetric, in the last term of \reff{functionalIto} we may replace $\pa_\o\pa_\mu f(s,\mu,\wh X_{\cd})$ with 
\bea
\label{symmetric}
\pa^{sym}_\o\pa_\mu f(s,\mu,\wh X) := {1\over 2}\Big[\pa_\o\pa_\mu f(s,\mu,\wh X) + [\pa_\o\pa_\mu f(s,\mu,\wh X)]^\top\Big].
\eea

\subsection{The restriction on the space of continuous paths}

\begin{defn}
\label{defn-C12}
(i) Let $C^{1,1,1}(\Th)$ denote the set of $f: \Th\to \dbR$ such that there exists $\wh f\in C^{1,1,1}(\wh\Th)$ satisfying $\wh f = f$ on $\Th$, and define, for all $(t,\mu, \o) \in \Th \times \O$,
\bea
\label{paf}
\left.\ba{c}
\dis \pa_t f(t, \mu) := \pa_t \wh f(t,\mu),~ \pa_\mu f(t, \mu, \o) := \pa_\mu \wh f(t,\mu, \o),\\
\dis \pa_\o\pa_\mu f(t,\mu, \o) := \pa_\o\pa_\mu \wh f(t,\mu, \o),\q  \pa^{sym}_\o\pa_\mu f(t,\mu, \o) := \pa^{sym}_\o\pa_\mu \wh f(t,\mu, \o).
\ea\right.
\eea
Moreover, we say $f\in C_b^{1,1,1}(\Th)$ if the extension $\wh f\in C^{1,1,1}_b(\wh\Th)$.

(ii) Let $\cP_L$ denote the subset of $\mu\in \cP_2$ such that  $\mu$ is a semimartingale measure with both the drift and diffusion characteristics bounded by $L$.
\end{defn}

The following result is a direct consequence of Theorem \ref{thm-Ito}.
\begin{thm}
\label{thm-Ito2}
Let $f\in C^{1,1,1}_b(\Th)$. 

(i) The derivatives $\pa_t f, \pa_\mu f, \pa_\o^{sym}\pa_\mu f$ do not depend on the choices of $\wh f$;

(ii) For any $L>0$ and  $\mu \in  \cP_L$, we have
\bea
\label{functionalIto2}
&&f(t,\mu)=f(0,\mu)+\int_0^t\pa_t f(s,\mu) ds\\
&&\qq+\dbE^\mu\Big[\int_0^t \pa_\mu f(s,\mu, X)\cd d  X_s +\frac{1}{2} \int_0^t \pa^{sym}_\o\pa_\mu f(s,\mu,X): d\la X\ra_s \Big].\nonumber
\eea
\end{thm}
\proof (ii) follows directly from Theorem \ref{thm-Ito} and \reff{paf}. To see (i), the uniqueness of $\pa_t f$ is obvious.  Now fix $(t,\mu)\in \Th$ and let $\wh f$ be an arbitrary extension. For any bounded $\cF_t$-measurable $\dbR^d$-valued random variable $b_t$,  let $\tilde \mu \in \cP_2$ be such that $\tilde \mu = \mu$ on $\cF_t$ and $X_s - X_t = b_t [s-t]$, $t\le s\le T$, $\tilde \mu$-a.s.  Following the same arguments as in Theorem \ref{thm-Ito}, for any $\d>0$ we have
\beaa
f(t+\d,\tilde \mu)- f(t,\mu)=\int_t^{t+\d} \pa_t f(s,\tilde \mu) ds +\dbE^\mu\Big[\int_t^{t+\d} \pa_\mu\wh f(s, \tilde \mu, X)\cd b_t d s  \Big].
\eeaa
Divide both sides by $\d$ and send $\d\to 0$, we obtain the uniqueness of $\dbE^\mu[\pa_\mu\wh f(t, \mu, X)\cd b_t]$. Since $b_t$ is arbitrary, we see that $\pa_\mu\wh f(t, \mu, X)$ is unique, $\mu$-a.s. Similarly, for any bounded $\cF_t$-measurable $\dbR^{d\times d}$-valued random variable $\si_t$,  let $\tilde \mu \in \cP_2$ be such that $\tilde \mu = \mu$ on $\cF_t$ and $X$ is a $\tilde \mu$-martingale on $[t, T]$ with diffusion coefficient $\si_t$. Then similarly  we can show that $\dbE^\mu[\pa^{sym}_\o\pa_\mu\wh f(t, \mu, X): \si_t\si^\top_t]$ is unique, which implies the $ \mu$-a.s. uniqueness of  $\pa^{sym}_\o\pa_\mu\wh f(t, \mu, X)$.
\qed

We remark that, under some stronger technical conditions,  as in  Cont \& Fournie \cite{CF} one can show that $\pa_\o\pa_\mu f$ also does not depend on the choices of $\wh f$. However, the analysis below will depend only on $\pa_\o^{sym}\pa_\mu f$, so we do not pursue such generality here.

{\color{black}
\begin{rem}
\label{rem-state}
{\rm Let $V\in C^{1,1,1}(\Th)$. If $V(t, \mu) = V(t, \mu_t)$ is state dependent, it is clear that $\pa_\mu V(t, \mu, \o) = \pa_\mu V(t, \mu, \o_t)$ also depends only on the current state $\o_t$. Then naturally we may consider $\pa_x \pa_\mu V$ instead of $\pa_\o \pa_\mu V$. Throughout the paper we shall take this convention in the state dependent case.
\qed}
\end{rem}
}

\section{\color{black} Parabolic master equations and some applications} 
\label{sect-classical}
\setcounter{equation}{0}
In this paper we are interested in the following so called master equation:
\bea
\label{master}
\sL V (t, \mu):= \pa_t V(t, \mu) + G\big(t, \mu,  V(t,\mu), \pa_\mu V(t,\mu,\cd), \pa_\o\pa_\mu V(t,\mu, \cd)\big) =0,~ (t,\mu) \in \Th.
\eea
where $G(t, \mu, y, Z, \G)\in \dbR$ is defined in the domain where $(t,\mu, y)\in \Th \times \dbR$, and $(Z, \G)\in  C^0(\O; \dbR^d) \times C^0(\O; \dbR^{d\times d})$  are $\cF_t$-measurable.  We remark that $G$ depends on the whole random variables $Z$ and $\G$, rather than their values. Such dependence is typically through $\dbE^\mu$ in the form: $G = G_1(t, \mu, y, \dbE^\mu[G_2(t, \mu, y, Z, \G)])$ for some deterministic functions  $G_1: \Th \times \dbR\times \dbR^k \to \dbR$ and $G_2: \Th \times \dbR \times \dbR^d \times \dbR^{d\times d}\to \dbR^k$  for some dimension $k$.

\begin{assum}
\label{assum-G}
(i) $G$ is continuous in $(t, \mu)$ and uniformly Lipschitz continuous in $y$ with a Lipschitz constant $L_0$.

(ii) $G$ is uniformly Lipschitz continuous in $(Z, \G)$ with a Lipschitz constant $L_0$ in the following sense: for any $(t,\mu, y)$ and any $\cF_t$-measurable random variables $Z_1, \G_1, Z_2, \G_2$, there exist $\cF_t$-measurable random variables $b_t, \si_t$ such that $|b_t|, {1\over 2}|\si_t|^2\le L_0$, and   
\bea
\label{GLip}
G(t,\mu, y, Z_1, \G_1) - G(t,\mu, y, Z_2, \G_2)= \dbE^\mu\Big[ b_t \cd [Z_1-Z_2]+ {1\over 2} \si_t \si_t^\top : [\G_1-\G_2]\Big].
\eea
\end{assum}
We remark that, while \reff{GLip} may look a little less natural,  one can  easily verify it for all the examples in this paper. Moreover, when $\mu$ is degenerate and thus $Z, \G$ becomes deterministic numbers rather than random variables, \reff{GLip} is equivalent to the standard Lipschitz continuity.  

\begin{rem}
\label{rem-parabolic}
{\rm By \reff{GLip}, it is clear that $G$ depends on $\G$ only through $\G^{sym} := {1\over 2}[\G + \G^\top]$, and $G$ is increasing in $\G^{sym}$. So \reff{master} depends on $\pa_\o \pa_\mu V$ only through $\pa^{sym}_\o\pa_\mu V$, which is unique (or say, well defined) by  Theorem \ref{thm-Ito2} (i). 
\qed}
\end{rem}

\begin{defn}
\label{defn-classical}
Let $V \in C^{1,1,1}(\Th)$. We say $V$ is a classical solution (resp. classical subsolution, classical supersolution) of the master equation \reff{master} if 
\beaa
\sL V (t, \mu) = ~(\mbox{resp.}~ \ge, \le)~ 0,\q\mbox{for all}~ (t,\mu) \in \Th.
\eeaa
\end{defn}

{\color{black}In the rest of this section we show several examples, which can be viewed as some typical applications  of our parabolic master equations. We remark that the smooth differentiability of the involved value functions are often very challenging (and in general may not be true), and thus the main focus of this paper is the viscosity solution. However, for illustration purpose, in this section we shall assume the value functions are smooth and verify they are classical solutions of the corresponding master equations. We shall also show in some special cases that  the value functions under consideration  are indeed smooth. 

\subsection{Stochastic optimization with deterministic controls}
\label{sect-deterministic}
While the value function of a control problem will automatically be path dependent if the coefficients are path dependent, in this subsection we present a state dependent example which endogenously induces a path dependent master equation.
Consider a standard control problem :
\bea
\label{deterministicV0}
\left.\ba{c}
\dis V_0 = \sup_{\a\in \cA} Y^\a_0,\q\mbox{where}\\
\dis X^\a_t = x_0 + \int_0^t b(s, X^\a_s, \a_s) ds + \int_0^t \si(s, X^\a_s, \a_s) dB_s,\\
\dis Y^\a_t = g(X^\a_T) + \int_t^T f(s, X^\a_s, Y^\a_s, Z^\a_s, \a_s) ds - \int_t^T Z^\a_s dB_s,\\ 
\ea\right.
\eea
Here $B$ is a $\dbP_0$-Brownian motion,  the control $\a$ takes values in an appropriate set $A$, and the coefficients $b, \si, f, g$ satisfy standard technical conditions which we shall not specify.  When $\cA$ is the set of $\dbF^B$ or $\dbF^{X^\a}$-progressively measurable processes, it is a classical result that $V_0 = u(0, x_0)$, where $u$ is the solution to an HJB equation, and the optimal control $\a^*$, if it exists, typically is feedback type: $\a^*_t = I(t, X^*_t)$ for some deterministic function $I$. 
 
 In practice, quite often one needs some time to analyze the information (including the time for numerical computation), and in operations management, one needs to place orders some time before the parts are actually used. Mathematically, this amounts to require $\a_t$ to be $\cF_{t-\d}$-measurable, for some information delay parameter $\d$. For simplicity let's assume $T\le \d$, then $\a$ becomes deterministic.  In the rest of this subsection, we shall consider the problem \reff{deterministicV0} where 
 \bea
 \label{deterministic}
 \mbox{the admissible controls $\a\in \cA$ are deterministic.}
 \eea
This seemingly simple problem is actually more involved, and to our best knowledge is not covered by the existing methods in the literature. The main difficulty is the time inconsistency. Indeed, if one natively defines $u(t,x)$ as the value of the optimization problem on $[t, T]$ with initial condition $X_t=x$,  then $u$ does not satisfy  the dynamic programming principle  and consequently it does not satisfy any PDE.
  }
  
 In Saporito \& Zhang \cite{SZ} we investigated this problem in the case $f = f(t, X^\a_t, \a_t)$. It turns out that in this case the optimal $\a^*$ takes the form: $\a^*_t = I(t, \cL_{X^*_t})$, which is deterministic.  To be precise, for any $(t, \mu)\in \Th$ and $\a\in \cA$ (deterministic),  let $\dbP^{t,\mu, \a}$ be the unique solution satisfying $\dbP^{t,\mu,\a}_{[0, t]} = \mu_{[0, t]}$ and,  for some  $\dbP^{t,\mu, \a}$-Brownian motion $B^\a$,
  \bea
  \label{Ptmua}
\dis X_s = X_t + \int_t^s b(r, X_r, \a_r) dr + \int_t^s \si(r, X_r, \a_r) dB^\a_r,\q s\in [t, T], \dbP^{t,\mu, \a}\mbox{-a.s.}
\eea
Define 
\bea
\label{deterministicV}
V(t,\mu) := \sup_{\a\in\cA} \dbE^{\dbP^{t,\mu, \a}}\Big[g(X_T) + \int_t^T f(s, X_s, \a_s)ds\Big].
\eea 
Then by \cite{SZ} we have the following result.

\begin{prop}
\label{prop-deterministic1} 
Assume $f = f(t,x,a)$, $b, \si, f, g$ satisfy standard technical conditions, and define $V$ by \reff{deterministicV} under \reff{deterministic}. Then

(i)  $V(t, \mu) = V(t, \mu_t)$ is state dependent and the dynamic programming principle holds: 
\bea
\label{deterministic-DPP1}
V(t_1, \mu_{t_1}) = \sup_{\a\in\cA} \Big[ V(t_2, \dbP^{t_1, \mu, \a}_{t_2}) + \int_{t_1}^{t_2} \dbE^{\dbP^{t_1,\mu, \a}} \big[ f(s, X_s, \a_s)\big] ds\Big],\q t_1 < t_2.
\eea

(ii) Assume $V\in C^{1,1,1}(\Th)$ (more precisely $C^{1,1,1}([0, T]\times \cP_2(\dbR^d))$ here, and also recalling Remark \ref{rem-state}), then $V$ is the classical solution to the following master equation:
\bea
\label{deterministic-master1}
\left.\ba{lll}
 \dis\pa_t V(t,\mu) + \sup_{a\in A} \dbE^\mu\Big[{1\over 2}  \pa_x \pa_\mu V(t, \mu, X_t) : \si\si^\top(t, X_t, a)+   \pa_\mu V(t, \mu, X_t) \cd b(t, X_t, a) \\
 \dis\qq\qq + f(t, X_t, a)\Big] =0,\qq\qq\qq  V(T, \mu) = \dbE^\mu[g(X_T)].
 \ea\right.
 \eea 
 
 (iii)  Assume further that the Hamiltonian in \reff{deterministic-master1} has an optimal argument $a^* = I(t, \mu_t)$, and the following McKean-Vlasov SDE has a solution:
 \bea
 \label{deterministic-X^*}
 X^*_t = x_0 + \int_0^t b(s, X^*_s, I(s, \cL_{X^*_s}) )ds + \int_0^t \si(s, X^*_s, I(s, \cL_{X^*_s}) )dB_s,\q \dbP_0\mbox{-a.s.}
 \eea
 Then $\a^*_t := I(t, \cL_{X^*_t})$ is an optimal control to the problem \reff{deterministicV0}.
 \end{prop} 
 We remark that the expectation involved in \reff{deterministic-master1}  is a function of $(t, \mu_t, a)$, so the optimal control $a^*$ takes the form $I(t, \mu_t)$ in (iii). 

{\color{black}
While induced endogenously, the master equation \reff{deterministic-master1} is still state dependent. 
We now consider \reff{deterministicV0} with nonlinear $f$, again with deterministic $\a$. The general case is quite involved, and  we consider only a special case here:
 $f = f(t, X_t, Y_t)$.  Given $(t,\mu)\in \Th$ and $\a\in \cA$, let $\dbP^{t,\mu,\a}$ be defined by \reff{Ptmua}, and consider the following BSDE:
 \bea
  \label{Ytmua}
\dis Y^{t,\mu,\a}_s = g(X_T) + \int_s^T f(r, X_r, Y^{t,\mu,\a}_r) dr  + M^{t,\mu,\a}_T - M^{t,\mu,\a}_s,\q s\in [0, T], ~\dbP^{t,\mu, \a}\mbox{-a.s.}
\eea
Here the component $M$ of the solution pair $(Y, M)$  is a $\dbP^{t,\mu,\a}$-martingale. 
If we set $V(t, \mu):= \sup_{\a\in \cA} \dbE^\mu[Y^{t,\mu,\a}_t]$ as in \reff{deterministicV}, then $V$ will still be state dependent, but in general the DPP in the spirit of \reff{deterministic-DPP1} does not hold, because of the nonlinearity of $f$. To keep the time consistency,  in this case we shall define the value function as:
\bea
\label{deterministic-V2}
V(t,\mu) := \sup_{\a\in \cA}  \dbE^\mu[Y^{t,\mu,\a}_0].
\eea
Note that $V(t,\mu)$ is path dependent, in particular, $V(T,\mu) = Y^{T,\mu}_0$, where $Y^{T,\mu}$ is the solution to BSDE \reff{Ytmua} under $\mu$. Then we  can extend Proposition \ref{prop-deterministic1} to this case.
\begin{thm}
\label{thm-deterministic2} 
Assume $f= f(t,x,y)$, $b, \si, f, g$ satisfy standard technical conditions, and define $V$ by \reff{deterministic-V2} under \reff{deterministic}. Then

(i)  The following dynamic programming principle holds: 
\bea
\label{deterministic-DPP2}
V(t_1, \mu) = \sup_{\a\in\cA}  V(t_2, \dbP^{t_1, \mu, \a}),\q t_1 < t_2.
\eea

(ii) Assume $V\in C^{1,1,1}(\Th)$, then $V$ satisfies the path dependent  master equation:
\bea
\label{deterministic-master2}
&\dis \pa_t V(t,\mu) +  \sup_{a\in A} \dbE^\mu\Big[{1\over 2} \pa_\o \pa_\mu V(t, \mu, X),  a) :  \si\si^\top(t, X, a)  +  \pa_\mu V(t, \mu, X) \cd b(t, X, a)\Big]  =0,\nonumber\\
& V(T, \mu) =  Y^{T,\mu}_0.
 \eea 
 
 (iii)  Assume further that the Hamiltonian in \reff{deterministic-master2} has an optimal argument $a^* = I(t, \mu_{[0,t]})$, and the following McKean-Vlasov SDE has a solution:
 \bea
 \label{deterministic-X^*2}
 X^*_t = x_0 + \int_0^t b(s, X^*_s, I(s, \cL_{X^*_{s\wedge \cd}}) )ds + \int_0^t \si(s, X^*_s, I(s, \cL_{X^*_{s\wedge \cd}}) )dB_s,\q \dbP_0\mbox{-a.s.}
 \eea
 Then $\a^*_t := I(t, \cL_{X^*_{t\wedge \cd}})$ is an optimal control to the problem \reff{deterministicV0}.
 \end{thm} 
\proof (i) We emphasize that, since the $\dbP^{t_1, \mu, \a}$ inside $V(t_2, \cd)$  is deterministic, the DPP \reff{deterministic-DPP2} does not require any regularity or even measurability of $V$. Indeed, denote the right side of \reff{deterministic-DPP2} as $\tilde V(t_1, \mu)$. For any $\a\in \cA$,  by the flow property of SDEs and BSDEs we have
\beaa
\dbP^{t_1,\mu,\a} = \dbP^{t_2, \dbP^{t_1, \mu, \a}, \a},\q\mbox{and thus}\q Y^{t_1, \mu,\a}_0 = Y^{t_2, \dbP^{t_1,\mu,\a}, \a}_0.
\eeaa
Note that $\dbP^{t_1,\mu,\a} = \mu$ on $\cF_0$. This implies that
\beaa
\dbE^\mu[Y^{t_1, \mu,\a}_0] = \dbE^{\dbP^{t_1,\mu,\a}}[Y^{t_2, \dbP^{t_1,\mu,\a}, \a}_0] \le V(t_2, \dbP^{t_1, \mu, \a}).
\eeaa
Then by \reff{deterministic-V2} we see that $V(t_1, \mu) \le \tilde V(t_1,\mu)$. To see the opposite inequality, for any $\a\in \cA$ and any $\e>0$, there exists $\a^\e\in \cA$ such that 
\beaa
  V(t_2, \dbP^{t_1, \mu, \a}) \le \dbE^\mu[Y^{t_2, \dbP^{t_1, \mu, \a}, \a^\e}_0] +\e.
\eeaa
Denote $\tilde \a^\e_s := \a_s \1_{[0, t_2)}(s) + \a^\e_s \1_{[t_2, T]}(s)$. Then clearly $\tilde \a^\e \in \cA$, $\dbP^{t_1, \mu, \tilde \a^\e}_{[0, t_2]} = \dbP^{t_1, \mu,  \a}_{[0, t_2]}$, and
\beaa
\dbE^\mu[ Y^{t_2, \dbP^{t_1, \mu, \a}, \a^\e}_0] =  \dbE^\mu[Y^{t_2, \dbP^{t_1, \mu, \tilde \a^\e}, \a^\e}_0]  =\dbE^\mu[  Y^{t_1, \mu, \tilde \a^\e}_0]  \le V(t_1,\mu).
 \eeaa
 This implies that $V(t_2, \dbP^{t_1, \mu, \a}) \le V(t_1, \mu) + \e$. Then it follows from the arbitrariness of $\a$ and $\e$ that $\tilde V(t_1, \mu) \le  V(t_1,\mu)$.
 
 (ii) By applying the functional It\^{o} formula \reff{functionalIto} on the right side of \reff{deterministic-DPP2} we obtain the master equation \reff{deterministic-master2} immediately. The terminal condition follows from the definitions. 
 
(iii) Denote $\mu^* := \cL_{X^*}$ and $\a^*_t := I(t, \mu^*_{[0,t]})$. Apply the functional It\^{o} formula \reff{functionalIto} on $V(t,\mu^*)$ we obtain:
\beaa
&&{d\over dt} V(t, \mu^*) \\
&&= \pa_t V(t, \mu^*) + \dbE^{\mu^*}\Big[ {1\over 2} \pa_\o \pa_\mu V(t, \mu^*, X): \si\si^\top(t, X_t, \a^*_t) + \pa_\mu V(t, \mu, X) \cd b(t, X_t, \a^*_t)\Big]\\
&&= \pa_t V(t,\mu^*) + \sup_{a\in A} \dbE^{\mu^*}\Big[{1\over 2} \pa_\o \pa_\mu V(t, \mu^*, X) :\si\si^\top(t, X_t, a) + \pa_\mu V(t, \mu^*, X) \cd b(t, X_t, a)\Big],
\eeaa
where the last equality thanks to the fact that $\a^*$ is an optimal argument of the Hamiltonian. By the master equation \reff{deterministic-master2} we obtain ${d\over dt} V(t, \mu^*)=0$. Thus, noting that $\mu^* = \dbP^{0, \d_{\{x_0\}}, \a^*}$, 
\beaa
V_0=V(0, \d_{\{x_0\}}) = V(T, \mu^*) =   Y^{0, \d_{\{x_0\}}, \a^*}_0.
\eeaa
That is, $\a^*$ is an optimal control.
\qed

\subsection{Mean field control problems}
\label{sect-meanfieldcontrol}
The mean field control problem is one major application of the master equations, and will be studied in more details in Section \ref{sect-HJB} below. Consider a system of $N$  controlled interacting particle system: $i=1,\cds, N$,
\bea
\label{Xia}
\left.\ba{c}
\dis X^{\a, i}_t = x_i + \int_0^t b(s, X^{\a,i}_s, \mu^N_s, \a_s(X^{\a,i})) ds +   \int_0^t \si(s, X^{\a,i}_s, \mu^N_s, \a_s(X^{\a,i})) dB^i_s,\\
\mbox{where}\q \mu^N_s := {1\over N} \sum_{i=1}^N \d_{\{X^{\a, i}_s\}}.
\ea\right.
\eea
Here $B^i$ are independent Brownian motions,  the control $\a$ is a closed loop control and is chosen by a central planner (and thus the same $\a$ for all $i$), and the interaction is through the empirical measure $\mu^N$.  Assume $\mu^N_0 :=  {1\over N} \sum_{i=1}^N \d_{\{x_i\}} \to \mu_0$, while highly nontrivial, under appropriate conditions one can show that, see e.g. Lacker \cite{Lacker} (for relaxed controls), the above system converges to the following controlled McKean-Vlasov SDE with initial distribution $\cL_{X_0} = \mu_0$:
\bea
\label{Xa0}
\dis X^{\a}_t = X_0 + \int_0^t b(s, X^{\a}_s, \cL_{X^\a_s}, \a_s(X^{\a})) ds +   \int_0^t \si(s, X^{\a}_s, \cL_{X^\a_s}, \a_s(X^{\a})) dB_s,~ \dbP_0\mbox{-a.s.}
\eea
In many applications, the dynamics could be path dependent (e.g. SDEs with delays), so at below we extend \reff{Xa0} to the path dependent equation. Moreover, we shall consider a dynamic setting. To be precise, fix $t$ and a process $\xi$ on $[0, t]$, for a control $\a$, let $X^\a_s := \xi_s$ for $s\in [0, t]$ and consider the following equation on $[t, T]$ under $\dbP_0$:
 \bea
\label{Xa1}
\dis X^{\a}_s = \xi_t + \int_t^s b(r, X^{\a}_{r\wedge \cd}, \cL_{X^\a_{r\wedge \cd}}, \a_r(X^{\a}_{r\wedge \cd})) dr +    \int_t^s \si(r, X^{\a}_{r\wedge \cd}, \cL_{X^\a_{r\wedge \cd}}, \a_r(X^{\a}_{r\wedge \cd})) dB_r.
\eea
Since we will only care about the law of $X^\a$, it is more convenient to use the weak formulation in the canonical setting. That is, instead of fix $\dbP_0$ and consider the controlled process $X^\a$, we fix the canonical process $X$ and consider the controlled probability $\dbP^\a$. Now  given $(t,\mu)\in \Th$ and a control $\a$, let $\dbP^{t,\mu,\a}\in \cP_2$ be such that $\dbP^{t,\mu,\a}_{[0, t]} = \mu_{[0,t]}$ and, for $s\in [t, T]$ and  for some $\dbP^{t,\mu,\a}$-Brownian motion $B^\a$, the following holds $\dbP^{t,\mu,\a}$-a.s.
 \bea
\label{control-Pa}
\dis X_s = X_t + \int_t^s b(r, X_{r\wedge \cd},  \dbP^{t,\mu,\a}_{[0, r]}, \a_r(X_{r\wedge \cd})) dr +    \int_t^s \si(r, X_{r\wedge \cd},  \dbP^{t,\mu,\a}_{[0, r]}, \a_r(X_{r\wedge \cd})) dB^\a_r.
\eea
Note that $\a$ becomes a standard $\dbF$-progressively measurable process now. Our admissible controls are: for some appropriate set $A$ and for any $t_0$,
\bea
\label{control-cA}
\left.\ba{c}
\cA_{t_0}:= \Big\{\a: [t_0, T] \times \O \to A: \mbox{for any $t\in [t_0, T]$ and }\\
\mbox{any $\dbP_0$-square integrable process $\xi$, SDE \reff{Xa1} has a unique weak solution}.\Big\}
\ea\right.
\eea
Then  \reff{control-Pa} has a unique solution $\dbP^{t,\mu,\a}$ for any $\a\in \cA_t$.

We are now ready to define our value function:
\bea
\label{control-Vtmu}
V(t,\mu) := \sup_{\a\in \cA_t} J(t,\mu,\a) := \sup_{\a\in \cA_t} \dbE^{\dbP^{t,\mu,\a}}\Big[g(X, \dbP^{t,\mu,\a}) + \int_t^T f(s, X, \dbP^{t,\mu,\a}, \a_s)ds\Big].
\eea
Similar to Theorem \ref{thm-deterministic2}, we have the following result.
\begin{thm}
\label{thm-meancontrol} 
Assume $b, \si, f, g$ satisfy standard technical conditions, in particular  they are $\dbF$-adapted both in $X$ and in $\mu$, and define $V$ by \reff{control-Pa}-\reff{control-Vtmu}. Then

(i)  The following dynamic programming principle holds: 
\bea
\label{meancontrol-DPP}
V(t_1, \mu) = \sup_{\a\in\cA_{t_1}} \Big[V(t_2, \dbP^{t_1, \mu, \a}) + \int_{t_1}^{t_2} \dbE^{\dbP^{t_1,\mu,\a}}[f(s, X, \dbP^{t_1,\mu,\a}, \a_s)] ds\Big],\q t_1 < t_2.
\eea

(ii) Assume $V\in C^{1,1,1}(\Th)$, then $V$ satisfies the path dependent  master equation:
\bea
\label{meancontrol-master}
&\dis \pa_t V(t,\mu) + \dbE^\mu\big[ \sup_{a\in A} G_2(t, \mu, X,  \pa_\mu V(t, \mu, X), \pa_\o \pa_\mu V(t, \mu, X),  a)\big],\nonumber\\
&\dis  V(T, \mu) = \dbE^\mu\big[g(X,\mu)\big], \\
&\dis \mbox{where}\q G_2(t,\mu,\o, z, \g, a) := {1\over 2}\g :\si\si^\top(t, \o,\mu, a) + z \cd b(t, \o,\mu, a) + f(t, \o, \mu, a).\nonumber
 \eea 
 
 (iii)  Assume further that the Hamiltonian   in \reff{meancontrol-master} has an optimal argument $a^* = I(t, \o, \mu)$,  and the following McKean-Vlasov SDE has a solution:
 \bea
 \label{meancontrol-X^}
 X^*_t &=& x_0 + \int_0^t b(s, X^*,  \cL_{X^*}, I(s,  X^*, \cL_{X^*}) )ds \\
 &&+ \int_0^t \si(s, X^*,  \cL_{X^*}, I(s,  X^*, \cL_{X^*}) )dB_s,\q \dbP_0\mbox{-a.s.}\nonumber
 \eea
 Denote $\a^*_t := I(t, X^*, \cL_{X^*})$. If $\a^*\in\cA_0$, then it is an optimal control to the problem $V(0, \d_{\{x_0\}})$ in  \reff{control-Vtmu}.
 \end{thm} 
 
 We remark that, since the control $\a$ is deterministic in Theorem  \ref{thm-deterministic2}, in \reff{deterministic-master2} the $\sup_{a\in A}$ is outside of the expectation $\dbE^\mu$ and thus the optimal control depends only on $\mu$, but not on $X$. Here, in \reff{meancontrol-master} the  $\sup_{a\in A}$ is inside of the expectation $\dbE^\mu$ and thus the optimal control depends  on $X$ as well. 
 
 \ms
 \proof The proof of (ii) and (iii) are almost the same as that of  Theorem  \ref{thm-deterministic2}, we thus omit it. The proof of (i) is also similar, but since the involvement of $\cA_t$ is quite subtle, as we will discuss in more details in Section \ref{sect-HJB}, we provide a detailed proof again. We emphasize that, even though $\a$ is random here, the $\dbP^{t_1, \mu, \a}$ inside $V(t_2,\cd)$ is still deterministic and the DPP \reff{meancontrol-DPP}  does not require the measurability of $V$. 
 
 The proof relies on the following two compatibility properties of $\cA_t$: for any $t_1 < t_2$,
 \bea
 \label{compatibility}
 \left.\ba{c}
 \mbox{for any $\a\in \cA_{t_1}$, we have $\a_{[t_2, T]} \in \cA_{t_2}$};\\
 \mbox{for any $\a^1 \in \cA_{t_1}$, $\a^2\in \cA_{t_2}$, we have}~ \a := \a^1 \1_{[t_1, t_2)} + \a^2\1_{[t_2, T]}\in \cA_{t_1}.
 \ea\right.
 \eea
 Now denote the right side of \reff{meancontrol-DPP} as $\tilde V(t_1, \mu)$. On one hand, for any $\a\in \cA_{t_1}$, denote $\tilde \a:= \a_{[t_2, T]}$ and $\tilde \mu := \dbP^{t_1,\mu, \a}$. Note that $\tilde \a\in \cA_{t_2}$, thanks to the first line of \reff{compatibility}. Then $\dbP^{t_1,\mu,\a} =\dbP^{t_2,\tilde \mu,\tilde\a}$, and thus
 \beaa
 J(t_1, \mu,\a) &=& J(t_2,\tilde \mu,\tilde\a) +\dbE^{\tilde \mu}\Big[ \int_{t_1}^{t_2} f(s, X, \tilde \mu, \a_s)ds\Big]\\
 &\le& V(t_2, \tilde \mu) + \dbE^{\tilde \mu}\Big[ \int_{t_1}^{t_2} f(s, X, \tilde \mu, \a_s)ds\Big] \le \tilde V(t_1,\mu).
 \eeaa
 This implies that $V(t_1,\mu) \le \tilde V(t_1, \mu)$. One the other hand, for any $\a\in \cA_{t_1}$ and any $\e>0$, there exists $\tilde\a \in \cA_{t_2}$ such that: again denoting $\tilde \mu := \dbP^{t_1,\mu, \a}$,
 \beaa
 V(t_2, \tilde \mu) \le J(t_2,\tilde \mu,\tilde\a)  +\e.
 \eeaa
 Now denote $\hat \a := \a \1_{[t_1, t_2)} + \tilde\a\1_{[t_2, T]}\in \cA_{t_1}$, thanks to the second line of \reff{compatibility}. Then $\dbP^{t_1,\mu, \hat \a} =\dbP^{t_2,\tilde \mu, \tilde\a}$, and thus
 \beaa
 &&V(t_2, \tilde \mu) + \dbE^{\tilde \mu}\Big[ \int_{t_1}^{t_2} f(s, X, \tilde \mu, \a_s)ds\Big]\\
 &\le& J(t_2,\tilde \mu,\tilde\a) +\dbE^{\tilde \mu}\Big[ \int_{t_1}^{t_2} f(s, X, \tilde \mu, \a_s)ds\Big]+\e=  J(t_1,\mu,\hat\a)+\e \le V(t_1,\mu) + \e.
 \eeaa
 This implies $\tilde V(t_1, \mu) \le V(t_1, \mu)$.
 \qed
}

For illustration purpose, in the rest of this subsection we show that $V$ is indeed smooth when there is no control, and hence the master equation is linear. For simplicity we assume $d=1$, $b=0$, $\si=1$, and $f, g$ do not depend on $\mu$ and thus the path dependence is only through $X$.   For this purpose, let $(t,\mu)\in \Th$, denote by $\dbP^{t,\mu}_0\in \cP_2$ be such that $\dbP^{t,\mu}_0 = \mu$ on $\cF_t$ and $X_{t, \cd}$ is a $\dbP^{t,\mu}_0$-Brownian motion on $[t, T]$ independent of $\cF_t$. For $g: \wh \O \to \dbR$, define $D_t g: \wh\O \to \dbR$ by:
\bea
\label{Dg}
 D_t g(\wh \o)  :=  \lim_{\e\to 0}{g(\wh \o + \e  \1_{[t, T]}) - g(\wh\o)\over \e},
\eea
and define $D^2_t g: \wh \O \to \dbR$ similarly. We note that $D_t g$ is essentially the Malliavin derivative, and in particular $D_t g = 0$ if $g$ is $\cF_s$-measurable for some $s<t$. 
\begin{eg}
\label{eg-heat}
Let $g\in C^0_b (\wh\O; \dbR)$ and $f\in C^0_b([0, T]\times \wh\O; \dbR)$. Assume $D_t g, D^2_t g, D_t f, D^2_t f$ exist and are bounded, and  $D_t g(\wh\o), D^2_tg(\wh\o)$ are jointly continuous in $(t, \wh\o)$ under the distance $d((t,\wh\o), (t',\wh \o')) := |t-t'| + \|\wh\o-\wh\o'\|$,  $D_t f(s, \wh \o), D^2_t f(s, \wh\o)$ are jointly continuous in $(t, s, \wh\o)$ under the distance $d((t,s,\wh\o), (t', s', \wh \o')) := |t-t'| + |s-s'|+ \|\wh\o_{s\wedge \cd}-\wh\o'_{s'\wedge \cd}\|$. Define
\bea
\label{heat-V}
V(t, \mu) := \dbE^{\dbP^{t,\mu}_0}\Big[ g(X) + \int_t^T f(s, X) ds\Big].
\eea
Then $V \in C^{1,1,1}_b(\Th)$ and satisfies the following linear master equation:
\bea
\label{linearMaster}
\pa_t V(t, \mu) + \dbE^\mu\Big[ {1\over 2}  \pa_\o\pa_\mu V(t, \mu, X) + f(t, X)\Big]=0,\q V(T, \mu) = \dbE^\mu \big[g(X)\big].
\eea
\end{eg}
\proof The proof  follows similar arguments as in Peng \& Wang \cite{PengW}, which deals with semilinear path dependent PDEs, so we shall only sketch it. We remark that the continuity of $f$ implies its $\dbF$-adaptedness.

 First it is clear that we can extend \reff{heat-V} to all $(t, \wh\mu) \in \wh \Th$ in an obvious way. Denote $(\wh \o \otimes_t \o)_s := \wh\o_s \1_{[0, t]}(s) + [\wh\o_t + \o_s - \o_t]\1_{(t, T]}(s)$ for all $\wh\o\in \wh\O$ and $\o\in \O$. Then 
\beaa
V(t, \wh\mu) = \dbE^{\wh\mu}[ u(t, \wh X)],\q \mbox{where}\q u(t, \wh \o) := \dbE^{\dbP_0}\Big[g(\wh\o\otimes_t X) + \int_t^T f(s, \wh\o\otimes_t X) ds\Big].
\eeaa
By straightforward computation, we have
\beaa
\pa_\mu V(t, \wh\mu, \wh\o) = \pa_\o u(t, \wh\o) =  \dbE^{\dbP_0}\Big[D_t g(\wh\o\otimes_t X) + \int_t^TD_t f(s, \wh\o\otimes_t X) ds\Big],
\eeaa
where $\pa_\o u$ is Dupire's path derivative as in \reff{paxpamuf}. We note that in this particular case $\pa_\mu V$ actually does not depend on $\mu$. Then
\beaa
\pa_\o \pa_\mu V(t, \wh\mu, \wh\o) = \pa_\o \pa_\o u(t,\wh\o) =   \dbE^{\dbP_0}\Big[D^2_t g(\wh\o\otimes_t X) + \int_t^TD^2_t f(s, \o\otimes_t X) ds\Big].
\eeaa
By our conditions, it is quite obvious that $V, \pa_\mu V, \pa_\o\pa_\mu V$ are continuous. 

On the other hand, note that
\beaa
V(t+\d, \wh\mu_{[0, t]}) - V(t, \wh\mu_{[0, t]}) = \dbE^{\wh \mu} \Big[u(t+\d, \wh X_{t\wedge \cd}) - u(t, \wh X) \Big]
\eeaa
Fix $t$ and $t+\d$, let $t=t_0<\cds<t_n= t+\d$. Recall \reff{Xst} and denote, for $0\le m\le n$, 
\beaa
X^{n,m}:= \wh\o_{t\wedge \cd} + \sum_{i=1}^m X_{t_{i-1}, t_i} \1_{[t_i, T]} + X_{t_n, \cd}\1_{[t_n, T]}.
\eeaa
Note that
\beaa
&\dis\wh\o \otimes_t X = \lim_{n\to\infty} \Big[\wh\o_{t\wedge \cd} + \sum_{i=1}^{n-1} X_{t,t_i} \1_{[t_i, t_{i+1})} + X_{t,\cd} \1_{[t_n, T]}\Big] = \lim_{n\to\infty} X^{n,n};&\\
&\dis\wh\o_{t\wedge \cd} \otimes_{t+\d} X = \wh\o_{t\wedge \cd}  + X_{t_n,\cd}\1_{[t_n, T]} = X^{n,0}.&
\eeaa
 Then, denoting $X^{n,m,\th}:=X^{n,m}+ \th X_{t_{m}, t_{m+1}} \1_{[t_{m+1}, T)}$,
\beaa
&& \dbE^{\dbP_0}\Big[ g(\wh\o \otimes_t X) - g(\wh\o_{t\wedge \cd} \otimes_{t+\d} X)\Big]\\
& =& \lim_{n\to\infty}\dbE^{\dbP_0} \Big[g(X^{n,n})-g(X^{n,0}) \Big]= \lim_{n\to\infty} \sum_{m=1}^{n} \dbE^{\dbP_0} \Big[g(X^{n,m}) - g(X^{n,m-1})\Big] \\
&=&   \lim_{n\to\infty} \sum_{m=1}^{n} \dbE^{\dbP_0} \Big[g(X^{n,m-1} + X_{t_{m-1}, t_m}\1_{[t_{m}, T]}) - g(X^{n,m-1})\Big]  \\
&=&\lim_{n\to\infty} \sum_{m=1}^{n} \dbE^{\dbP_0} \Big[D_{t_m} g(X^{n,m-1}) X_{t_{m-1}, t_m}  + {1\over 2}D^2_{t_m} g(X^{n,m-1}) X_{t_{m-1}, t_m}^2 \\
&&\qq \qq+ {1\over 2} \big[D^2_{t_m} g(X^{n,m-1,\th_m}) -D^2_{t_m} g(X^{n,m-1})\big]   X_{t_{m-1}, t_m}^2\Big],
\eeaa
for some random variable $\th_m$ taking values in $[0,1]$.
Note that, under $\dbP^0$, $X_{t_{m-1}, t_m}$ and $X^{n,m-1}$ are independent. Then
\beaa
&&\dbE^{\dbP_0} \Big[D_{t_m} g(X^{n,m-1}) X_{t_{m-1}, t_m} \Big]=0,\\
& &\dbE^{\dbP_0} \Big[D^2_{t_m} g(X^{n,m-1}) X^2_{t_{m-1}, t_m}\Big]= \dbE^{\dbP_0} \Big[D^2_{t_m} g(X^{n,m-1})\Big] [t_m- t_{m-1}],
\eeaa
and
\beaa
&&\Big|\dbE^{\dbP_0} \Big[ \big[D^2_{t_m} g\big(X^{n,m-1, \th_m} \big) -D^2_{t_m} g(X^{n,m-1})\big]    X_{t_{m-1}, t_m}^2\Big]\Big|\\
&\le& C \Big( \dbE^{\dbP_0} \big[ \sup_{0\le \th \le 1}\big|D^2_{t_m} g\big(X^{n, m-1,\th} \big) -D^2_{t_m} g(X^{n, m-1})\big|^2\big]\dbE^{\dbP_0}[|X_{t_{m-1}, t_m}|^4]\Big)^{1\over 2}\\
&\le&  C \Big( \dbE^{\dbP_0} \Big[\sup_{0\le \th \le 1} \big|D^2_{t_m} g\big(X^{n,m-1,\th} \big) -D^2_{t_m} g(X^{n,m-1})\big|^2\Big]\Big)^{1\over 2} [t_m-t_{m-1}].
\eeaa
Then, by the assumed regularity and the dominated convergence theorem, we can easily show that
\beaa
&&\dis\dbE^{\dbP_0}\Big[ g(\wh\o \otimes_t X) - g(\wh\o_{t\wedge \cd} \otimes_{t+\d} X)\Big]  ={1\over 2} \lim_{n\to\infty}  \sum_{m=1}^{n}  \dbE^{\dbP_0} \Big[D^2_{t_m} g(X^{n,m-1})\Big] [t_m- t_{m-1}] \\
&&\dis={1\over 2} \int_t^{t+\d}  \dbE^{\dbP_0} \Big[D^2_s g\Big(\wh\o_{t\wedge \cd} + X_{t,\cd} \1_{[t, s]} + X_{t,s}\1_{[s, T]} + X_{t+\d, \cd}\1_{[t+\d, T]}\Big)\Big]ds.
\eeaa
This implies
\beaa
\lim_{\d\to 0} {1\over \d} \dbE^{\dbP_0}\Big[ g(\wh\o \otimes_t X) - g(\wh\o_{t\wedge \cd} \otimes_{t+\d} X)\Big] = {1\over 2}  \dbE^{\dbP_0} \Big[D^2_t g(\wh\o\otimes_t X)\Big].
\eeaa
Similar results hold for $f$. Then
\beaa
&&\pa_t u(t, \wh \o) := \lim_{\d\to 0} {u(t+\d, \wh \o_{\cd\wedge t}) - u(t,\wh\o)\over \d} \\
&&=  - \dbE^{\dbP_0} \Big[{1\over 2}D^2_t g(\wh\o\otimes_t X) +{1\over 2} \int_t^T D^2_t f(s,\wh\o\otimes_t X) ds + f(t, \wh\o\otimes_t X) \Big].
\eeaa
Note that $\pa_t V(t, \wh \mu) = \dbE^{\wh \mu}[\pa_t u(t, \wh X)].$ Then one can easily verify the result.
\qed

{\color{black}
\subsection{Stochastic control under probability distortion}
In this subsection we study another application of the parabolic master equation. Probability distortion is an important tool in behavioral finance, in particular the prospect theory, see the survey paper Zhou \cite{Zhou} and the references therein. We say a function $\k: [0, 1]\to [0,1]$ is a {\it probability distortion function} if $\k$ is continuous, strictly increasing, and $\k(0) = 0, \k(1) = 1$. Given a random variable $\xi\ge 0$, introduce a nonlinear expectation:
\bea
\label{cEk}
\cE[\xi] := \int_0^\infty \k\big(\dbP(\xi \ge x)\big) dx.
\eea
The following properties are straightforward:

$\bullet$ If $\k(p) = p$, then $\cE[\xi] = \dbE[\xi]$.

$\bullet$ In general, $\cE$ is nonlinear: $\cE[\xi_1 + \xi_2] \neq \cE[\xi_1] + \cE[\xi_2]$.

$\bullet$ $\cE$ is law invariant: if $\cL_{\xi_1} = \cL_{\xi_2}$, then $\cE[\xi_1] = \cE[\xi_2]$.

\no In prospect theory, typically $\k$ is in reverse $S$-shape, namely concave around $0$ and convex around $1$. Indeed, assume $\k$ is smooth and $\xi$ has density $f(x)$, then it follows from the integration by parts formula that
\beaa
\cE[\xi] = \int_0^\infty \k'\big(\dbP(\xi \ge x)\big)   f(x) ~ x dx. 
\eeaa
Note that $\k'(p)$ is large for $p$ around $0$ and $1$, so at above integration the probability density  $f(x)$ is amplified by $\k'$ when $x$ is around $0$ and $\infty$, which is referred as probability distortion. 

Mathematically, the main challenge in this framework is the time inconsistency in the following sense. Assume $X$ is a Markovian process, $g$ is a positive function, and denote
\bea
\label{distortionu}
u(t,x) :=  \cE[ g(X^{t,x}_T)] :=  \int_0^\infty \k\big(\dbP(g(X_T) \ge y | X_t=x)\big) dy.
\eea
Then the flow property (hence the DPP when controls are involved) fails:
\beaa
\cE[g(X_T)] \neq \cE[ u(t, X_t) ].
\eeaa
In particular, the above function $u$ does not satisfy any PDE.
}

One remedy for the above time inconsistency is to consider $\cL_{X_t}$, instead of $X_t$, as the state variable. Then the expected PDE becomes a master equation. To be precise, assume $d=1$ and recall the $\dbP^{t,\mu}_0$ in Example \ref{eg-heat} and recall Remark \ref{rem-state}. 

\begin{eg}
\label{eg-distortion}
Assume the distortion function $\k\in C^1([0,1])$ and $g\in C^0_b(\dbR; \dbR_+)$. Define
\bea
\label{distortion-V}
V(t,\mu) := \int_0^\infty \k\Big( \dbP^{t,\mu}_0(g(X_T) \ge y)\Big) dy,\q (t,\mu)\in\Th.
\eea
Then $V$ is state dependent: $V(t,\mu) = V(t, \mu_t)$, and  $V \in C^{1,1,1}([0, T] \times \cP_2(\dbR))$ satisfies the following master equation:
\bea
\label{distortionMaster}
\pa_t V(t, \mu) + {1\over 2}  \dbE^\mu\Big[ \pa_x\pa_\mu V(t, \mu, X_t))\Big]=0,\q V(T, \mu) =  \int_0^\infty \k\Big(\mu(g(X_T) \ge y)\Big)dy.
\eea
\end{eg}
\proof It is clear that
\beaa
V(t,\mu) = \int_0^\infty \k\Big( \dbE^\mu[I(t, X_t, y)] \Big)dy, \q I(t,x, y) :=\int_{g^{-1}([y, \infty))} {1\over \sqrt{2\pi(T-t)}} e^{-{(x-z)^2\over 2(T-t)}} dz. 
\eeaa
One can easily check that
\beaa
&&\pa_t V(t,\mu) =  \int_0^\infty \k' \Big( \dbE^\mu[I(t, X_t, y)] \Big)   \dbE^\mu[ \pa_t I(t, X_t, y)] dy;\\
&&\pa_\mu V(t, \mu, x) = \int_0^\infty \k' \Big( \dbE^\mu[I(t, X_t, y)] \Big) \pa_x I(t, x, y) dy;\\
&&\pa_x\pa_\mu V(t, \mu, x) = \int_0^\infty \k' \Big( \dbE^\mu[I(t, X_t, y)] \Big) \pa_{xx} I(t, x, y) dy.
\eeaa
It is clear that $\pa_t I(t, x, y) + {1\over 2} \pa_{xx} I(t,x,y) = 0$. This implies  \reff{distortionMaster} straightforwardly.
\qed

 \begin{rem}
 \label{rem-distortion}
 {\rm (i) While $\cE$ is a nonlinear function, the master equation \reff{distortionMaster} is actually linear. The nonlinearity is only in the terminal condition: the mapping $\mu \mapsto V(T,\mu)$ is nonlinear in the sense that $V(T, \cL_{\xi_1 + \xi_2}) \neq V(T, \cL_{\xi_1}) + V(T, \cL_{\xi_2})$. 
  
 (ii)  In Ma, Wong, \& Zhang \cite{MWZ}, we introduced a dynamic distortion function $\k(t, x, p)$ to recover the flow property for the corresponding $u$ in \reff{distortionu} in some special cases. In Example \ref{eg-distortion},  we instead raise the "dimension" of the state space from $\dbR$ to $\cP_2(\dbR)$ so as to recover the flow property. We remark that this approach works for many time inconsistent problems, including those in Subsection \ref{sect-deterministic}.   However, in practice it may not be reasonable to use $V(t,\mu)$ as one's utility at time $t$, because by that time one observes a path of $X_{t\wedge \cd}$, then it is not reasonable to consider the whole distribution of $X_{t\wedge \cd}$ which involves other paths. Nevertheless, when one observes the value $X_0 = x_0$ at time $t=0$, the master equation \reff{distortionMaster} provides a nice characterization for the value $V(0, \d_{\{x_0\}})$.
 \qed}
 \end{rem}

{\color{black}
We next extend the above discussion to control problems under probability distortion, which to our best knowledge is new in the literature. Recall the $\cA$ in \reff{control-cA}, and similarly as \reff{control-Pa} we determine $\dbP^{t,\mu,\a}$  by the following controlled SDE on $[t, T]$:
 \bea
\label{distortion-Pa}
\dis X_s = X_t + \int_t^s b(r, X, \a_r(X)) dr +    \int_t^s \si(r, X, \a_r(X)) dB^\a_r,  \q  \dbP^{t,\mu,\a}\mbox{-a.s.}
\eea
where $b, \si, \a$ are all $\dbF$-adapted. Our value function is: given $g: \O\to [0, \infty)$,
\bea
\label{distortion-V1}
V(t,\mu) := \sup_{\a\in \cA}  \int_0^\infty \k\Big( \dbP^{t,\mu,\a}(g(X) \ge y)\Big) dy,\q (t,\mu)\in\Th.
\eea
Note that $\tilde g(\mu):= \int_0^\infty \k\Big( \mu(g(X_\cd) \ge y)\Big) dy$ is actually a deterministic function of $\mu$.  Then by considering $f=0$ and terminal condition $\tilde g$ in Theorem \ref{thm-meancontrol}, we obtain
\begin{cor}
\label{cor-distortion} 
Assume $b, \si,  g$ satisfy standard technical conditions, $\k$ is a probability distortion function, and define $V$ by \reff{distortion-Pa} and \reff{distortion-V1}. Then

(i)  The following dynamic programming principle holds: 
\bea
\label{distortion-DPP}
V(t_1, \mu) = \sup_{\a\in\cA} V(t_2, \dbP^{t_1, \mu, \a}),\q t_1 < t_2.
\eea

(ii) Assume $V\in C^{1,1,1}(\Th)$, then $V$ satisfies the path dependent  master equation:
\bea
\label{distortion-master1}
\left.\ba{c}
\dis \pa_t V(t,\mu) + \dbE^\mu\Big[ \sup_{a\in A} \big[{1\over 2} \pa_\o \pa_\mu V(t, \mu, X) \si^2(t, X, a) + \pa_\mu V(t, \mu, X)  b(t, X,  a)\big]\Big] =0,  \\
\dis V(T, \mu) = \int_0^\infty \k\Big( \mu(g(X) \ge y)\Big) dy.
 \ea\right.
 \eea 
 
 (iii)  Assume further that the Hamiltonian  in \reff{distortion-master1} has an optimal argument $a^* = I(t, \o, \mu)$, where $I$ is uniformly Lipschitz continuous in $\o$, and the following McKean-Vlasov SDE has a solution:
 \bea
 \label{distortion-X}
 X^*_t = x_0 + \int_0^t b(s, X^*, I(s,  X^*, \cL_{X^*}) )ds + \int_0^t \si(s, X^*,   I(s,  X^*, \cL_{X^*}) )dB_s.
 \eea
 Then $\a^*_t := I(t, X^*, \cL_{X^*})$ is an optimal control to the problem $V(0, \d_{\{x_0\}})$ in  \reff{distortion-V1}.
 \end{cor} 
}

\section{Viscosity solution of master equations}
\label{sect-Viscosity}
\setcounter{equation}{0}
We emphasize again that the smoothness of $V$ required in Theorem \ref{thm-meancontrol} is very difficult to verify. 
In this section we propose a notion of viscosity solution for master equation \reff{master}, which requires less regularity, and establish its basic properties.

\subsection{Definition of viscosity solutions}
For $(t,\mu)\in \Th$ and constant $L>0$,  {\color{black}let $\cP_L(t,\mu)$ denote the set of $\dbP\in\cP_2$ such that $\dbP_{[0,t]}=\mu_{[0,t]}$ and $X_{[t, T]}$ is a $\dbP$-semimartingale with drift and diffusion characteristics bounded by $L$, in the spirit of the $\wh\cP_L$ introduced in the beginning of Subsection \ref{sect-FIto}.}  Note that we do not require $X$ to be a $\mu$-semimartingale on $[0, t]$. {\color{black}The following simple estimates will be used frequently in the paper:  for any $(t,\mu)\in \Th$, $\d\in [0, T-t]$, and $L>0$, $p\ge 1$,
\bea
\label{pLest}
 \sup_{\dbP\in \cP_L(t,\mu)} \dbE^\dbP\big[\sup_{t\le s\le t+\d}|X_{t,s}|^p\big] \le C_{p, L} \d^{p\over 2}.
\eea
}

The following compactness result is the key for our viscosity theory.

\begin{lem}
\label{lem-wcpt}
For any $(t,\mu)\in \Th$ and $L>0$, the set $[t, T]\times \cP_L(t,\mu)$ is  compact under $\cW_2$. 
\end{lem}
\proof We first show that  $\cP_L(t,\mu)$ is compact. Let $\{\dbP^n\}_{n\ge 1}\subset \cP_L(t,\mu)$. By Zheng \cite{W.A.Zheng85} Theorem 3, $\cP_L(t,\mu)$ is weakly compact, then there exist a convergent subsequence, and without loss of generality we assume $\dbP^n \to \dbP\in \cP_L(t,\mu)$ weakly. Note that 
\beaa
\|X\| \le \|X_{t\wedge \cd}\| + \sup_{t\le s\le T}|X_{t,s}| \le 2\Big[[\|X_{t\wedge \cd}\|]~ \vee ~[ \sup_{t\le s\le T}|X_{t,s}|]\Big].
\eeaa
{\color{black} Since $\dbP^n = \mu$ on $\cF _t$, $\|X_{t\wedge \cd}\|$ has the same distribution under $\dbP^n$ and  $\mu$. Moreover, since $\dbP_n\in \cP_L(t,\mu)$,  for any $R>0$ and any $n$, by \reff{pLest} (with $\d = T-t$) we have}
\beaa
&& \dbE^{\dbP_{n}}\Big[\|X\|^2\1_{\{\|X\|\ge R\}}\Big] \le  4 \dbE^{\dbP_{n}}\Big[\|X_{t\wedge \cd}\|^2\1_{\{\|X_{t\wedge \cd}\|\ge {R\over 2}\}} +  \sup_{t\le s\le T}|X_{t,s}|^2\1_{\{\sup_{t\le s\le T}|X_{t,s}|\ge {R\over 2}\}}\Big]\\
 &&=4 \dbE^\mu\Big[\|X_{t\wedge \cd}\|^2\1_{\{\|X_{t\wedge \cd}\|\ge {R\over 2}\}}\Big] + 4 \dbE^{\dbP_{n}}\Big[\sup_{t\le s\le T}|X_{t,s}|^2\1_{\{\sup_{t\le s\le T}|X_{t,s}|\ge {R\over 2}\}}\Big]\\
 &&\le 4 \dbE^\mu\Big[\|X\|^2\1_{\{\|X\|\ge {R\over 2}\}}\Big] + {8\over R} \dbE^{\dbP_{n}}\Big[\sup_{t\le s\le T}|X_{t,s}|^3\Big] \le 4 \dbE^\mu\Big[\|X\|^2\1_{\{\|X\|\ge {R\over 2}\}}\Big] + {C_L\over R}.
 \eeaa
 Thus, by the dominated convergence theorem under $\mu$,
 \beaa
 \lim_{R\to\infty} \sup_{n\ge 1}  \dbE^{\dbP_{n}}\Big[\|X\|^2\1_{\{\|X\|\ge R\}}\Big]  \le  4 \lim_{R\to\infty} \dbE^\mu\Big[\|X\|^2\1_{\{\|X\|\ge {R\over 2}\}}\Big] = 0.
 \eeaa
Then it follows from  Carmona \& Delarue \cite{CD1} Theorem 5.5  that $\lim_{n\to\infty}\cW_2(\dbP^{n}, \dbP) =0$.

Next, let $(t_n, \dbP_n)\in [t, T]\times \cP_L(t,\mu)$.  By the compactness of $[t, T]$ and $\cP_L(t,\mu)$, we may assume without loss of generality that $t_n \to t^*$ and $\dbP_n\to \dbP$. Then
\beaa
&&\cW_2\Big((t_n, \dbP_n), (t^*, \dbP)\Big) \le \cW_2\Big((t_n, \dbP_n), (t^*, \dbP_n)\Big) + \cW_2\Big((t^*, \dbP_n), (t^*, \dbP)\Big)\\
&&\le \Big( |t_n-t^*| + {\color{black}\dbE^{\dbP_n}[\|X_{t_n\wedge \cd} - X_{t^*\wedge \cd}\|^2]}\Big)^{1\over 2}+ \cW_2( \dbP_n, \dbP)\\
&&\le C|t_n-t^*|^{1\over 2} + \cW_2( \dbP_n, \dbP)\to 0,\q\mbox{as}~n\to\infty.
\eeaa
This implies that $[t, T]\times \cP_L(t,\mu)$ is also compact.
\qed

For the viscosity theory, another crucial thing is the functional  It\^{o} formula \reff{functionalIto2}. For this purpose, we shall weaken the regularity requirement for the test functions, which will make the theory more convenient. 

\begin{defn}
\label{defn-C12Ito}
Let $0\le t_1< t_2\le T$ and $\cP \subset \cP_2$ such that $X$ is a semimartingale on $[t_1, t_2]$ under each $\dbP\in \cP$. We say $V\in C^{1,1,1}([t_1, t_2]\times \cP)$ if $V\in C^0([t_1, t_2]\times \cP)$ and there exist $\pa_t V \in C^0([t_1, t_2]\times \cP), \pa_\mu V, \pa_\o\pa_\mu V \in C^0([t_1, t_2]\times \cP\times \O)$ with appropriate dimensions, such that the functional It\^{o} formula \reff{functionalIto2} holds true on $[t_1, t_2]$ under every $\dbP\in \cP$.

 Moreover, let $C^{1,1,1}_b([t_1, t_2]\times \cP)$ denote the subset of $C^{1,1,1}([t_1, t_2]\times \cP)$ such that $\pa_t V$ is bounded and, for some constants $C\ge 0$,
 \beaa
 |\pa_\mu V(t, \mu, \o)| + |\pa_\o\pa_\mu V(t, \mu, \o)| \le C [ 1+\|\o\|],\q \mbox{for $\dbP$-a.e. $\o$, and for all $\dbP\in \cP$}.
 \eeaa
  \end{defn}

\begin{rem}
\label{rem-C12Ito}
{\rm (i) By  Theorem \ref{thm-Ito2},  $C^{1,1,1}_b(\Th) \subset C^{1,1,1}([t_1, t_2]\times \cP_L(t_1, \mu))$ for all $(t_1, t_2)$, $L$, and $\mu\in \cP_2$, and the derivatives $\pa_t V, \pa_\mu V, \pa^{sym}_\o\pa_\mu V$ are consistent.

(ii) Following the same arguments as in Theorem \ref{thm-Ito2} (i), for $V\in C^{1,1,1}([t_1, t_2]\times \cP_L(t_1, \mu))$, $\pa_t V, \pa_\mu V, \pa^{sym}_\o\pa_\mu V$ are unique. Since by Remark \ref{rem-parabolic}, $G$ depends on $\G$ only through $\G^{sym}$, so the uniqueness of $ \pa^{sym}_\o\pa_\mu V$ is sufficient for our purpose.

(iii) When $\cP$ is compact, e.g. $\cP = \cP_L(t, \mu)$, the continuity implies uniform continuity as well as boundedness. In particular, in this case $V$ and $\pa_t V$ are automatically bounded and  the linear growth of $\pa_\mu V, \pa_\o\pa_\mu V$ in $\o$ is also a mild requirement. 
\qed}
\end{rem}
For a function $V:\Th\to\dbR$, we now introduce the following set of test functions:
\bea
\label{cA}
\left.\ba{lll}
&&\cA^L_\d V(t,\mu):= \big\{\f\in C^{1,1,1}_b\big( [t, t+\d] \times \cP_L(t,\mu)\big): (\f-V)(t,\mu)=0\big\};\\
&&\dis \underline{\cA}^L V(t,\mu):= \bigcup_{0<\d\le T-t}\big\{\f\in \cA^L_\d V(t,\mu): \inf_{(s, \dbP) \in [t, t+\d] \times \cP_L(t,\mu)}(\f-V )(s,\dbP) = 0\big\};\\
&&\dis\overline{\cA}^L V(t,\mu):=\bigcup_{0<\d\le T-t}\big\{\f\in \cA^L_\d V(t,\mu): \sup_{(s, \dbP) \in [t, t+\d] \times \cP_L(t,\mu)}(\f-V )(s,\dbP) =0\big\}.
\ea\right.
\eea

\begin{defn}
\label{defn-viscosity}
Let $V \in C^0(\Th)$. 

(i) We say $V$ is an $L$-viscosity subsolution (resp. supersolution) of \reff{master} if  $\sL\f(t,\mu) \ge (\mbox{resp.} ~\le)  ~ 0$ for all $(t,\mu)\in \Th$ and all $\f\in \underline{\cA}^L V(t,\mu)$ (resp. $\overline{\cA}^L V(t,\mu)$).

(ii)  We say $V$ is an $L$-viscosity solution of \reff{master} if it is both an $L$-viscosity subsolution and an $L$-viscosity supersolution, and $V$ is a viscosity solution if it is  an $L$-viscosity solution for some $L>0$.
\end{defn}

\begin{rem}
\label{rem-cPL}
{\rm (i) Our main idea here is to use $\cP_L(t, \mu)$ in \reff{cA}, which by Lemma \ref{lem-wcpt} is compact under $\cW_2$ and in the meantime is large enough in most applications we are interested in.  This is in the same spirit as our notion of viscosity solutions for path dependent PDEs, see Ekren, Keller, Touzi, \& Zhang \cite{EKTZ} and Ekren, Touzi, \& Zhang \cite{ETZ1, ETZ2}.

(ii) When $V$ is state dependent: $V(t, \mu) = V(t, \mu_t)$, the above definition still works. However, in this case it is more convenient to change the test functions $\f$ to be state dependent only. In particular, we shall revise \reff{cA}  as follows:

$\bullet$ $ \underline{\cA}^L V(t,\mu)$ and $\overline{\cA}^L V(t,\mu)$ become  $ \underline{\cA}^L V(t,\mu_t)$ and $\overline{\cA}^L V(t,\mu_t)$;

$\bullet$ $\cP_L(t, \mu)$ becomes $\cP_L(t, \mu_t)$ where the initial constraint is relaxed to $\dbP_t = \mu_t$;

$\bullet$ the extremum is about $[\f-V](s, \dbP_s)$ for $(s, \dbP) \in [t, t+\d]\times \cP_L(t, \mu_t)$. 

{\color{black}
(iii) In the state dependent case, if we work on torus $\dbT^d$ instead of $\dbR^d$ (namely the state process $X$ takes values in $\dbT^d$), then the following $\d$-neighborhood   is compact under $\cW_2$:
\bea
\label{Dd}
D_\d(t,\mu_t) :=\big\{(s, \dbP_s)\in [t, t+\d]\times \cP_2(\dbT^d): \cW_2(\dbP_s, \mu_t)\le \d\big\}.
\eea
and we expect the main results in this paper will remain true by replacing $\cP_L(t, \mu_t)$ with $D_\d(t,\mu_t)$. However, we lose such compactness on $\cP_2(\dbR^d)$,  for example, $ \mu_n:= {1\over n} \d_{\{n^2\}} + [1-{1\over n}]\d_{\{0\}}\in \cP_2(\dbR)$ converges to $\d_{\{0\}}$ weakly, but not under $\cW_2$. So our definition of viscosity solution is novel even in the state dependent case.}
\qed}
\end{rem}

\subsection{Some equivalence results}
\begin{thm}[Consistency]
\label{thm-consistency}
Let Assumption \ref{assum-G} hold and $V \in C^{1,1,1}_b(\Th)$. Then $V$ is a viscosity solution (resp. subsolution, supersolution) of master equation \reff{master} if and only if it is a classical solution  (resp. subsolution, supersolution) of master equation \reff{master}.
\end{thm}
\proof We shall only prove the equivalence of the subolution property. If $V$ is a viscosity subsolution, note that $V$ itself is in  $\underline{\cA}^L V(t,\mu)$, then clearly $\sL V(t, \mu) \ge 0$ and thus is a classical subsolution. Now assume $V$ is a classical subsolution. Fix $(t, \mu)\in \Th$ and $\f\in \underline{\cA}^L V(t,\mu)$ for some $L \ge L_0$, where $L_0$ is the Lipschitz constant  in Assumption \ref{assum-G}. Given $\cF_t$-measurable random variables $b_t, \si_t$ with $|b_t|, {1\over 2}|\si_t|^2\le L$,  let $\dbP\in \cP_L(t, \mu)$ be such that $X_{t, \cd}$ is a $\dbP$-semimartingale with drift $b_t$ and volatility $\si_t$.  Then, denoting $\psi := \f - V$,
\beaa
0 &\le& \psi(t+\d, \dbP) - \psi(t, \mu) \\
&=& \int_t^{t+\d} \Big[\pa_t \psi(s, \dbP)  + \dbE^\dbP\big[ b_t\cd \pa_\mu \psi(s, \dbP, X) + {1\over 2} \si_t\si^\top_t : \pa_\o \pa_\mu \psi(s, \dbP, X)\big]\Big] ds.
\eeaa
Divide both sides by $\d$ and send $\d\to 0$, we obtain
\beaa
0 &\le& \pa_t \psi(t, \mu)  + \dbE^\mu\Big[ b_t\cd \pa_\mu \psi(t, \mu, X) + {1\over 2} \si_t\si^\top_t : \pa_\o \pa_\mu \psi(t, \mu, X)\Big].
\eeaa
Set $y := V(t,\mu) = \f(t, \mu)$, $Z_1 :=  \pa_\mu \f(t, \mu, \cd)$, $Z_2 :=  \pa_\mu V(t, \mu, \cd)$, $\G_1 :=  \pa_\o \pa_\mu \f(t, \mu, \cd)$, and $\G_2 :=  \pa_\o \pa_\mu V(t, \mu, \cd)$. Let $b_t$ and $\si_t$ be as in \reff{GLip}, then
\beaa
&0\le  \pa_t \f(t, \mu) - \pa_t V(t, \mu)  + G(t, \mu, y, Z_1, \G_1) -G(t, \mu, y, Z_2, \G_2) = \sL \f(t, \mu) - \sL V(t, \mu),  &
\eeaa
and thus $\sL \f(t, \mu) \ge \sL V(t, \mu) \ge 0$. That is, $V$ is a viscosity subsolution.
\qed

As in the standard viscosity theory,  we may alternatively define viscosity solutions via semi-jets.  For $t\in [0, T]$, $y\in \dbR$, $v\in \dbR$ and $\cF_t$-measurable $Z, \G \in C^0(\O; \dbR^d) \times C^0(\O; \dbR^{d\times d})$ with $|Z(\o)|+|\G(\o)|\le C[ 1+\|\o\|]$ for some $C>0$, define paraboloids as follows: 
 \bea
 \label{PPDE-general-phi}
 \phi^{t,y, v,Z,\G}(s, \dbP)   :=  y+  v[s- t] + \dbE^\dbP \Big[ Z \cd X_{t,s}  + {1\over 2}  \G: [ X_{t,s}X_{t,s}^\top ] \Big], \q (s, \dbP)\in [t, T]\times \cP_2.
 \eea
  For any $(t,\mu)\in \Th$,  it is clear that $ \phi^{t,y, v,Z,\G} \in C^{1,1,1}_b([t, T]\times \cP_L(t,\mu))$ with: 
\bea
\label{phi0}
\mbox{for}~\f:=  \phi^{t, V(t,\mu), v,Z,\G}:\q \pa_t \f (t, \mu) = v,~\pa_\mu \f (t, \mu, \cd) = Z,~ \pa^{sym}_\o \pa_\mu \f(t, \mu, \cd)  = \G^{sym}.
 \eea
We then introduce the corresponding subjets and superjets: for $L>0$,
 \bea
 \label{PPDE-general-jets}
 \left.\ba{c}
\dis \ol{\cJ}^L V(t,\mu) :=  \Big\{ (v,Z,\G):  \phi^{t,V(t,\mu), v,Z, \G} \in \ol\cA^L V(t,\mu) \Big\};\ms\\
 \dis \ul{\cJ}^L V(t,\mu) :=  \Big\{ (v,Z, \G):   \phi^{t,V(t,\mu), v, Z, \G} \in \ul\cA^L V(t,\mu) \Big\}.
 \ea\right.
 \eea

\begin{thm}
\label{thm-jets}
Let Assumption \ref{assum-G} hold and $V \in C^0(\Th)$.  Then $V$ is an $L$-viscosity supersolution (resp. subsolution) of master equation \reff{master} if and only if: for any $(t,\mu)\in \Th$,
 \bea
 \label{equiv}
\!\!\!\! v +G(t,\mu, V(t, \mu), Z,\G)\le ~(\mbox{resp.} ~\ge )~ 0,~ \forall  (v, Z, \G)\in \ol \cJ^L V(t,\mu) ~(\mbox{resp.} ~\ul{\cJ}^L V(t,\mu)).
 \eea
\end{thm}

\proof "$\Longrightarrow$" Assume $V$ is an $L$-viscosity supersolution at $(t,\mu)$. For any  $(v, Z, \G)\in \ol \cJ^L V(t,\mu)$, since $\phi^{t,V(t,\mu), v, Z, \G} \in \ol \cA^L V(t,\mu)$, then it follows from the viscosity property of $V$ and \reff{phi0} that 
$
0\ge \sL\f(t, \mu)  = v +G(t,\mu, V(t,\mu), Z, \G).
$

"$\Longleftarrow$"  Assume \reff{equiv} holds at $(t,\mu)$ and $\f\in \ol \cA^L V(t,\mu)$ with  corresponding $\d$. Denote
\bea
\label{PPDE-general-e}
v := \pa_t V(t, \mu),~~ v_\e:= v-\e(1+2L),~~ Z:= \pa_\mu \f(t, \mu,\cd),~~ \G:= \pa_\o\pa_\mu \f(t,\mu,\cd),\q\forall \e>0.
\eea
Then, for any $(s, \dbP) \in [t, t+\d] \times \cP_L(t,\mu)$,
\beaa
&&\phi^{t, V(t, \mu), v_\e, Z, \G}(s, \dbP) - \f(s,\dbP) = \int_t^s  [v_\e - \pa_t \f(r, \dbP)] dr \\
&&\q+ \dbE^\dbP \Big[\int_t^s [Z + \G X_{t,r} - \pa_\mu \f(r, \dbP,\cd)]  \cd dX_r  + {1\over 2}\int_t^s [\G - \pa_{\o}\pa_\mu \f(r, \dbP, \cd)]: d\la X\ra_r\Big].
\eeaa
By choosing $\d>0$ small, we may assume without loss of generality that 
\bea
\label{PPDE-general-smallerror}
|\pa_t \f(s, \dbP) -v| \le \e, ~ \dbE^\dbP\Big[|\pa_\mu \f(s,\dbP) -Z-\G X_{t,s}|\Big]\le \e,~   \dbE^\dbP\Big[ |\pa_\o\pa_\mu \f(s, \dbP) -\G|\Big]\le \e,
\eea
for all $(s, \dbP) \in [t, t+\d] \times \cP_L(t,\mu)$. Then, 
\beaa
 \phi^{t, V(t,\mu), v_\e,Z, \G}(s, \dbP) - \f(s,\dbP) \le [s-t]\Big[ v_\e - v + \e + L \e +  L\e\Big] = 0.
\eeaa
Since  $\f \in  \ol\cA^L V(t,\mu)$, this implies immediately that  $(v_\e, Z, \G) \in \ol\cJ^L V(t,\mu)$. By our assumption we have
$v_\e + G(t,\mu, V(t,\mu), Z, \G) \le 0$. Send $\e\to 0$, we obtain $\sL\f (t,\mu) = v+ G(t,\mu, V(t,\mu),  Z,\G) \le 0$. That is, $V$ is an $L$-viscosity supersolution at $(t,\mu)$.
\qed

\begin{rem}
\label{rem-jets}
{\rm Technically speaking, since we can use the semi-jets to define viscosity solutions, our viscosity theory does not require the functional It\^{o} formula. Instead, it is sufficient to have the It\^{o} formula for the paraboloids in \reff{PPDE-general-phi}. But nevertheless the functional It\^{o} formula is crucial for classical solutions and is interesting in its own right.
\qed}
\end{rem}

Finally, the following change variable formula is also important for comparison principle.

\begin{thm}
\label{thm-change}
Let Assumption \ref{assum-G} hold and $V \in C^0(\Th)$. For any constant $\l\in \dbR$, define 
\bea
\label{change}
\tilde V(t, \mu) := e^{\l t} V(t, \mu),\q \tilde G(t, \mu, y, Z, \G) := e^{\l t} G(t, \mu, e^{-\l t} y, e^{-\l t} Z, e^{-\l t} \G).
\eea
Then $V$ is an $L$-viscosity solution (resp. subsolution, supersolution) of master equation \reff{master} if and only if $\tilde V$ is an $L$-viscosity solution (resp. subsolution, supersolution) of the following master equation:
\bea
\label{changeMaster}
\pa_t \tilde V(t, \mu)  - \l \tilde V(t, \mu) + \tilde G(t, \mu, \tilde V, \pa_\mu \tilde V, \pa_\o\pa_\mu \tilde V)=0.
\eea
\end{thm}
\proof We shall only prove that the viscosity subsolution property of $V$ implies the viscosity subsolution property of $\tilde V$. The other implications follow the same arguments.  

Assume $V$ is an $L$-viscosity subsolution of \reff{master}. Let $(\tilde v, \tilde Z, \tilde \G) \in  \ul \cJ^L\tilde V(t,\mu)$ with corresponding $\d_0>0$.  Then, for any $(s, \dbP)\in [t, t+\d_0]\times \cP_L(t,\mu)$,
\beaa
\tilde V(t,\mu) + \tilde v [s-t] + \dbE^\dbP \Big[ \tilde Z \cd X_{t,s}  + {1\over 2}  \tilde \G: [ X_{t,s}X_{t,s}^\top ] \Big] \ge \tilde V(s, \dbP).
\eeaa
Thus
\beaa
&V(t,\mu) +  v [s-t] + \dbE^\dbP \Big[ Z \cd X_{t,s}  + {1\over 2}  \G: [ X_{t,s}X_{t,s}^\top ] \Big] \ge e^{\l(s-t)} V(s, \dbP),&\\
&\mbox{where}\q v:= e^{-\l t} \tilde v,\q Z:= e^{-\l t} \tilde Z,\q \G := e^{-\l t} \tilde \G.&
\eeaa
Note that $V$ is continuous and  $[t, t+\d_0]\times \cP_L(t,\mu)$ is compact, then $V$ is bounded and uniformly continuous. Thus
\beaa
e^{\l(s-t)} V(s, \dbP) = [1 + \l (s-t) + o(s-t)] V(s, \dbP) = V(s, \dbP) + \l V(t,\mu) [s-t] + o(s-t). 
\eeaa
Therefore, for any $\e>0$, there exists $\d\in (0, \d_0)$ such that, for $(s, \dbP)\in [t, t+\d]\times \dbP_L(t,\mu)$, 
\beaa
V(t,\mu) +  \big[v - \l V(t,\mu) + \e\big] [s-t] + \dbE^\dbP \Big[ Z \cd X_{t,s}  + {1\over 2}  \G: [ X_{t,s}X_{t,s}^\top ] \Big] \ge  V(s, \dbP).
\eeaa
This implies that $(v - \l V(t,\mu) + \e, Z, \G) \in  \ul \cJ^L V(t,\mu)$, and thus
\beaa
v - \l V(t,\mu) + \e + G(t, \mu, V(t, \mu), Z, \G) \ge 0.
\eeaa
Since $\e>0$ is arbitrary, we have
\beaa
v - \l V(t,\mu) + G(t, \mu, V(t, \mu), Z, \G) \ge 0.
\eeaa
This implies immediately that 
\beaa
\tilde v - \l\tilde V(t,\mu) + \tilde G(t, \mu, \tilde V(t, \mu),\tilde Z, \tilde \G) \ge 0.
\eeaa
That is, $\tilde V$ is an $L$-viscosity subsolution of \reff{changeMaster}.
\qed

\subsection{Stability}

 For any $(t, \mu, y, Z, \G)$ and $\d>0$,  denote
\bea
\label{stability-Od}
O^L_\d(t, \mu, y, Z, \G) := \Big\{ (s, \dbP, \tilde y, \tilde Z, \tilde \G): (s, \dbP)\in [t, t+\d]\times \cP_L(t,\mu),\nonumber\\
  |\tilde y-y|\le \d, \q \dbE^\dbP\big[|\tilde Z-  Z|^2 + |\tilde G - G|^2\big] \le \d^2\Big\}.
\eea

\begin{thm}
\label{thm-stability}  Let $L>0$, $G$ satisfy Assumption \ref{assum-G}, and $V \in C^0(\Th)$. Assume

(i) for any $\e>0$, there exist $G^\e$ and $V^\e\in C^0(\Th)$ such that $G^\e$ satisfies Assumption \ref{assum-G} and $V^\e$ is an $L$-viscosity subsolution of master equation \reff{master} with generator $G^\e$;

(ii) as $\e\to 0$, $(G^\e, V^\e)$ converge to $(G, V)$  locally uniformly in the following sense: for any $(t, \mu, y, Z, \G)$, there exists $\d>0$ such that, 
\bea
\label{localuniform}
\lim_{\e\to 0} \sup_{(s, \dbP, \tilde y, \tilde Z, \tilde \G) \in O^L_\d(t, \mu, y, Z, \G)}\Big[|[G^\e - G](s, \dbP, \tilde y, \tilde Z, \tilde \G) |+ |[V^\e - V] (s, \dbP)| \Big] = 0.
\eea
Then $V$ is an $L$-viscosity subsolution  of master equation \reff{master}  with generator $G$.
\end{thm}
\proof  Let $\f\in {\overline\cA}^{L} V(t,\mu)$  with corresponding $\d_0$. By \reff{localuniform} we may choose $\d_0>0$ small enough  such that 
\bea
\label{localuniform2}
&\dis\lim_{\e\to 0} \rho(\e,\d_0)=0, \q\mbox{where, denoting} ~( y_0, Z_0, \G_0):=(\f(t, \mu), \pa_\mu \f(t, \mu,\cd), \pa_\o\pa_\mu \f(t,\mu,\cd),&\nonumber\\
&\dis\rho(\e,\d ) :=  \sup_{(s, \dbP, y, Z, \G) \in O^L_\d(t, \mu, y_0, Z_0, \G_0) }\Big[|[G^\e - G](s, \dbP, y, Z, \G)] + |[V^\e - V](s, \dbP)|\Big].&
\eea
For $0 <\d \le \d_0$, denote $\f_\d(s,\dbP) := \f(s,\dbP) + \d [s-t]$. Then
\beaa
[\f_\d-V](t,\mu)&=& [\f-V](t,\mu)= 0 \\
&\le& \inf_{\dbP \in \cP_L(t,\mu)}   [\f-V](t+\d,\dbP) <  \inf_{\dbP \in \cP_L(t,\mu)}   [\f_\d-V](t+\d,\dbP).
\eeaa 
By \reff{localuniform2}, there exists $\e_\d>0$ small enough such that, for any $\e \le \e_\d$, 
\bea
\label{stabilityest1}
[\f_\d-V^\e](t,\mu)  <  \inf_{\dbP \in \cP_L(t,\mu)}   [\f_\d-V^\e](t+\d,\dbP).
\eea
Then there exists $(t^*, \dbP^*) \in [t, t+\d) \times \cP_L(t, \mu)$, which may depend on $(\e, \d)$, such that
\beaa
 c^*:= \inf_{(s,\dbP) \in [t, t+\d]\times \cP_L(t,\mu)}   [\f_\d-V^\e](s,\dbP) = [\f_\d-V^\e](t^*,\dbP^*)
 \eeaa
This implies immediately that 
\beaa
\f^\e_\d := \f_\d - c^*  \in   {\ul\cA}^{L} V^\e(t^*, \dbP^*).
\eeaa
Since $V^\e$ is a viscosity $L$-subsolution of master equation \reff{master} with generator $G^\e$, we have 
\bea
\label{stabilityest2}
0 &\le& \Big[\pa_t \f^\e_\d + G^\e(\cd, \f^\e_\d, \pa_\mu \f^\e_\d, \pa_\o\pa_\mu \f^\e_\d)\Big](t^*, \dbP^*)\nonumber\\
&=&  \Big[\pa_t \f + \d + G^\e(\cd, V^\e, \pa_\mu \f, \pa_\o\pa_\mu \f)\Big](t^*, \dbP^*)\nonumber\\
&\le&  \Big[\pa_t \f +  G(\cd, V^\e, \pa_\mu \f, \pa_\o\pa_\mu \f)\Big](t^*, \dbP^*) +\d+ \rho(\e, \d_0),
\eea
for $\e$ and $\d$ small enough. Now send $\d\to 0$, we get
\beaa
0 &\le& \Big[\pa_t \f +  G(\cd, V^\e, \pa_\mu \f, \pa_\o\pa_\mu \f)\Big](t, \mu) +\rho(\e, \d_0).
\eeaa
Send further $\e\to 0$ and then $\d_0\to 0$, we obtain the desired viscosity subsolution property of $V$ at $(t, \mu)$.
 \qed

\subsection{Partial comparison principle}
\begin{thm}[Partial Comparison Principle]
\label{thm-partial}
Let Assumption \ref{assum-G} hold,  $V^1$ be a viscosity subsolution and $V^2$ a viscosity supersolution of \reff{master}. If $V^1(T,\cd)\le V^2(T,\cd)$ and either $V^1\in C^{1,1,1}_b(\Th)$ or $V^2\in C^{1,1,1}_b(\Th)$, then $V^1\le V^2$.
\end{thm}
\proof We shall prove by contradiction. Denote $\D V := V^1- V^2$. Assume without loss of generality that $V^2\in C^{1,1,1}_b(\Th)$ and that   $c:=\D V(t,\mu)>0$ for some $(t,\mu)\in \Th$.  Define
\bea
\label{partial-c*}
c^*:= \sup_{(s, \dbP) \in [t, T] \times \cP_L(t,\mu)}\Big[\D V(s,\dbP)-\frac{c}{2(T-t)}(T-s)\Big].
\eea
 Note that $\D V$ is continuous and $[t, T]\times \cP_L(t,\mu)$ is compact, then there exists $(t^*, \dbP^*) \in  [t, T] \times \cP_L(t,\mu)$ such that 
 \beaa
\D V(t^*,\dbP^*)-\frac{c}{2(T-t)}(T-t^*) = c^*.
\eeaa
By considering $s=t$ in \reff{partial-c*} it is clear that $c^* \ge {c\over 2} >0$. Moreover, by the boundary condition that $\D V(T,\cd)\le 0$, we see that $t^* < T$.   Define
\beaa
\f(s,\dbP):=V^2(s,\dbP)+c^\ast+\frac{c}{2(T-t)}(T-s).
\eeaa
Then $\f(t^\ast,\dbP^\ast)=V^1(t^\ast,\dbP^\ast)$. Since $\cP_L (t^\ast, \dbP^\ast)\subseteq \cP_L(t,\mu)$, for any $s\ge t^\ast$ and $\dbP\in\cP_L({t^\ast}, \dbP^\ast)$, we have $\f(s,\dbP)\ge V^1(s,\dbP)$. This implies that $\f\in\underline{\cA}^LV^1(t^\ast,\dbP^\ast)$, and thus
\beaa
0 &\le& \sL \f(t^*, \dbP^*) = \pa_t \f(t^*, \dbP^*) + G\big(t^*, \dbP^*, \f(t^*, \dbP^*),\pa_\mu \f(t^*, \dbP^*,\cd), \pa_\o\pa_\m \f(t^*, \dbP^*, \cd) \big)\\
&=& \pa_t V^2(t^*, \dbP^*) - {c\over 2(T-t)} + G\big(t^*, \dbP^*, \f(t^*, \dbP^*),\pa_\mu V^2(t^*, \dbP^*,\cd), \pa_\o\pa_\m V^2(t^*, \dbP^*, \cd)\big).
\eeaa
By Theorem \ref{thm-change}, we can assume without loss of generality that $G$ is decreasing in $y$. Then, since $\f(t^*, \dbP^*) >V^2(t^*, \dbP^*) + c^* > V^2(t^*, \dbP^*)$, we have
\beaa
0 &\le&  \pa_t V^2(t^*, \dbP^*) - {c\over 2(T-t)} + G\big(t^*, \dbP^*, V^2(t^*, \dbP^*),\pa_\mu V^2(t^*, \dbP^*,\cd), \pa_\o\pa_\m V^2(t^*, \dbP^*, \cd)\big)\\
&=& \sL V^2(t^*, \dbP^*) - {c\over 2(T-t)} \le -{c\over 2(T-t)},
\eeaa
thanks to the classical supersolution property of $V^2$. This is a desired contradiction.
\qed

\subsection{Comparison principle}
\label{sect-comparison}
Given $g\in C^0(\cP_2, \dbR)$, define
\bea
\label{Vbar}
 \overline V(t,\mu) := \inf\big\{ \psi(t,\mu):  \psi \in \overline \cU_g\big\},\q \underline V(t,\mu) := \sup\big\{ \psi(t,\mu):  \psi  \in \ul \cU_g\big\},
\eea
where
\bea
&&\dis \cU:= \Big\{\psi: \Th\to \dbR~   \mbox{adapted, continuous in $\mu$,  \cad ~ in $t$,  and $\exists~ 0=t_0<\cds<t_n=T$ }\nonumber\\
&&\dis\qq\qq \mbox{such that}~ \psi \in C^{1,1,1}_b([t_i, t_{i+1})\times \cP_L(t_i, \mu))~\mbox{for any}~t_i, \mu\in \cP_2, L>0\Big\};\nonumber\\
&&\dis \overline\cU_g := \Big\{\psi \in \cU: \psi(T, \cd)\ge g,~ \mbox{and for the corresponding}~\{ t_i\},  \psi_{t_i} \le \psi_{t_i-},\\
&&\dis \mbox{and $\psi$ is a classical supersolution of master equation \reff{master} on each $[t_{i-1}, t_i)$}   \Big\};\nonumber\\
&&\dis \underline\cU_g := \Big\{\psi \in \cU: \psi(T, \cd)\le g,~ \mbox{and for the corresponding}~ \{t_i\},  \psi_{t_i} \ge \psi_{t_i-},\nonumber\\
&&\dis \mbox{and $\psi$ is a classical subsolution of master equation \reff{master} on each $[t_{i-1}, t_i)$}   \Big\}.\nonumber
\eea

Under mild conditions, for example when $g$ and $G(t,\mu, 0,0,0)$ are bounded, one can easily see that $\ol \cU$ and $\ul \cU$ are not empty.

\begin{prop}
\label{prop-olVvis}
Let Assumption \ref{assum-G} hold, $g\in C^0(\cP_2, \dbR)$, and $\ul \cU_g \neq \emptyset$. If $\ul V \in C^0(\Th)$, then $\ul V$   is a viscosity subsolution of master equation \reff{master}.
\end{prop}
\proof Fix $(t, \mu)\in\Th$.  Let $\f\in \ul \cA^L \ul V(t, \mu)$ with corresponding $\d>0$. For any $\e>0$, let   $\psi^\e\in \ul \cU_g$ be such that $\psi^\e(t, \mu) \ge \ul V(t,\mu) - \e$. It is clear that $\psi^\e(s, \dbP) \le \ul V(s,\dbP)$ for all $(s, \dbP)\in [t, T]\times \cP_L(t,\mu)$. Denote $\f_\d(s, \dbP) := \f(s, \dbP) + \d [s-t]$. For $\e < \d^2$ and any $\dbP\in \cP_L(t,\mu)$, we have
\beaa
[\f_\d - \psi^\e](t, \mu) &=& [\ul V- \psi^\e] (t, \mu) \le \e  < \d^2  =  [\f_\d - \f](t+\d, \dbP) \\
&\le& [\f_\d - \ul V](t+\d, \dbP) \le [\f_\d - \psi^\e](t+\d, \dbP).
\eeaa
Then there exists $(t^*, \dbP^*)\in [t, t+\d) \times \cP_L(t,\mu)$ such that
\beaa
[\f_\d - \psi^\e](t^*, \dbP^*) = c^* := \inf_{(s, \dbP)\in [t, t+\d] \times \cP_L(t,\mu)}[\f_\d - \psi^\e](s, \dbP).
\eeaa
This implies that $\f^\e_\d := \f_\d + c^* \in \ul \cA^L \psi^\e(t^*, \dbP^*)$. By Theorem \ref{thm-change}, we may assume without loss of generality that $G$ is increasing in $y$. Then by Theorem \ref{thm-consistency} we have
\beaa
0 \le \sL \f^\e_\d (t^*, \dbP^*) &=& \pa_t \f (t^*, \dbP^*) + \d + G\big(t^*, \dbP^*, \psi^\e(t^*, \dbP^*), \pa_\mu \f (t^*, \dbP^*,\cd), \pa_\o \pa_\mu \f (t^*, \dbP^*,\cd)\big)\\
&\le&\pa_t \f (t^*, \dbP^*) + \d + G\big(t^*, \dbP^*, \ul V(t^*, \dbP^*), \pa_\mu \f (t^*, \dbP^*,\cd), \pa_\o \pa_\mu \f (t^*, \dbP^*,\cd)\big).
\eeaa
Send $\d\to 0$, we  have $(t^*, \dbP^*) \to (t, \mu)$. Then the above inequality implies $\sL \f(t,\mu)\ge 0$. 
\qed

\begin{thm}
\label{thm-comparison}
Let Assumption \ref{assum-G} hold and $g\in C^0(\cP_2; \dbR)$. Assume $V_1$ and $V_2$ are viscosity subsolution and viscosity supersolution of master equation \reff{master} with $V_1(T,\cd) \le g \le V_2(T,\cd)$. Assume further that $\ul \cU_g $ and $\ol \cU_g$ are not empty and 
\bea
\label{olV=ulV}
\ol V = \ul V =: V.
\eea
Then $V_1 \le  V  \le V_2$ and $V$  is the  unique   viscosity solution  of master equation \reff{master}.
\end{thm}
\proof First one can easily show that $\ol V$ is lower semicontinuous and $\ul V$ is upper semicontinuous. Then by \reff{olV=ulV} $V$ is continuous, and thus it follows from Proposition \ref{prop-olVvis} that $V$ is a  viscosity solution  of master equation \reff{master}. 

To see the comparison principle, which implies immediately the uniqueness, we fix an arbitrary $\psi \in \ol\cU_g$.  First notice that $V_1(T,\cd) \le g \le \psi(T,\cd)$. Since $V_1$ is continuous and $\psi(T,\cd) \le \psi(T-,\cd)$, we have $V_1(T-,\cd) = V_1(T,\cd) \le  \psi(T,\cd) \le \psi(T-, \cd)$. Now apply the partial comparison principle  Theorem \ref{thm-partial}, one can easily see that  $V_1(t, \cd) \le \psi(t,\cd)$ for $t\in [t_{n-1}, t_n)$. Repeat the arguments backwardly in time we can prove $V_1 \le \psi$ on $\Th$. Since $\psi \in \ol\cU_g$ is arbitrary, we have $V_1 \le \ol V$. Similarly, one can show that $V_2 \ge \ul V$. Then it follows from  \reff{olV=ulV}  that $V_1 \le V \le V_2$.  
\qed 

The following result is a direct consequence of the above theorem.

\begin{thm}
\label{thm-approx}
Let Assumption \ref{assum-G} hold and $g\in C^0(\cP_2; \dbR)$. Assume there exist $(\ol G^n, \ol g^n)$ and $(\ul G^n,  \ul g^n)$  such that, for each $n$,

(i)  $\ol G^n, \ul G^n$ satisfy Assumption \ref{assum-G} and  $\ol g^n,  \ul g^n\in C^0(\cP_2; \dbR)$;

(ii) the master equation \reff{master} with generator $ \ol G^n$ (resp. $\ul G^n$) and terminal condition $\ol g^n$ (resp. $\ul g^n$) has a classical solution $\ol V^n$ (resp. $\ul V^n$);

(iii)  $\ol G_n \le G \le \ol G_n$, $\ul g_n \le g\le \ol g_n$;

(iv) $\lim_{n\to \infty} \ol V^n =  \lim_{n\to \infty} \ul V^n =: V$.

\no Then comparison principle holds for master equation \reff{master} with generator $G$ and terminal condition $g$, and $V$ is its unique viscosity solution.  
\end{thm}
\proof Clearly $\ol V^n$ is a classical supersolution of master equation \reff{master} with generator $G$ and terminal condition $g$, and it satisfies $\ol V^n \ge g$. Then $\ol V^n \ge \ol V$. Similarly $\ul V^n \le \ul V$. Then (iv) implies \reff{olV=ulV} and thus the statements follow from Theorem \ref{thm-comparison}. 
\qed

\subsection{\color{black}Some examples}
In this subsection we provide two examples for which we have the complete result for the comparison principle. While only for these special cases, the results are new in the literature, to our best knowledge. The comparison principle for more general master equations, especially the verification of condition \reff{olV=ulV}, is very challenging and we shall leave it for future research. 
\begin{eg}
\label{eg-distortion2}
Consider the setting in Example \ref{eg-distortion}, but relax the regularity of $\k$ to be only continuous. Then the $V$ defined by \reff{distortion-V} is in 
$ C^0([0, T] \times \cP_2(\dbR^d))$ and is the unique viscosity solution of the master equation \reff{distortionMaster}.
\end{eg} 
\proof (i) One can easily verify that $V$ is continuous and the DPP \reff{meancontrol-DPP} becomes:
\bea
\label{distortionDPP}
V(t, \mu_t) = V\big(t+\d, (\dbP^{t,\mu}_0)_{t+\d} \big),\q (t,\mu)\in \Th.
\eea
Denote $\nu := \mu_t$. Now let $L\ge 1$ and  $\f\in \underline{\cA}^L V(t,\nu)$. Clearly $\dbP^{t,\mu}_0 \in \cP_L(t,\nu)$. Then
\beaa
\f(t,\nu) =  V(t,\nu)  = V\big(t+\d, (\dbP^{t,\mu}_0)_{t+\d}\big) \le \f\big(t+\d, (\dbP^{t,\mu}_0)_{t+\d}\big).
\eeaa
Apply the It\^{o} formula, this implies
\beaa
0 &\le& \int_t^{t+\d} \pa_t  \f\big(t+\d, (\dbP^{t,\mu}_0)_s\big) ds \\
&&+ \dbE^{\dbP^{t,\mu}_0}\Big[\int_t^{t+\d} \pa_\mu \f(s,  (\dbP^{t,\mu}_0)_s, X_s) d X_s + {1\over 2} \int_t^{t+\d} \pa_x \pa_\mu \f(s,  (\dbP^{t,\mu}_0)_s, X_s) d \la X\ra_s\Big]\\
&=& \int_t^{t+\d} \Big[\pa_t  \f\big(t+\d,  (\dbP^{t,\mu}_0)_s\big) + {1\over 2} \dbE^{\dbP^{t,\mu}_0}\big[ \pa_x \pa_\mu \f(s,  (\dbP^{t,\mu}_0)_s, X_s) \big]\Big]ds.
\eeaa
Divide both sides by $\d$ and send $\d\to 0$, we obtain
\beaa
\pa_t  \f(t, \nu) + {1\over 2} \dbE^\nu\big[ \pa_x \pa_\mu \f(t, \nu, X_t) \big] \ge 0.
\eeaa
That is, $V$ is a viscosity subsolution at $(t,\nu)$. Similarly one can show that $V$ is a viscosity supersolution at $(t,\nu)$, hence a viscosity solution.

(ii) We next prove the comparison principle, which implies the uniqueness.  Assume $|g|\le C_0$. Then \reff{distortion-V} can be rewritten as:
\beaa
V(t,\mu) := \int_0^{C_0} \k\Big( \dbP^{t,\mu}_0(g(X_T) \ge y)\Big) dy.
\eeaa
Since $\k$ is continuous on $[0, 1]$, it is uniformly continuous, then there exists a smooth molifier $\k_n$ such that $\k_n$ is strictly increasing and $|\k_n - \k|\le {1\over n}$. Denote $\ol \k_n := \k_n + {1\over n}$, $\ul \k_n := \k_n - {1\over n}$, and define
\beaa
\ol V_n (t,\mu) := \int_0^{C_0} \ol\k_n\Big( \dbP^{t,\mu}_0(g(X_T) \ge y)\Big) dy,\q \ul V_n (t,\mu) := \int_0^{C_0} \ul\k_n\Big( \dbP^{t,\mu}_0(g(X_T) \ge y)\Big) dy
\eeaa
We remark that $\ol \k_n$ and $\ul \k_n$ does not satisfy the boundary conditions: $\k(0) =0, \k(1) = 1$. Nevertheless, following the same arguments in Example \ref{eg-distortion}, one can easily see that $\ol V_n$ and $\ul V_n$ are classical solutions of master type heat equation \reff{distortionMaster}, with terminal conditions 
\beaa
\ol V_n(T,\mu) := \int_0^{C_0} \ol\k_n\Big( \mu(g(X_T) \ge y)\Big) dy,\q \ul V_n(T,\mu) := \int_0^{C_0} \ul\k_n\Big( \mu(g(X_T) \ge y)\Big) dy,
\eeaa
respectively. It is clear that $\ul V_n \le V \le \ol V_n$ and $\lim_{n\to\infty} \ol V_n = \lim_{n\to\infty} \ul V_n = V$. Then the result follows from Theorem \ref{thm-approx} immediately.
\qed

The next example considers the following nonlinear (state dependent) master equation, which can be viewed as a special case of \reff{deterministic-master1} (see \cite{SZ}):
 \bea
\label{nonlinearMaster}
\pa_t V (t,\mu) \!+ \! {1\over 2} \dbE^\mu\big[\pa_x \pa_\mu V(t, \mu, X_t)\big]  \!+ \!  G_1\big(\dbE^\mu\big[\pa_\mu V(t, \mu, X_t)\big]\big)\!=\! 0,  V(T,\mu)\! = \! \dbE^\mu\big[g(X_T)\big].
\eea

\begin{eg}
\label{eg-semilinear2} Assume

(i) $g$ is Lipschitz continuous with Lipschitz constant $L_0$, and $G_1\in C^0\big([-L_0, L_0]\big)$;

(ii) Either $g$ is convex and $G_1$ is concave, or $g$ is concave and $G_1$ is convex;

\no Then the master equation  \reff{nonlinearMaster}  has a unique viscosity solution $V\in C^0([0, T]\times \cP_2(\dbR))$. 
\end{eg} 
\proof Let $G^n_1$ and $g_n$ be smooth mollifiers of $G_1$ and $g$, respectively, such that $|G_1^n - G_1|\le {1\over n}, |g_n-g|\le {1\over n}$. Denote $\ol G_1^n := G_1^n + {1\over n}$, $\ul G_1^n := G_1^n - {1\over n}$, $\ol g_n := g_n + {1\over n}$, $\ul g_n := g_n - {1\over n}$. Then  $(\ol G_1^n, \ol g_n)$ and $(\ul G_1^n, \ul g_n)$ are smooth and still satisfy (i) and (ii)  with the same $L_0$. By Saporito \& Zhang \cite[Theorem 3.1]{SZ} the corresponding master equations \reff{nonlinearMaster} have a classical solution $\ol V_n$ and $\ul V_n$, respectively. 

Now by Theorem \ref{thm-approx}, it suffices to show that $\ol V_n$ and $\ul V_n$ converge to the same limit. Without loss of generality, we assume $G_1$ is convex (and $g$ is concave). Denote 
\beaa
b(a) := \sup_{y\in [-L_0, L_0]} [ ay - G_1(y)],\q b_n(a) := \sup_{y\in [-L_0, L_0]} [ ay - G_1^n(y)],\q a\in \dbR.
\eeaa
By \cite{SZ} (or following similar arguments as in Section \ref{sect-HJB} below), we have
\beaa
\ol V_n(t, \mu) = \sup_{a\in \dbR} \dbE^{\dbP^{t,\mu}_0}\Big[g\big(X_T + [b_n(a) - {1\over n}][T-t] \big)\Big],\\
\ul V_n(t, \mu) = \sup_{a\in \dbR} \dbE^{\dbP^{t,\mu}_0}\Big[g\big(X_T+ [b_n(a) + {1\over n}][T-t]\big)\Big].
\eeaa
It is clear that $|b_n-b|\le {1\over n}$. Then it is straightforward to show that 
\beaa
\lim_{n\to\infty} \ol V_n(t, \mu) = \lim_{n\to\infty} \ol V_n(t, \mu) = V(t,\mu) := \sup_{a\in \dbR} \dbE^{ \dbP^{t,\mu}_0}\Big[g\big(X_T + b(a) [T-t] \big)\Big].
\eeaa
Now the result follows directly from Theorem \ref{thm-approx}.
\qed

\section{McKean-Vlasov SDEs with closed-loop controls}
\label{sect-HJB}
\setcounter{equation}{0} 
{\color{black}In this section we apply our viscosity theory to the mean field control problem introduced in Subsection \ref{sect-meanfieldcontrol}.  Recall \reff{Xa1}, \reff{control-Pa}, and \reff{control-Vtmu}, we shall assume

\begin{assum}
\label{assum-standing} $b, \si, f$ are $\dbF$-progressively measurable in all variables $(t,\o,\mu,a)\in [0, T]\times \O\times \cP_2\times A$ (and in particular $\dbF$-adapted in both $\o$ and $\mu$), and $g$ is progressively measurable in $(\o,\mu) \in \O\times \cP_2$. Moreover,

(i) $b, \si$ are bounded by a constant $C_0$, continuous in $a$, and uniform Lipschitz continuous in $(\o, \mu)$ with a Lipschitz constant $L_0$:
\beaa
|(b,\si)(t,\o, \mu,a)-(b,\si)(t,\o', \mu',  a)|\le L_0 \big[\|\o_{t\wedge \cd} - \o'_{t\wedge \cd}\|+ \cW_2(\mu_{[0, t]}, \mu'_{[0, t]})\big];
\eeaa

(ii) $f(t, 0, \d_{\{0\}}, a)$ is bounded by a constant $C_0$, $f$ is continuous in $a$, and $f$ and $g$ are uniformly continuous in $(\o, \mu)$ with a modulus of continuity function $\rho_0$:
\beaa
&|f(t,\o, \mu,a)-f(t,\o', \mu',  a)|\le \rho_0 \big(\|\o_{t\wedge \cd} - \o'_{t\wedge \cd}\|+ \cW_2((t,\mu),(t,\mu'))\big);\\
&|g(\o, \mu)-g(\o', \mu')|\le \rho_0 \big(\|\o - \o'\|+ \cW_2((T,\mu),(T,\mu'))\big).
\eeaa

(iii) $ \f= b, \si, f$ is locally uniformly continuous in $t$ in the following sense:
\beaa
  |\f(s, \o_{t\wedge \cd}, \mu_{[0, t]}, a) - \f(t, \o, \mu, a)| \le C\Big[ 1+ \|\o_{t\wedge \cd}\| +  \cW_2(\mu_{[0, t]}, \d_{\{0\}})\Big]\rho_0(s-t),~ t<s.
\eeaa

(iv) $\si \si^\top$ is  positive definite.
\end{assum} 
We remark that one sufficient condition for (iii) is that $A$ is compact, and the nondegeneracy of $\si$ in (iv) is used in Lemma \ref{lem-hatVreg} below, but we do not need uniform  nondegeneracy.

The choice of the admissible controls is very subtle, with \reff{control-cA} as one example. We shall discuss alternative choices in details at below. One basic requirement is that  the corresponding value function should satisfy the DPP. 
}

\subsection{Open-loop controls}
In this subsection, we consider open-loop controls, namely $\a_t = \a_t(B_\cd)$ depending on $B$, where $B$ is a Brownian motion in a probability space  $(\O, \cF, \dbP_0)$. There are two natural choices: (i) $\cA_t^{1}$,  where $\a_s=\a(s,(B_{t,r})_{t\le r\le s})$ is adapted to the shifted filtration of $B$; and (ii) $\cA_t^{2}$, where  $\a_s=\a(s,(B_r)_{0\le r\le s})$ is adapted to the full filtration of $B$. For the standard control problems, they would induce the same value function. However, in our setting the issue is quite subtle.  To be precise, for an {\color{black} $\dbF^B$-progressively measurable process $\xi$} on $[0, t]$ and a control $\a$, let $X^{t,\xi,\a}_s := \xi_s$, $s\in [0, t]$, and 
\beaa
X^{t,\xi,\a}_s = \xi_t + \int_t^s b(r, X^{t,\xi,\a}, \cL_{X^{t,\xi, \a}}, \a_r) dr +  \int_t^s \si(r, X^{t,\xi,\a}, \cL_{X^{t,\xi, \a}}, \a_r) dB_r,\q \dbP_0\mbox{-a.s.}
\eeaa
which has a unique strong solution under Assumption \ref{assum-standing}. Introduce the values functions:
\bea
\label{openVi}
V_i(t,\xi) := \sup_{\a\in\cA^{i}_t}  \dbE^{\dbP_0} \Big[g(X^{t,\xi,\a}, \cL_{X^{t,\xi, \a}}) + \int_t^T f(s, X^{t,\xi,\a}, \cL_{X^{t,\xi, \a}}, \a_s)ds\Big], \q i=1,2. 
\eea
The following example shows that $\cA^1_t$ is  not a good choice. 
\begin{eg}
\label{eg-neq}
Let $d=1$, $A=[-1,1]$, $b(t,\o, \mu, a) = a$, $\si\equiv 1$, $f\equiv 0$, and $g(\o, \mu) = g(\mu) = -Var^\mu (X_T)$. 

(i) $V_1(t,\xi) < V_2(t,\xi)$ when $\xi_t = (T-t) \mbox{sign}(B_t)$ and $T-t > 1$.
 
 (ii) $V_1$ does not satisfy the DPP:  $ V_1(t_1,\xi) \neq \sup_{\a\in \cA^1_{t_1}} V_1(t_2, X^{t_1, \xi, \a})$.
 \end{eg}
\proof (i) For any $(t,\xi)$ and  $\a\in \cA^{1}_t$, notice that $\xi_t$ is independent of $\a$ and thus is also independent of $X^{t,\xi,\a}_T - \xi_t$. Then
\bea
\label{cAo1}
Var(X^{t,\xi,\a}_T) = Var(\xi_t) +  Var(X^{t,\xi,\a}_T - \xi_t),~ \mbox{ thus}~ V_1(t, \xi) \le - Var(\xi_t)= - [T-t]^2.
\eea
On the other hand, set $\a_s := -\mbox{sign}(B_t)$, $s\in [t, T]$. Then $\a\in \cA^{2}_t$, $X^{t,\xi,\a}_T = B_{t,T}$, and thus
\beaa
 V_2(t,\xi) \ge Var(B_{t, T})=-[T-t] >  - [T-t]^2 =  V_1(t, \xi). 
 \eeaa

 (ii) Denote 
 $
\dis h(t) := \sup_{\a\in \cA^{1}_0} \big[-Var\big( \int_0^t \a_s ds + B_t\big)\big].
 $
 Then by \reff{cAo1}  one can easily see that 
 \beaa
V_1(t,\xi)=  h(T-t) -Var(\xi_t).
\eeaa
 Assume by contradiction that DPP holds. Then, for any $0<t<T$, 
 \beaa
h(T) &=& V_1(0, \d_{\{0\}}) =  \sup_{\a\in \cA^{1}_0} V_1(t, X^{0,0,\a}_t) \\
&=& \sup_{\a\in \cA^{1}_0} \big[ h(T-t)- Var(X^{0,0,\a}_t)\big] = h(t) + h(T-t).
 \eeaa
 Following the same arguments we see that $h$ is linear in $t$. Since $|\a|\le 1$, it is clear that 
{\color{black}
\beaa
Var\big( \int_0^t \a_s ds + B_t\big) &=& \dbE\Big[ \big(\int_0^t \a_s ds + B_t - \dbE[\int_0^t \a_s ds]\big)^2\Big] \\
&=& \dbE\Big[ \big(B_t +O(t)]\big)^2\Big] = \dbE[|B_t|^2] + o(t) = t + o(t).
\eeaa
}
Then
 \beaa
 \lim_{t\to 0} {h(t) \over t} =-1,\q\mbox{and thus}\q h(t) =- t.
 \eeaa
 On the other hand,  fix $t\in (0, T)$  and set $\a_s := \Big[ (-1) \vee (-{B_{t}\over  T-t}) \wedge 1\Big]1_{[t, T]}(s)$. Then
 \beaa
 \int_0^T \a_s ds + B_T =  (t-T) \vee (- B_{t}) \wedge (T-t) + B_{t} + B_{t,T}. 
 \eeaa
 Thus
 \beaa
- h(T) &\le& Var\Big(\int_0^T \a_s ds + B_T\Big) = Var\Big( (t-T) \vee (- B_{t}) \wedge (T-t) + B_{t}  \Big) + T-t\\
 &=& \dbE^{\dbP_0}\Big[ \big([|B_t| - [T-t]]^+\big)^2\Big] +T-t <  \dbE^{\dbP_0}\Big[ |B_t|^2\Big] + T-t  = t+T-t= T.
 \eeaa
 This is a desired contradiction.
 \qed

Technically, the choice of $\cA^{2}_t$ would work, see e.g. Bayraktar, Cosso, \&  Pham \cite{BCP}. The following results can be proved easily, in particular, the viscosity property in (iii) follows similar arguments as in Theorem \ref{thm-controlviscosity} below, and thus we omit the proofs. 

\begin{prop}
\label{prop-strong}
Let Assumption \ref{assum-standing} hold and define $V_2(t,\xi)$ by \reff{openVi}. Then

(i) $V_2$ satisfies the following DPP:
{\color{black}\beaa
V_2(t_1, \xi) = \sup_{\a\in\cA^2_{t_1}} \Big[V_2(t_2, X^{t_1, \xi, \a}) + \int_{t_1}^{t_2} \dbE^{\dbP_0}[f(s, X^{t,\xi,\a}, \cL_{X^{t,\xi,\a}}, \a_s)] ds\Big],\q t_1 < t_2.
\eeaa
}

(ii) $V_2(t,\xi)$ is law invariant and thus we may define $V_2'(t,\mu)$  by $V_2(t, \xi) =  V_2'(t, \cL_\xi)$. 

(iii) $ V_2' \in C^0(\Th)$  and  is a viscosity solution of the HJB type of master equation \reff{meancontrol-master}.
\end{prop}

Despite the above nice properties, in many applications the state process $X$ is observable while the Brownian motion $B$ is used to model the distribution of $X$ and may not  be observable. Then it is not reasonable to have the controls relying on $B$. The issue becomes more serious when one considers games instead of control problems. We refer to Zhang \cite{Zhang} Section 9.1 for detailed discussions on these issues. Therefore, in the next subsection we shall turn to closed-loop controls. 

\subsection{Closed-loop controls}
We now assume $\a$ depends on the state process $X^{t,\xi,\a}$. One choice is to use the (state dependent) feedback controls:  $\a_s=\a(s,X^{t,\xi,\a}_s)$, see, e.g. Pham \& Wei \cite{PW}. However, we prefer not to use this for several reasons:

$\bullet$ In practice it is not natural to assume the players cannot use past information;

$\bullet$  It seems difficult to have regularity of $V(t,\mu)$ without strong constraint on $\a$;

$\bullet$ It fails to work in non-Markovian models, which are important in applications.

\no We shall assume $\a$ is $\dbF^{X^{t,\xi,\a}}$-measurable, namely $\a_s = \a(s, (X^{t,\xi,\a}_r)_{0\le r\le s})$, and thus we are considering \reff{Xa1}. As mentioned in Subsection \ref{sect-meanfieldcontrol}, in this case it is more convenient to use weak formulation. That is, we shall use the canonical setting in Subsection \ref{sect-canonical}, and consider the optimization problem   \reff{control-Pa} and \reff{control-Vtmu}. However, under closed-loop controls, the regularity of $V(t,\mu)$ is rather technical. In this section we content ourselves with the following piecewise constant control process:
\bea
\label{cAt}
\left.\ba{c}
\dis\cA_t:=\Big\{\a:  \exists ~\, n ~\mbox{and } t=t_0<\cdots < t_n=T~ \mbox{ such that }
\a_s=\sum_{i=0}^{n-1}h_i1_{[t_i,t_{i+1})}(s),\\ 
\dis \mbox{where } h_i: \O \to A \mbox{ is $\cF_{t_i}$-measurable for}~ i=0,\cds, n-1\Big\}.
\ea\right.
\eea
{\color{black}We emphasize that here we are abusing the notation $\cA$ with \reff{control-cA}. So throughout this section, our optimization problem will always be \reff{control-Pa}-\reff{cAt}-\reff{control-Vtmu}.
}

\begin{rem}
\label{rem-closed}
{\rm (i) Each $\a\in \cA_t$  here also satisfies the requirement in \reff{control-cA}, and thus \reff{control-Pa} has a unique (strong) solution $\dbP^{t,\mu,\a}\in \cP_L(t,\mu)$, where $L\ge C_0 \vee [ {1\over 2} C_0^2]$ for the bound $C_0$ in Assumption \ref{assum-standing} (i). {\color{black} In particular, $\dbP^{t,\mu,\a}$ satisfies  the uniform estimate \reff{pLest}.

(ii) Obviously the $\cA_t$ in \reff{cAt} also satisfies \reff{compatibility}. Then following the same arguments as in Theorem \ref{thm-meancontrol} (i) we see that, under Assumption \ref{assum-standing},  $V$ satisfies the  DPP \reff{meancontrol-DPP}. }
\qed}
\end{rem}

\begin{rem}
\label{rem-strong}
{\rm (i)  Although \reff{Xa1} (and \reff{control-Pa}) has a strong solution, the formulation \reff{control-Vtmu}  is still different from the $V_2(t,\xi)$ in \reff{openVi}. Indeed, by the piecewise constant structure, one can easily see that $\dbF^{X^{t,\xi,\a}}$ is the same as the filtration generated by the process $\tilde B_s:= \xi_s 1_{[0, t]}(s) + [\xi_t + B_{t,s}]1_{(t, T]}(s)$, and thus one may rewrite $\a(s, X^{t,\xi,\a}_{[0, s]})$ as $\tilde \a(s, \tilde B_{[0,s]})$ for some measurable function $\tilde \a$. However, note that $\tilde B_{[0,t]} = \xi_{[0,t]} \neq B_{[0,t]}$, so this control is still not in $\cA^{2}_t$.  Indeed, in many practical situations, at time $t$, one can observe the state process $\xi_{\cd\wedge t}$,  but not necessarily observe an underlying Brownian motion path in the past. That is the main reason we consider the closed-loop controls in this paper.

(ii) The regularity of $V_2$ and $V_2'$ in Proposition \ref{prop-strong} (iii) is straightforward. However, the above subtle difference makes the regularity of $V$ in \reff{control-Vtmu} quite involved, as we will see in Example \ref{eg-discontinuous} and Subsection \ref{sect-reg} below.
\qed}
\end{rem}

\begin{eg}
\label{eg-discontinuous}
Let $d=1$, $A=[-1,1]$, $T=1$, $b\equiv 0$, $\sigma(t,\mu,a) = 1+a^2$, $f\equiv 0$,  $g(\o,\mu) = g(\mu)=\frac{1}{3}\dbE^\mu[X_1^4]-(\dbE^\mu[X^2_1])^2$, and $\cA^0_0$ consist of constant controls: $\a_t\equiv\a_0(X_0), \forall t\in[0,1]$. Then $V^0_0(\mu) := \sup_{\a\in \cA^0_0} g(\dbP^{0,\mu, \a})$ is discontinuous in $\mu\in \cP_2$.
\end{eg}
\proof Let $\mu_0 := \d_{\{0\}}$ and $\mu_\e := {1\over 2} [\d_{\{\e\}} + \d_{\{-\e\}}]$. It is clear that $\lim_{\e\to 0}\cW_2(\mu_\e, \mu_0) = 0$.
For any $\a\in \cA^0_0$, we have $\a_t = \a_0(0)$ and $X_1 = [1+|\a_0(0)|^2] B^\a_1$, $\dbP^{0,\mu, \a}$-a.s. Then, denoting $c := 1+|\a_0(0)|^2$ and $\dbP^\a := \dbP^{0,\mu, \a}$, we have
\beaa
g(\dbP^{0,\mu, \a}) = {1\over 3}\dbE^{\dbP^\a}[c^4|B^\a_1|^4]-(\dbE^{\dbP^\a}[c^2|B^\a_1|^2])^2=0,\q\mbox{and thus}\q V^0_0(\mu_0) = 0.
\eeaa
On the other hand, for each $\e>0$, set $\a_t := \a_0(X_0) := \1_{\{X_0>0\}}$. Then
\beaa
X_1 = \big[\e + 2 B^\a_1\big] 1_{\{X_0=\e\}} + \big[-\e +  B^\a_1\big] 1_{\{X_0=-\e\}},\q \dbP^{0, \mu_\e, \a}\mbox{-a.s.}
\eeaa
Thus, denoting $\dbP^{\e}:= \dbP^{0, \mu_\e, \a}$,
\beaa
 g( \dbP^{0, \mu_\e, \a}) &=& \frac{1}{6}\dbE^{\dbP^{\e}}\big[(2B^\a_1+\eps)^4+(B^\a_1-\eps)^4\big]-\Big(\frac{1}{2}\dbE^{\dbP^{\e}}[(2B^\a_1+\eps)^2+(B^\a_1-\eps)^2]\Big)^2\\
 &=& {1\over 6}[51 + 18 \e^2 + 2\e^4] - \Big({1\over 2}[5+2\e^2] \Big)^2 = \frac{9}{4}-2\e^2 -\frac{2}{3}\eps^4.
\eeaa
Therefore, for all $\e>0$ small,
\beaa
V^0_0(\mu_\e) \ge  \frac{9}{4}-2\e^2 -\frac{2}{3}\eps^4 \ge 2 > 0 = V^0_0(\mu_0).
\eeaa
 This implies that $V^0_0$ is discontinuous at $\mu_0$.
\qed

Nevertheless, by using piecewise constant controls $\cA_t$, we have 
\begin{thm}
\label{thm-regmu}
Under Assumption \ref{assum-standing}, there exists a modulus of continuity function $\rho$ such that
\bea
\label{Vreg}
|V(t_1, \mu) - V(t_2, \nu)| \le C\rho\big(\cW_2((t_1, \mu), (t_2, \nu) \big) + C[1+\cW_2(\mu_{[0, t_1]}, \d_{\{0\}})] [t_2-t_1].
\eea
Assume further that $f$ is bounded, then $V$ is uniformly continuous in $(t,\mu)$. 
\end{thm}
The proof of this theorem is quite involved, so we defer it to the next subsection.

Given the above regularity, we can easily verify the viscosity property.
\begin{thm}
\label{thm-controlviscosity}
Under Assumption \ref{assum-standing}, $V$  is a viscosity solution of the HJB type of master equation  \reff{meancontrol-master}.
\end{thm}
\proof  Fix $L>0$ such that $|b|, {1\over 2}|\sigma|^2\le L$.  We shall show that $V$ is an $L$-viscosity solution.

{\it Step 1.} We first verify its the viscosity subsolution property. Assume by contradiction that $V$ is not an $L$-viscosity subsolution at $(t,\mu)$, then there exists $(v, Z, \G) \in\underline{\cJ}^L V(t,\mu)$ with corresponding $\d$, such that
\bea
\label{cL1}
-c:= \sL \f(t,\mu) = v+ \dbE^\mu\big[\sup_{a\in A} G_2(t, \mu, X, Z, \G, a)\big] <0,\q\mbox{where}~ \f:= \phi^{t, V(t,\mu), v, Z, \G}.
\eea
 For any $\a\in \cA_t$, applying the functional It\^{o} formula  we have
\bea
\label{HJB-Ito1}
\f(t+\d, \dbP^{t,\mu,\a})  - \f(t,\mu) = \int_t^{t+\d} \sL^\a\f(s, \dbP^{t,\mu,\a})ds,
\eea
where, abbreviating $\dbP^\a := \dbP^{t,\mu,\a}$,
\bea
\label{cLa}
\sL^\a \f(s, \dbP^\a) = a +  \dbE^{\dbP^\a}\Big[ \big[b(s, X, \dbP^\a, \a_s) \cd [Z+ \G X_{t,s}]  + {1\over 2} \G: \si\si^\top (s, X, \dbP^\a, \a_s)   \big]\Big].
\eea
Note that
 \bea
  \label{HJB-I12}
  \sL^\a \f(s, \dbP^\a) -  \sL \f(t, \mu) = I_1(s) + I_2(s) - \dbE^{\dbP^\a} [f(s, X, \dbP^\a, \a_s)],
  \eea
  where
  \beaa
&&I_1(s) := \dbE^{\dbP^\a}\Big[ G_2(t, \mu, X, Z, \G, \a_s) - \sup_{a\in A}G_2(t, \mu, X, Z, \G, a) \Big];\\
&&I_2(s) := \dbE^{\dbP^\a}\Big[ Z\cd \big[b(s, X, \dbP^\a, \a_s)- b(t, X,  \mu, \a_s) \big] + b(s, X, \dbP^\a, \a_s) \cd \G X_{t,s}  \\
&&\q + {1\over 2} \G: \big[\si\si^\top (s, X, \dbP^\a, \a_s)  -\si\si^\top (t, X, \mu, \a_s) \big] + \big[f(s, X, \dbP^\a, \a_s) - f(t, X, \mu, \a_s)\big]\Big].
  \eeaa
  It is clear that $I_1(s) \le 0$.  By Assumption \ref{assum-standing} (ii) and (iii), we have, for $s\in [t, t+\d]$,
  \beaa
  &&\dbE^{\dbP^\a}\Big[ f(s, X, \dbP^\a, \a_s) - f(t, X, \mu, \a_s)\Big]\\
  &&=\dbE^{\dbP^\a}\Big[ \big[f(s, X, \dbP^\a, \a_s) - f(s, X_{t\wedge \cd}, \mu_{[0,t]}, \a_s)\big]+\big[ f(s, X_{t\wedge \cd}, \mu_{[0, t]}, \a_s) - f(t, X, \mu, \a_s)\big]\Big]\\
  &&\le C\dbE^{\dbP^\a}\Big[ \rho_0\big(\|X_{s\wedge \cd} - X_{t\wedge \cd}\| + \cW_2( \dbP^\a_{[0, s]}, \mu_{[0, t]}) \big)\Big] + C\Big(\dbE^\mu\big[ 1+ \|X_{t\wedge \cd}\|^2]\big]\Big)^{1\over 2}  \rho_0(\d).
\eeaa
  Since $b, \si$ are bounded, by \reff{pLest} one can easily see that
  \beaa
  \lim_{\d\to 0} \sup_{\a\in \cA_t} \dbE^{\dbP^\a}\Big[ \rho_0\big(\|X_{s\wedge \cd} - X_{t\wedge \cd}\| + \cW_2( \dbP^\a_{[0, s]}, \mu_{[0, t]}) \big)\Big] = 0.
  \eeaa
 Then, for some $\d = \d(t,\mu)>0$ small enough,   we have
 \beaa
 \dbE^{\dbP^\a}\Big[ f(s, X, \dbP^\a, \a_s) - f(t, X, \mu, \a_s)\Big] \le {c\over 8},
  \eeaa
  for all $s\in [t, t+\d]$ and all $\a\in \cA_t$. Similarly, recalling that by definition $Z, \G$ have linear growth in $\o$, we may have the desired estimates for the other terms in $I_2(s)$, and thus $I_2(s) \le {c\over 2}$. Therefore, \reff{HJB-I12} implies that
  \beaa
   \sL^\a \f(s, \dbP^\a) + \dbE^{\dbP^\a} [f(s, X, \dbP^\a, \a_s)] =  \sL \f(t, \mu) + I_1(s) + I_2(s) \le -{c\over 2}.
   \eeaa
    Plug this into \reff{HJB-Ito1} and recall \reff{cA}, we get
  \beaa
  &&V(t+\d, \dbP^{t,\mu,\a})  + \int_t^{t+\d} \dbE^{\dbP^\a} [f(s, X, \dbP^\a, \a_s)] ds - V(t,\mu) \\
  &&\le  \f(t+\d, \dbP^{t,\mu,\a})  - \f(t,\mu)  + \int_t^{t+\d} \dbE^{\dbP^\a} [f(s, X, \dbP^\a, \a_s)] ds \\
  && = \int_t^{t+\d} \Big[ \sL^\a \f(s, \dbP^\a) + \dbE^{\dbP^\a} [f(s, X, \dbP^\a, \a_s)] \Big] ds \le -{c\d \over 2},\q\forall \a\in \cA_t.
  \eeaa
    Take supremum over $\a\in\cA_t$, this contradicts with the DPP \reff{meancontrol-DPP}, see Remark \ref{rem-closed} (ii).

  {\it Step 2.} We next verify its  viscosity supersolution property. Assume by contradiction that $V$ is not an $L$-viscosity supersolution at $(t,\mu)$, then there exists $(v, Z, \G) \in\overline{\cJ}^L V(t,\mu)$ with corresponding $\d$, such that
\bea
\label{cL2}
c:= \sL \f(t,\mu) = v+ \dbE^\mu\big[\sup_{a\in A} G_2(t, \mu, X, Z, \G, a)\big] > 0,\q\mbox{where}~ \f:= \phi^{t, V(t,\mu), v, Z, \G}.
\eea
Note that $G_2(t, \mu, X, Z, \G, a)$ is $\cF_t$-measurable, there exists an $\cF_t$-measurable $A$-valued random variable $\a_t$ such that
\bea
\label{cLa2}
 v+ \dbE^\mu\big[G_2(t, \mu, X, Z, \G, \a_t)\big]  \ge {c\over 2}.
\eea
Now let $\a_s := \a_t$, $s\in [t, t+\d]$ and denote $\dbP := \dbP^{t,\mu, \a}$. Clearly $a\in \cA_t$. Applying the functional It\^{o} formula  we have
\bea
\label{HJB-Ito2}
\f(t+\d, \dbP)  - \f(t,\mu) = \int_t^{t+\d} \sL^\a\f(s, \dbP)ds.
\eea
where $\sL^\a$ is the same as \reff{cLa}.  Similar to the estimate of $I_2(s)$ in Step 1, for $\d>0$ small enough we have
\beaa
&& v+ \dbE^\mu\big[G_2(t, \mu, X, Z, \G, \a_t)\big] - \Big[\sL^\a\f(s, \dbP)+ \dbE^{\dbP} [f(s, X, \dbP, \a_t)] \Big] \le {c\over 4}.
\eeaa
Then, by \reff{cLa2},
\beaa
\sL^\a\f(s, \dbP)+ \dbE^{\dbP} [f(s, X, \dbP, \a_t) ]\ge {c\over 4},\q  s\in [t, t+\d].
\eeaa
This implies 
  \beaa
  &&V(t+\d, \dbP) +  \int_t^{t+\d}\dbE^{\dbP} [f(s, X, \dbP, \a_t) ] ds  - V(t,\mu) \\
  &&\ge  \f(t+\d, \dbP) +  \int_t^{t+\d}\dbE^{\dbP} [f(s, X, \dbP, \a_t) ] ds - \f(t,\mu)  \\
  &&=\int_t^{t+\d} \Big[\sL^\a\f(s, \dbP)+ \dbE^{\dbP} [f(s, X, \dbP, \a_t) ]\Big]ds  \ge {c\d \over 4}.
  \eeaa
 Again this contradicts with the DPP \reff{meancontrol-DPP}.
\qed

We remark again that the comparison principle for  HJB master equation \reff{meancontrol-master}  is quite challenging and we shall leave it for future research.

\begin{rem}
\label{rem-open-close}
{\rm Under nice conditions, in particular when the comparison principle for the master equation \reff{meancontrol-master} holds, by Proposition \ref{prop-strong} and Theorem \ref{thm-controlviscosity} we see that $V = V'_2$, for the $V_2'$ defined by \reff{openVi} and Proposition \ref{prop-strong} (ii). This is well known for standard control problems, and is also known in state dependent McKean-Vlasov setting, see Lacker \cite{Lacker}.  

 However, for zero-sum games, the open-loop controls and closed-loop controls are quite different, see e.g. Pham \& Zhang \cite{PZ}, Sirbu \cite{Sirbu},  and  Possamai, Touzi, \& Zhang \cite{PTZ} in the standard setting.  While in this paper we consider only the control problem, we expect our arguments will work for zero-sum game problems with closed-loop controls in McKean-Vlasov setting. We note that such game problem is studied in recent work Cosso \&  Pham \cite{CP} by using strategy versus open-loop controls.
\qed}
\end{rem}

\begin{rem}
\label{rem-optimal}
{\rm The restriction to piecewise constant controls makes it essentially impossible to obtain optimal controls. As we understand such restriction is mainly for the regularity of the value function $V$. In Possamai, Touzi, \& Zhang \cite{PTZ}, we studied the zero sum games under general closed-loop controls (but without involving the measures) and faced similar regularity issues. However, in \cite{PTZ} we obtained the desired regularity  when $b$ and $\si$ do not depend on the path and then proved the verification theorem for optimal controls. It will be interesting to remove the piecewise constant constraint in this framework when $b$ and $\si$ do not depend on $\mu$.
\qed}
\end{rem}

\subsection{Regularity of $V$}
\label{sect-reg}
In this subsection we prove Theorem \ref{thm-regmu}. To simplify the notation, in this subsection we assume $d=1$. But the proof can be easily extended to the multidimensional case. Introduce
\bea
\label{hatcAt}
&\dis V_0(t,\mu):=\sup_{\a\in\cA^0_t} J(t,\mu,\a),\q\mbox{where $J$ is defined in \reff{control-Vtmu} and }&\\
&\dis\cA^0_t:=\Big\{\a = \sum_{i=0}^{n-1} h_i 1_{[t_i, t_{i+1})} \in \cA_t: ~\mbox{there exist $0\le s_1<\cds< s_{m}\le t$ such that }\nonumber\\
&\dis h_i = h_i(X_{s_1},\cdots,X_{s_m},X_{[t,t_i]})~\mbox{for}~i=0,\cds,n-1\Big\}.&\nonumber
\eea
That is, $h_i$ depends on $X_{[0, t]}$ only discretely.  Since $\cA^0_t\subset\cA_t$, clearly $V_0(t,\mu)\le V(t,\mu)$. We will actually prove $V_0 = V$, then it suffices to establish the regularity of $V_0$. 

To see the idea, let's first observe the following simple fact. Given an arbitrary probability space $(\tilde \O, \tilde \cF, \tilde \dbP)$ and a random variable $\zeta$ with continuous distribution, then for any other random variable $\tilde \zeta$, there exists a deterministic function $\f$ such that 
\bea
\label{fzeta}
\cL_{\f(\zeta)} = \cL_{\tilde \zeta},
\eea
where $\cL$ denotes the distribution under $\tilde \dbP$. 
  Indeed, denoting by $F$ the cumulative distribution function, then $\f := F_{\tilde \zeta}^{-1} \circ F_{\zeta}$ serves our purpose. In Example \ref{eg-discontinuous}, assume $\cL_{\zeta} = \mu_0$ and $\cL_{\tilde \zeta} = \mu_\e$. The discontinuity of $V^0_0$ at $\mu_0$ is exactly because there is no function $\f$ such that \reff{fzeta} holds. The next lemma is crucial for overcoming such difficulty. Recall the $\cP(\mu, \nu)$ and the product space $(\O\times \O, \cF\times \cF)$ in \reff{cW21}, and denote the canonical process as $ (X, X')$. Moreover, for a partition $\pi$: $0\le s_1<\cds<s_m\le t$,  $\mu \in \cP_2$, and two processes $\xi, \eta$ on a probability space $(\tilde \O, \tilde \cF, \tilde \dbP)$, we introduce the notations:
  \bea
  \label{pi}
  \left.\ba{c}
\dis  \mu_\pi := \mu \circ (X_{s_1},\cds, X_{s_m})^{-1},\q \xi_\pi := (\xi_{s_1},\cds, \xi_{s_m}),\\
\dis \|\xi-\eta\|_{\tilde \dbP,\pi} := \|\xi_\pi-\eta_\pi\|_{\tilde \dbP} :=  \Big(\dbE^{\tilde \dbP}\big[\max_{1\le j\le n} |\xi_{s_j}-\eta_{s_j}|^2\big]\Big)^{1\over 2}.
\ea\right.
  \eea

\begin{lem}
\label{lem-construct}
Let $0< t< T$, $\mu, \nu \in\cP_2$, $\overline \dbP\in\cP(\mu,\nu)$. Then for any $\e>0$, $\d>0$, and any partition  $\pi: 0\le s_1<\cdots<s_m\le t$, there exist a probability space $(\tilde \O, \tilde \cF, \tilde \dbP)$, two continuous processes $(\xi, \eta)$, and a Brownian motion $\tilde{B}$ on $(\tilde \O, \tilde \cF, \tilde \dbP)$ such that: 
\begin{enumerate}
\item[(i)] $\cL_\xi=\mu, \cL_\eta = \nu$, and $\eta$ is independent of $\tilde B$;
\item[(ii)] $\xi_\pi$ is measurable to the $\si$-algebra $\sigma(\eta_\pi,\tilde{B}_{[0,\d]})$.
\item[(iii)] $\|\xi-\eta\|_{\tilde\dbP,\pi}\le \|X - X'\|_{\overline \dbP, \pi} +\eps$.
\end{enumerate}
\end{lem}
\proof We prove the lemma in several cases, depending on the joint distribution $\nu_\pi$.  Fix an arbitrary process $\eta$ with $\cL_\eta = \nu$. Note that we shall extend the space $(\tilde \O, \tilde \cF, \tilde \dbP)$  whenever needed, and we still denote this process as $\eta$.

\textit{\underline{\smash{Case 1}}:} $\nu_\pi$ is degenerate, namely  $\nu_\pi = \d_{\{(x_1,\cds, x_m)\}}$ for some $(x_1,\cds, x_m)\in\dbR^m$, and thus  $\eta_{s_j} = x_j$, $\tilde \dbP$-a.s.  Pick a Brownian motion $\{\tilde{B}_s\}_{s\in[0,\d]}$ independent of $\eta$ (which is always doable by extending the probability space if necessary). In the spirit of \reff{fzeta},  one can easily construct a $m$-dimensional random vector  $\tilde \xi_\pi =(\tilde \xi_1, \cds, \tilde \xi_m)$ such that  $\cL_{\tilde \xi_\pi} = \mu_\pi$ and $\tilde\xi_\pi$ is measurable to the $\si$-algebra $\si\big(\tilde B_{{(j-1)\d\over m}, {j\d\over m}}, j=1,\cds, m\big)\subset \si(\tilde B_{[0,\d]})$. Moreover, by otherwise extending the probability space further, it is straightforward to extend $\tilde \xi_\pi$ to a continuous process $\xi$ such that $\cL_\xi = \mu$ and $\xi_{s_j} = \tilde \xi_{s_j}$, $j=1,\cds, m$, $\tilde \dbP$-a.s. Finally, since $\nu_\pi$ is degenerate, we have
\beaa
\|\xi-\eta\|_{\tilde\dbP,\pi}^2 = \dbE^{\tilde \dbP}\big[\max_{1\le j\le n} |\xi_{s_j}-x_j|^2\big] =  \dbE^{\mu}\big[\max_{1\le j\le n} |X_{s_j}-x_j|^2\big]  = \|X - X'\|_{\overline \dbP, \pi}^2.
\eeaa
 This verifies all the requirements in (i)-(iii). 
 
\textit{\underline{\smash{Case 2}}:} $\nu_\pi$ is discrete, namely  $\nu_\pi =\sum_{i\ge 1} p_i\d_{\{(x^i_1,\cds, x^i_m)\}}$, with $p_i> 0$ and $\sum_{i\ge 1} p_i=1$. Fix a partition $\{O_i\}_{i\ge 1}\subset \cB(\dbR^{m})$ of $\dbR^m$ such that $(x^i_1,\cds,x^i_m)\in O_i$. Let $\tilde B^i_{[0,\d]}$ be a sequence of independent Brownian motions such that they are all independent of $\eta$. For each $i$,  define a conditional probability: 
\beaa
\mu^i(E) := {1\over p_i} \overline \dbP\Big( X_\pi \in E, X'_\pi \in O_i\Big),\q E\in \cB(\dbR^m).
\eeaa
Then by Case 1, one may construct a random vector $\tilde \xi^i_\pi =  \f_i(\tilde B^i_{[0,\d]})$ measurable to $\si(\tilde B^i_{[0,\d]})$ such that $\cL_{\tilde \xi^i_\pi} = \mu^i$. Define
\beaa
\tilde B := \sum_{i\ge 1} \tilde B^i 1_{O_i}(\eta_\pi),\q  \tilde \xi_\pi := \sum_{i\ge 1} \tilde \xi^i_\pi 1_{O_i}(\eta_\pi).
\eeaa
We now verify the desired properties. First, 
since all $\tilde B^i$ are independent of $\eta$, then $\tilde B$ is also a $\tilde \dbP$-Brownian motion. Moreover, for any $\tilde \pi: 0=t_0<\cds<t_n =\d$ and any $E, \tilde E \in \cB(\dbR^n)$, 
\beaa
\tilde \dbP(\tilde B_{\tilde \pi} \in E, \eta_{\tilde \pi} \in \tilde E) \!\!\!\!&=& \!\!\!\!\sum_{i\ge 1} \tilde \dbP\big(\tilde B^i_{\tilde \pi} \in E, \eta_{\tilde \pi} \in \tilde E, \eta_{\pi}\in O_i\big) = \sum_{i\ge 1} \tilde \dbP\big(\tilde B^i_{\tilde \pi} \in E\big) \tilde \dbP\big(\eta_{\tilde \pi} \in \tilde E, \eta_{\pi}\in O_i\big) \\
\!\!\!\! &=&\!\!\!\! \sum_{i\ge 1} \tilde \dbP\big(\tilde B_{\tilde \pi} \in E\big) \tilde \dbP\big(\eta_{\tilde \pi} \in \tilde E, \eta_{\pi}\in O_i\big) = \tilde \dbP\big(\tilde B_{\tilde \pi} \in E\big) \tilde \dbP\big(\eta_{\tilde \pi} \in \tilde E\big).
\eeaa
That is, $\tilde B$ is also independent of $\eta$. Next, since $O_i$ is a partition, we see that $\tilde \xi_\pi := \sum_{i\ge 1}\f_i(\tilde B^i_{[0,\d]}) 1_{O_i}(\eta_\pi) =  \sum_{i\ge 1} \f_i(\tilde B_{[0,\d]}) 1_{O_i}(\eta_\pi)$ and thus $\tilde \xi_\pi$ is measurable to $\si(\eta_\pi, \tilde B_{[0,\d]})$. Moreover, note that $\tilde \xi^i$s are also independent of $\eta$, then
\beaa
\tilde \dbP(\tilde \xi_\pi \in E) &=& \sum_{i\ge 1} \tilde \dbP(\tilde \xi^i_\pi \in E, \eta_\pi \in O_i) = \sum_{i\ge 1} \tilde \dbP(\tilde \xi^i_\pi \in E) \tilde \dbP(\eta_\pi \in O_i)\\
&=&\sum_{i\ge 1}\mu^i(E)  p_i = \sum_{i\ge 1}\overline \dbP\Big( X_\pi \in E, X'_\pi \in O_i\Big) = \overline \dbP\Big( X_\pi \in E\Big) = \mu(X_\pi\in E).
\eeaa
That is, $\cL_{\tilde\xi_\pi} = \mu_\pi$. Then similar to Case 1, by extending the space if necessary, we may construct $\xi$ such that $\cL_\xi = \mu$ and $\xi_\pi = \tilde \xi_\pi$, $\tilde\dbP$-a.s. Finally, 
\beaa
&&\|\xi-\eta\|_{\tilde\dbP,\pi}^2 = \dbE^{\tilde \dbP}\big[\max_{1\le j\le n} |\xi_{s_j}-\eta_{s_j}|^2\big] = \sum_{i\ge 1}   \dbE^{\tilde \dbP}\big[\max_{1\le j\le n} |\xi_{s_j}-\eta_{s_j}|^2 1_{O_i}(\eta_\pi)\big] \\
&&=\sum_{i\ge 1}   \dbE^{\tilde \dbP}\big[\max_{1\le j\le n} |\tilde \xi^i_{s_j}-x^i_j|^2 1_{O_i}(\eta_\pi)\big] = \sum_{i\ge 1}   \dbE^{\tilde \dbP}\big[\max_{1\le j\le n} |\tilde \xi^i_{s_j}-x^i_j|^2\big] \tilde \dbP(\eta_\pi\in O_i)\\
&&= \sum_{i\ge 1}   \dbE^{\mu^i}\big[\max_{1\le j\le n} |X_{s_j}-x^i_j|^2\big] p_i  =\sum_{i\ge 1} \dbE^{\overline \dbP}\big[\max_{1\le j\le n} |X_{s_j}-x^i_j|^2 1_{O_i}(X'_\pi) \big]  \\
&&= \sum_{i\ge 1} \dbE^{\overline \dbP}\big[\max_{1\le j\le n} |X_{s_j}-X'_{s_j}|^2 1_{O_i}(X'_\pi) \big] =\dbE^{\overline \dbP}\big[\max_{1\le j\le n} |X_{s_j}-X'_{s_j}|^2 \big] = \|X - X'\|_{\overline \dbP, \pi}^2.
\eeaa

\textit{\underline{\smash{Case 3}}:} We now consider the general case. Let $\{O_i\}$ be a countable partition of $\dbR^m$ such that for each $i$, the diameter of $O_i$ is less than $\eps/2$. For each $i$, fix an arbitrary $x^i\in O_i$ and denote $p_i := \nu_\pi(O_i)$. By otherwise eliminating some $i$, we may assume $p_i >0$ for all $i$. Denote $\tilde \eta_\pi := \sum_{i\ge 1} x^i 1_{O_i}(\eta_\pi)$ and $\tilde X'_\pi :=  \sum_{i\ge 1} x^i 1_{O_i}(X'_\pi)$.  By Case 2, there exist a $\tilde \dbP$-Brownian motion $\tilde B_{[0, \d]}$  and a continuous process $\xi$ such that 

$\bullet$ $\cL_\xi = \mu$ and $\tilde B$ is independent of $\tilde \eta_\pi$. Moreover, from the arguments we may assume further that  $\tilde B$ is independent of $\eta$;

$\bullet$ Each $\xi_{s_j}$ is measurable to $\si\big(\tilde \eta_\pi, \tilde B_{[0,\d]}\big) \subset \si\big(\eta_\pi, \tilde B_{[0,\d]}\big)$;
 
$\bullet$ $\dbE^{\tilde \dbP} \big[\max_{1\le j\le m} |\xi_{s_j} - \tilde \eta_{s_j}|^2\big] = \dbE^{\overline \dbP} \big[\max_{1\le j\le m} |X_{s_j} - \tilde X'_{s_j}|^2\big]$.

\no This verifies (i) and (ii). To see (iii), note that $|\eta_{s_j} - \tilde \eta_{s_j}|\le {\e\over 2}$, $|X'_{s_j} - \tilde X'_{s_j}|\le {\e\over 2}$. Then 
\beaa
\|\xi-\eta\|_{\tilde\dbP,\pi} &=& \Big(\dbE^{\tilde \dbP}\big[\max_{1\le j\le n} |\xi_{s_j}-\eta_{s_j}|^2\big] \Big)^{1\over 2}\le  \Big(\dbE^{\tilde \dbP}\big[\max_{1\le j\le n} |\xi_{s_j}-\tilde \eta_{s_j}|^2 \big]\Big)^{1\over 2} + {\e\over 2} \\
&=&\Big( \dbE^{\overline \dbP} \big[\max_{1\le j\le m} |X_{s_j} - \tilde X'_{s_j}|^2\big]\Big)^{1\over 2} + {\e\over 2}   \le \Big( \dbE^{\overline \dbP} \big[\max_{1\le j\le m} |X_{s_j} -  X'_{s_j}|^2\big]\Big)^{1\over 2} + \e\\  
&=& \|X - X'\|_{\overline \dbP, \pi} +\e.
\eeaa
This completes the proof.
\qed

\begin{rem}
\label{rem-construct}
{\rm (i) As mentioned right before the lemma, the main difficulty of establishing the regularity of $V$ at $\nu$ is due to the possible degeneracy of $\nu$, and thus in the above lemma one may not be able to write $\xi$ as a function of $\eta$. Our trick here is to introduce the independent Brownian motion $\tilde B_{[0,\d]}$ (which always has continuous distribution) and then Lemma \ref{lem-construct} (ii) holds. 

(ii) The construction of $\xi$, which relies on \reff{fzeta}, works only for finite dimensional random vectors. It is not clear to us how to generalize this result to the case where the $m$-tuple $(s_1,\cds,s_m)$ is replaced by the uncountable interval $[0,t]$. This is why we need to consider value function $ V_0$ first.
\qed}
\end{rem}

\begin{lem}
\label{lem-hatVreg}
Under Assumption \ref{assum-standing}, $V_0$ is uniformly continuous in $\mu$, uniformly in $t$. That is, there exists a modulus of continuity function $\rho$ such that
\bea
\label{hatVLip}
|V_0(t,\mu) - V_0(t, \nu)| \le \rho\big( \cW_2(\mu_{[0,t]}, \nu_{[0,t]})\big),\q \mbox{for all}~t\in [0, T], \mu, \nu \in \cP_2.
\eea
\end{lem}
\proof Let's fix $t\in [0, T]$,  $\mu, \nu \in \cP_2$, $\a\in \cA^0_t$, and $\e, \d>0$. Choose  $\overline \dbP\in \cP(\mu,\nu)$ such that 
\bea
\label{oldbP}
\dbE^{\ol \dbP}[\|X_{t\wedge \cd}-X'_{t\wedge \cd}\|^2] \le c_0^2  + \e^2,\q \mbox{where}\q c_0:= \cW_2(\mu_{[0,t]}, \nu_{[0,t]}).
\eea
 Our idea is to construct some $\tilde \a\in \cA^0_t$ such that $\dbP^{t,\nu, \tilde \a}$ is close to $\dbP^{t,\mu,  \a}$ in certain way.

By \reff{hatcAt}, we assume  $\a = \sum_{i=0}^{n-1} h^0_i(X_{\pi_0}, X_{[t, t_i]})\1_{[t_i, t_{i+1}]}$, where $\pi_0: 0\le s^0_1<\cds<s^0_{m_0} = t$ and $t=t_0<\cds<t_n = T$. We shall fix $n$, and assume $\d <  \min_{1\le i\le n}[t_i-t_{i-1}]$. But to obtain a desired approximation, we shall consider finer partitions $\pi: 0\le s_1<\cds<s_{m} = t$ such that $\pi_0 \subset \pi$. Clearly, we may rewrite $\a = \sum_{i=0}^{n-1} h_i(X_{\pi}, X_{[t, t_i]})\1_{[t_i, t_{i+1}]}$. Let $(\tilde \O, \tilde \cF, \tilde \dbP)$, $\tilde B$, $\xi$ and $\eta$ be as in Lemma \ref{lem-construct}, corresponding to $(t,\mu, \nu, \pi, \e, \d, \overline \dbP)$. Denote $B'_s := \tilde B_{s-t}$, $B^\d_s := \tilde B_{\d, s-t+\d} $, $s\in [t, T]$.  It is clear that $\dbP^{t,\mu,\a} = \tilde \dbP \circ (X^\a)^{-1}$, where $X^\a_{[0,t]}= \xi_{[0,t]}$ and, for $i=0,\cds, n-1$ and $s\in (t_i, t_{i+1}]$,
\bea
\label{Xai}
X^\a_s = X^\a_{t_i} + \int_{t_i}^s b(r, X^\a, \cL_{X^\a}, h_i(\xi_\pi, X^\a_{[t, t_i]})) dr +  \int_{t_i}^s \si(r, X^\a, \cL_{X^\a}, h_i(\xi_\pi, X^\a_{[t, t_i]}) dB^\d_r.
\eea

{\it Step 1.} We first construct $\tilde \a \in \cA^0_t$ and $\tilde X := X^{t, \eta, \tilde \a}$ satisfying $\tilde X_{[0,t]} := \eta_{[0,t]}$ and
\bea
\label{tildeX}
d \tilde X_s :=  b(s, \tilde X, \cL_{\tilde X},  \tilde \a_s) ds +  \si(s, \tilde X, \cL_{\tilde X}, \tilde \a_s) dB'_s,\q s\ge t,\q\tilde \dbP\mbox{-a.s.}
\eea
The corresponding partitions for $\tilde \a$ will be $\pi$ and $t=t_0 < t_0+\d<t_1+\d<\cds<t_{n-1}+\d < t_n = T$. First, fix an arbitrary $a_0\in A$ and set $\tilde\a_s := a_0$ for $s\in [t_0, t_0+\d)$. Then we may determine $\tilde X$ on $[t_0, t_0+\d]$ by \reff{tildeX} with initial condition $\tilde X_{t_0} = \eta_{t_0}$.  Since the SDE \reff{tildeX} has a strong solution and $\si$ is non-degenerate, we know the $\si$-algebras $\si(\eta_\pi, \tilde B_{[0,\d]}) = \si(\eta_\pi, \tilde X_{[t, t_0+\d]})$ (abusing the notation $\si$ here!). Then, by Lemma \ref{lem-construct} (ii), $\xi_\pi = \f(\eta_\pi, \tilde X_{[t, t_0+\d]})$ for some function $\f$. Set $\tilde h_0(\eta_\pi,  \tilde X_{[t, t_0+\d]}) := h_0( \f(\eta_\pi, \tilde X_{[t, t_0+\d]})) = h_0(\xi_\pi)$. Then, for $s\in [t_0+\d, t_1+\d)$, setting $\tilde \a_s := \tilde h_0(\eta_\pi,  \tilde X_{[t, t_0+\d]}) = h_0(\xi_\pi)$,  we may determine $\tilde X$ further on $[t_0+\d, t_1+\d]$ by \reff{tildeX}. Next, again since $\si$ is nondegenerate, we see that $X^\a_{[t, t_1]}$ is measurable to 
\beaa
\si(\xi_\pi, B^\d_{[t_0, t_1]}) \subset \si(\eta_\pi, \tilde B_{[0, \d+t_1-t_0]}) \subset \si(\eta_\pi, \tilde X_{[t, t_1+\d]}).
\eeaa
 Then $h_1(\xi_\pi, X^\a_{[t, t_1]})= \tilde h_1(\eta_\pi, \tilde X_{[t, t_1+\d]})$ for some function $\tilde h_1$. For $s\in (t_1+\d, t_2+\d]$, set $\tilde \a_s :=  \tilde h_1(\eta_\pi, \tilde X_{[t, t_1+\d]}) = h_1(\xi_\pi, X^\a_{[t, t_1]})$. Repeat the arguments, we may construct $\tilde \a\in \cA^0_t$ such that, for the corresponding $\tilde X$ determined by \reff{tildeX}, 
\bea
\label{tildea}
\dis \tilde \a_s = \tilde h_i(\eta_\pi, \tilde X_{[t, t_i+\d]}) = h_i(\xi_\pi, X^\a_{[t, t_i]}),\q s\in [t_i+\d, t_{i+1}+\d).
\eea

{\color{black}{\it Step 2.} We next estimate the difference between $X^\a$ and $\tilde X$. Denote 
\bea
\label{OSC}
\D X_s := \tilde X_s - X^\a_s,~ \D_\d X_s:= \tilde X_{s+\d} - X^\a_s,~ OSC_\d(\tilde X) := \sup_{t\le r_1< r_2 \le T, r_2 - r_1 \le \d} |\tilde X_{r_1,r_2} |.
\eea
Note that, for $s\in [t, T-\d]$, 
\beaa
 \|\tilde X_{(s+\d)\wedge\cd } - X^\a_{s\wedge \cd}\| \le \|\xi_{t\wedge \cd} - \eta_{t\wedge \cd}\| + \sup_{t\le r\le s} |\D_\d X_r| + OSC_\d(\tilde X).
 \eeaa
 By Assumption \ref{assum-standing} (i) and (iii),  for $\f=b, \si, f$, and $s\in [t_i, t_{i+1}]$, 
\bea
\label{tildeX1}
&&\dbE^{\tilde \dbP}\Big[\big|\f(s+\d, \tilde X, \cL_{\tilde X}, h_i(\xi_\pi, X^\a_{[t, t_i]}) ) - \f(s, X^\a, \cL_{X^\a}, h_i(\xi_\pi, X^\a_{[t, t_i]}) )\big|^2\Big]\nonumber\\
&&\le C\dbE^{\tilde \dbP}\Big[\big|\f(s+\d, \tilde X, \cL_{\tilde X}, h_i(\xi_\pi, X^\a_{[t, t_i]}) ) - \f(s+\d, X^\a_{s\wedge \cd}, \cL_{X^\a_{[0, s]}}, h_i(\xi_\pi, X^\a_{[t, t_i]}) )\big|^2\Big]\nonumber\\
&&\q + C\dbE^{\tilde \dbP}\Big[\big| \f(s+\d, X^\a_{s\wedge \cd}, \cL_{X^\a_{[0, s]}}, h_i(\xi_\pi, X^\a_{[t, t_i]}) ) -\f(s, X^\a, \cL_{X^\a}, h_i(\xi_\pi, X^\a_{[t, t_i]}) )\big|^2\Big]\nonumber\\
&&\le C\dbE^{\tilde \dbP}\Big[\big| \|\tilde X_{(s+\d)\wedge \cd} - X^\a_{s\wedge \cd}\|^2 + \big[1+ \|X^\a_{s\wedge \cd}\|^2\big] \rho^2_0(\d)\Big] \nonumber\\
&&\le C\dbE^{\tilde \dbP}\Big[ \|\|\xi_{t\wedge \cd} - \eta_{t\wedge \cd}\|^2 +  \sup_{t\le r\le s} |\D_\d X_r|^2  + OSC^2_\d(\tilde X)  \Big] + C_\mu\rho^2_0(\d),
\eea
Note that we may rewrite \reff{tildeX} as, for $s\in [t_i, t_{i+1}]$,
\beaa
\tilde X_{s+\d}=  \tilde X_{t_i+\d} + \int_{t_i}^{s} b(r+\d, \tilde X, \cL_{\tilde X}, h_i(\xi_\pi, X^\a_{[t, t_i]}) )dr +  \int_{t_i}^s \si(r+\d, \tilde X, \cL_{\tilde X}, h_i(\xi_\pi, X^\a_{[t, t_i]}) ) dB^\d_r.
\eeaa
Compare this with \reff{Xai},  then it follows from standard arguments that 
\bea
\label{tildeX2}
\dbE^{\tilde \dbP}\Big[\sup_{t\le s\le T-\d} |\D_\d X_s| ^2\Big] \le C\dbE^{\tilde \dbP}\Big[ |\D_\d X_t|^2 + \|\xi_{t\wedge \cd} - \eta_{t\wedge \cd}\|^2   + OSC^2_\d(\tilde X)  \Big] + C_\mu\rho^2_0(\d)
\eea
Note that, for $s\in [t, T-\d]$,
\beaa
&&|\D_\d X_t| = |\tilde X_{t+\d} - \tilde X_t| + |\xi_t - \eta_t| \le  \|\xi_{t\wedge \cd} - \eta_{t\wedge \cd}\|   + OSC_\d(\tilde X);\\
&&|\D X_s| \le  |\D_\d X_s| + OSC_\d(\tilde X),\q s\in  [t, T-\d];\\
&& |\D X_s| \le |\D_\d X_{T-\d}| + |\tilde X_s - \tilde X_T| + |X^\a_s - X^\a_{T-\d}| \\
&&\qq\q \le |\D_\d X_{T-\d}| + OSC_\d(\tilde X)+ OSC_\d(X^\a),~ s\in [T-\d, T],
\eeaa
where $OSC_\d(X^\a)$ is defined similar to \reff{OSC}. Then \reff{tildeX2} leads to
\beaa
\dbE^{\tilde \dbP}\Big[\sup_{t\le s\le T} |\D X_s| ^2\Big] \le C\dbE^{\tilde \dbP}\Big[\|\xi_{t\wedge \cd} - \eta_{t\wedge \cd}\|^2   + OSC^2_\d(\tilde X) + OSC^2_\d(X^\a)  \Big] + C_\mu\rho^2_0(\d).
\eeaa
Since $|b|, |\si|\le C_0$, by Revuz \& Yor \cite[Chapter I, Theorem 2.1]{RY}
one can easily see that
\bea
\label{tildeX3}
\dbE^{\tilde \dbP}\Big[  |OSC_\d(\tilde X)|^2+ |OSC_\d(X^\a)|^2\Big]\le C\sqrt{\d} \le C\rho^2_0(\d).
\eea
Here we assume without loss of generality that $\rho_0(\d) \ge \d^{1\over 4}$ (otherwise replace $\rho_0$ with $\rho_0(\d) \vee \d^{1\over 4})$. 
Then, noting that $\|\D X\|\le \sup_{t\le s\le T}|\D X_s|+ \|\xi_{t\wedge \cd} - \eta_{t\wedge \cd}\|$, 
\bea
\label{tildeX4}
\dbE^{\tilde \dbP}\Big[\|\D X\| ^2\Big] \le C\dbE^{\tilde \dbP}\Big[ \|\xi_{t\wedge \cd} - \eta_{t\wedge \cd}\|^2\Big]  + C_\mu \rho_0^2(\d).
\eea

{\it Step 3.} We now estimate $V_0(t,\mu) - V_0(t, \nu)$. By Assumption \ref{assum-standing} (ii) and (iii), we have
\bea
\label{tildeX5}
&&J(t, \mu, \a)  - V_0(t, \nu) \le J(t,\mu,\a) - J(t,\nu,\tilde \a)\nonumber \\
&&=  \dbE^{\tilde \dbP} \Big[ g(X^\a, \cL_{X^\a})- g(\tilde X, \cL_{\tilde X})  + \sum_{i=0}^{n-1} \int_{t_i}^{t_{i+1}} f(s, X^\a, \cL_{X^\a}, h_i(\xi_\pi, X^\a_{[t, t_i]})) ds \nonumber\\
&&\q -\sum_{i=0}^{n-2} \int_{t_i}^{t_{i+1}} f(s+\d, \tilde X, \cL_{\tilde X}, h_i(\xi_\pi, X^\a_{[t, t_i]})) ds \nonumber\\
&&\q - \int_t^{t+\d} f(s, \tilde X, \cL_{\tilde X}, a_0) ds -  \int_{T-\d}^T f(s, \tilde X, \cL_{\tilde X}, \tilde \a_s) ds\Big]\nonumber\\
&&\le C\dbE^{\tilde \dbP}\Big[\rho_0\big(\|\D X\|  + \cW_2(\cL_{X^\a}, \cL_{\tilde X})\big)\nonumber \\
&& + [1+\|\tilde X\| + \|X^\a\|+\cW_2(\cL_{\tilde X}, \d_{\{0\}}) +\cW_2(\cL_{X^\a}, \d_{\{0\}}) ] \rho_0(\d)\Big] \nonumber\\
&&\le C\dbE^{\tilde \dbP}\Big[\rho_0\big(\|\D X\|  + \cW_2(\cL_{X^\a}, \cL_{\tilde X})\big)\Big] + C_{\mu,\nu} \rho_0(\d).
\eea
}
Note that we may assume without loss of generality that $\rho_0$ has linear growth. Then, 
\bea
\label{DJ}
&&J(t, \mu, \a)  - V_0(t, \nu) \le  C\dbE^{\tilde \dbP}\Big[\rho_0\big(\|\D X\|  + \cW_2(\cL_{X^\a}, \cL_{\tilde X})\big) \1_{\{\|\D X\|>c_0\}}\Big] \nonumber\\
&&\qq\qq + C\rho_0(c_0+ \cW_2(\cL_{X^\a}, \cL_{\tilde X}))  + C_{\mu,\nu} \rho_0(\d)\\
&&\le {C\over c_0} \dbE^{\tilde \dbP}[\|\D X\|^2] + C\rho_0\big(c_0+  (\dbE^{\tilde \dbP}[\|\D X\|^2])^{1\over 2}\big)  + C_{\mu,\nu} \rho_0(\d),\nonumber
\eea
where $c_0$ is defined by \reff{oldbP}.  Note further that, denoting $|\pi| := \min_{1\le j\le m} |s_j-s_{j-1}|$,
\beaa
 \|\xi_{t\wedge \cd} - \eta_{t\wedge \cd}\| \le OSC_{|\pi|}(\xi_{[0,t]})+OSC_{|\pi|}(\eta_{[0,t]}) + \max_{0\le j\le m}|\xi_{s_j} - \eta_{s_j}|.
\eeaa
Plug this into \reff{tildeX4}, by Lemma \ref{lem-construct} (iii) and \reff{oldbP} we have
\beaa
&&\dbE^{\tilde \dbP}\Big[\|\D X\| ^2\Big] \le C\dbE^{\tilde \dbP}\Big[OSC^2_{|\pi|}(\xi_{[0,t]}) + OSC^2_{|\pi|}(\eta_{[0,t]})\Big] + \|\xi-\eta\|_{\tilde \dbP, \pi}^2  + C_\mu \rho_0^2(\d)\\
&&\le C\Big[\dbE^{\mu}[OSC^2_{|\pi|}(X_{[0,t]})] + \dbE^{\nu}[OSC^2_{|\pi|}(X_{[0,t]})]  + \|X_{t\wedge \cd} -X'_{t\wedge \cd}\|^2_{\ol\dbP,\pi}  + \e^2  \Big] +C_\mu \rho_0^2(\d)\\
&&\le C\Big[\dbE^{\mu}[OSC^2_{|\pi|}(X_{[0,t]})] + \dbE^{\nu}[OSC^2_{|\pi|}(X_{[0,t]})]  + c_0^2  + \e^2  \Big] + C_\mu\rho_0^2(\d).
\eeaa
Plug this into \reff{DJ}, and note that $\a$ depends on $\pi_0$, but not $\pi$. Then, by sending $\d\to 0$, $\e\to 0$,  and $|\pi|\to 0$, we obtain:
\beaa
J(t, \mu, \a)  - V_0(t, \nu) &\le&  {C\over c_0}\cW^2_2(\mu_{[0,t]}, \nu_{[0,t]}) + C\rho_0\big(c_0+  \cW_2(\mu_{[0,t]},\nu_{[0,t]})\big)  \\
&=&C\cW_2(\mu_{[0,t]}, \nu_{[0,t]}) + C\rho_0\big(2 \cW_2(\mu_{[0,t]},\nu_{[0,t]})\big). 
\eeaa
Now by the arbitrariness of $\a\in \cA^0_t$, we obtain 
\beaa
 V_0(t, \mu) - V_0(t, \nu) \le C\cW_2(\mu_{[0,t]}, \nu_{[0,t]}) + C\rho_0\big(2 \cW_2(\mu_{[0,t]},\nu_{[0,t]})\big).
 \eeaa
 Following the same arguments we also have the estimate for $V_0(t, \nu) - V_0(t, \mu)$, and thus complete the proof.
 \qed

\begin{lem}
\label{lem-VhatV}
 Under Assumption \ref{assum-standing}, we have $V=V_0$.
\end{lem}
\proof By definition, it is clear that $V_0 \le V$. To prove the opposite inequality, we fix $(t,\mu) \in \Th$ and  $\a := \sum_{i=0}^{n-1} h_i(X_{[0, t_i]}) 1_{[t_i, t_{i+1})}\in \cA_t$ as in \reff{cAt}. Fix $\dbP_0$, $B$, and $\xi$ being such that $\dbP_0 \circ (\xi_{[0, t]})^{-1} = \mu$, and let $X^\a= X^{t,\xi,\a}$ be defined by \reff{Xa1}. We shall prove $J(t,\mu,\a)  \le V_0(t,\mu)$ in two steps.

{\it Step 1.} We first assume all the functions $h_i: C([0, t_i])\to \dbR$ are continuous. For each $m\ge 1$, consider the partition $\pi_m: 0=s^m_0<\cds<s^m_m=t$ be such that $s^m_i = {i\over m} t$. Define 
\bea
\label{hm1}
\left.\ba{c}
\dis h^m_i(\eta_{\pi_m}, \eta_{[t, t_i]}) := h_i(\eta^m_{[0,t_i]}), \\
\dis\mbox{where $\eta^m_{[0,t]}$ is the linear interpolation of $\eta_{\pi_m}$ and $\eta^m_{[t, t_i]} := \eta_{[t, t_i]}$}.
\ea\right.
\eea
Denote $\a^m := \sum_{i=0}^{n-1} h^m_i( X_{\pi_m}, X_{[t, t_i]}) 1_{[t_i, t_{i+1})} \in \cA^0_t$, and define $X^m := X^{t,\xi, \a^m}$ in an obvious way. We shall estimate $\D X^m := X^m - X^\a$. 

Clearly $\D X^m_s = 0$ for $s\in [0, t]$. For $ s\in [t_0, t_1]$, we have
\beaa
X^\a_s &=& \xi_t + \int_{t_0}^s b(r,  X^\a, \cL_{X^\a}, h_0(\xi_{[0, t]}) dr + \int_{t_0}^s \si(r, X^\a, \cL_{X^\a}, h_0(\xi_{[0, t]}) dB_r;\\
X^m_s &=& \xi_t + \int_{t_0}^s b(r, X^m, \cL_{X^m}, h^m_0(\xi_{\pi_m}, \xi_t)) dr + \int_{t_0}^s \si(r, X^m, \cL_{X^m}, h^m_0(\xi_{\pi_m}, \xi_t)) dB_r.
\eeaa
Since $h_0$ is continuous, it is clear that 
\beaa
\lim_{m\to\infty} \dbE^{\dbP_0}\Big[\Big|h^m_0(\xi_{\pi_m}, \xi_t)  - h_0(\xi_{[0,t]})\Big|^2\wedge 1\Big] = 0.
\eeaa
By  Assumption \ref{assum-standing} (i) and (iii), it follows from standard arguments that 
\bea
\label{DXmconv}
\lim_{m\to\infty} \dbE^{\dbP_0}\Big[\|\D X^m_{t_1\wedge \cd}\|^2\Big] = 0.
\eea
Next, for $s\in [t_1, t_2]$, we have
\beaa
&X^\a_s = X^\a_{t_1}+ \int_{t_1}^s b(r, X^\a, \cL_{X^\a}, h_1(\xi_{[0, t]}, X^\a_{[t, t_1]})) dr + \int_{t_1}^s \si(r, X^\a,\cL_{X^\a}, h_1(\xi_{[0, t]}, X^\a_{[t, t_1]}) dB_r;\\
&X^m_s = X^m_{t_1} + \int_{t_1}^s b(r, X^m,\cL_{X^m}, h^m_1(\xi_{\pi_m}, X^m_{[t, t_1]})) dr + \int_{t_0}^s \si(r, X^m,\cL_{X^m}, h^m_1(\xi_{\pi_m}, X^m_{[t, t_1]})) dB_r.
\eeaa
Since $h_1$ is continuous, by \reff{hm1} and  \reff{DXmconv} we have
\beaa
\lim_{m\to\infty} \dbE^{\dbP_0}\Big[\Big|h^m_1(\xi_{\pi_m}, X^m_{[t, t_1]})  - h_1(\xi_{[0, t]}, X^\a_{[t, t_1]})\Big|^2 \wedge 1\Big] = 0.
\eeaa
Then, similar to \reff{DXmconv} we have $\dis\lim_{m\to\infty} \dbE^{\dbP_0}\Big[\|\D X^m_{t_2\wedge \cd}\|^2\Big] = 0$. Repeat the arguments we obtain 
\beaa
\lim_{m\to\infty} \dbE^{\dbP_0}\Big[\|\D X^m\|^2\Big] = 0,~  \lim_{m\to\infty} \dbE^{\dbP_0}\Big[\Big|h^m_i(\xi_{\pi_m}, X^m_{[t, t_i]})  - h_i(\xi_{[0, t]}, X^\a_{[t, t_i]})\Big|^2 \wedge 1\Big] = 0,  i<n.
\eeaa
Now by the regularity of $f$ and $g$ in  Assumption \ref{assum-standing} (ii) and (iii), we have 
\bea
\label{JV0}
J(t,\mu,\a) =  \lim_{m\to\infty} J(t, \mu, \a^m)  \le V_0(t,\mu).
\eea

{\color{black}
{\it Step 2.} We now consider the general Borel measurable functions $h_i$.  We shall construct $\a^m = \sum_{i=0}^{n-1} h^m_i(X_{[0, t_i]}) \1_{[t_i, t_{i+1})} \in \cA_t$ such that each $h^m_i$ is continuous and, for the corresponding $X^m:= X^{t,\xi, \a^m}$ and denoting $\D X^m := X^m-X^\a$, 
\bea
\label{amXm}
\lim_{m\to\infty} \dbE^{\dbP_0}[|h^m_i(X^m_{[0, t_i]}) - h_i(X^\a_{[0, t_i]})|^2 \wedge 1] = 0,\q \lim_{m\to\infty} \dbE^{\dbP_0} [\|\D X^m\|^2] =0.
\eea
Then by Step 1 we have $J(t, x, \a^m) \le V_0(t,\mu)$, and similar to \reff{JV0} we can easily show that $ J(t,\mu,\a) =  \lim_{m\to\infty} J(t, \mu, \a^m) \le V_0(t,\mu)$. 

We now construct $h^m_i$ recursively in $i$.  First, denote $X^m_{[0, t]} := \xi_{[0, t]}$. Then $\|\D X^m_{t_0\wedge}\|=0$. Assume by induction that we have constructed $X^m_{[0, t_i]}$ such that 
\bea
\label{DXm0}
\lim_{m\to\infty} \dbE^{\dbP_0} [\|\D X^m_{t_i \wedge \cd}\|^2] = 0.
\eea
 For $h_i: C([0, t_0]) \to A$, by Lusin's lemma, there exist continuous functions $\tilde h^m_i: C([0, t_i]) \to A$ and  closed sets $K^m_i \subset C([0, t_i])$ such that 
\bea
\label{DXm}
\tilde h^m_i = h_i~\mbox{on}~ K^m_i\q\mbox{and}\q \lim_{m\to \infty}\dbP_0(X^\a_{[0,t_i]}\notin K^m_i) =0.
\eea
For each $m$, since $\tilde h^m_i$ is continuous, by \reff{DXm} we have 
\beaa
\lim_{k\to \infty} \dbE^{\dbP_0}\Big[|\tilde h^m_i(X^{k}_{[0, t_i]} ) -\tilde h^m_i(X^\a_{[0, t_i]})|^2\wedge 1\Big] = 0.
\eeaa
Then there exists $k_m$ such that 
\beaa
 \dbE^{\dbP_0}\Big[|\tilde h^m_i(X^{k_m}_{[0, t_i]} ) - \tilde h^m_i(X^\a_{[0, t_i]})|^2\wedge 1\Big] \le {1\over m},\q \forall m.
\eeaa
This implies that
\beaa
&& \dbE^{\dbP_0}\Big[|\tilde h^m_i(X^{k_m}_{[0, t_i]} ) - h_i(X^\a_{[0, t_i]})|^2\wedge 1\Big] \\
 &&\le C\dbE^{\dbP_0}\Big[|\tilde h^m_i(X^{k_m}_{[0, t_i]} ) - \tilde h^m_i(X^\a_{[0, t_i]})|^2\wedge 1\Big]  + C\dbE^{\dbP_0}\Big[|\tilde h^m_i(X^\a_{[0, t_i]}) - h_i(X^\a_{[0, t_i]})|^2\wedge 1\Big] \\
 &&\le {C\over m} + C\dbP_0(X^\a_{[0,t_i]}\notin K^m_i)  \to 0,\q \mbox{as}~m\to \infty.
\eeaa
By considering the subsequence $k_m$ and set $h^{k_m}_i := \tilde h^m_i$, we obtain
\bea
\label{amXm2}
\lim_{m\to\infty} \dbE^{\dbP_0}[|h^{k_m}_i(X^{k_m}_{[0, t_i]}) - h_i(X^\a_{[0, t_i]})|^2 \wedge 1] = 0.
\eea
By choosing the subsequence $k_m$, and for notational simplicity, we assume $k_m=m$, then we constructed the desired $h^m_i$ under assumption \reff{DXm0}.

Next, for $s\in [t_i, t_{i+1}]$ and for $\f = b, \si$, denote
\beaa
\D \f^m_s := \f\big(s, X^\a, \cL_{X^\a}, h^m_i(X^m_{[0, t_i]})\big)- \f\big(s, X^\a, \cL_{X^{\a}}, h_i(X^\a_{[0, t_i]})\big).
\eeaa
Since $b$ and $\si$ are bounded and  uniformly Lipschitz continuous in $(\o, \mu)$,  we have
\beaa
&&\dbE^{\dbP_0}\Big[\Big|\f(s, X^m, \cL_{X^m}, h^m_i(X^m_{[0, t_i]}))-\f(s, X^\a, \cL_{X^{\a}}, h_i(X^\a_{[0, t_i]})\Big|^2\Big]\\
&&=\dbE^{\dbP_0}\Big[\Big|\big[\f(s, X^m, \cL_{X^m}, h^m_i(X^m_{[0, t_i]}))-\f(s, X^\a, \cL_{X^{\a}}, h^m_i(X^m_{[0, t_i]})) \big]+ \D \f^m_s\Big|^2\Big]\\
&&\le C\dbE^{\dbP_0}\Big[ \|\D X^m_{s\wedge \cd}\|^2  + |\D \f^m_s|^2 \Big].
\eeaa
Note that
\beaa
&\dis X^{\a}_s = X^{\a}_{t_i}+ \int_{t_i}^s b(r, X^{\a}, \cL_{X^{\a}}, h_i(X^{\a}_{[0, t_i]})) dr + \int_{t_i}^s \si(r, X^{\a}, \cL_{X^{\a}}, h_i(X^{\a}_{[0, t_i]})) dB_r;\\
&\dis X^{m}_s = X^{m}_{t_i}+ \int_{t_i}^s b(r, X^{m}, \cL_{X^{m}}, h^m_i(X^{m}_{[0, t_i]})) dr + \int_{t_i}^s \si(r, X^{m}, \cL_{X^{m}}, h^m_i(X^{m}_{[0, t_i]})) dB_r.
\eeaa
By standard arguments one can easily see that
\beaa
\dbE^{\dbP_0}[\|\D X^m_{t_{i+1}\wedge \cd}\|^2] \le C\dbE^{\dbP_0}\Big[\|\D X^m_{t_i\wedge \cd}\|^2 + \int_{t_i}^{t_{i+1}} [|\D b^m_s|^2 + |\D \si^m_s|^2] ds\Big].
\eeaa
By \reff{amXm2} (with $k_m = m$) and the dominated convergence theorem, we have 
\beaa
\lim_{m\to\infty}\dbE^{\dbP_0}\Big[ \int_{t_i}^{t_{i+1}} [|\D b^m_s|^2 + |\D \si^m_s|^2] ds\Big] =0.
\eeaa
This, together with  \reff{DXm0}, implies that $\lim_{m\to\infty}\dbE^{\dbP_0}[\|\D X^m_{t_{i+1}\wedge \cd}\|^2] =0$. Then the induction procedure can continue, and by possibly choosing a subsequence, we construct the desired $h^m_i$ satisfying \reff{amXm} for all $i$, hence completing the proof.
\qed
}

\bs
\no{\bf Proof of Theorem \ref{thm-regmu}.} First, by Lemmas \ref{lem-hatVreg} and \ref{lem-VhatV}, we see that $V$ is uniformly continuous in $\mu$ with certain modulus of continuity function $\rho$. Now let $t_1<t_2$ and $\mu, \nu \in \cP_2$. 
{\color{black}By DPP \reff{meancontrol-DPP} (see Remark \ref{rem-closed} (ii)) and noting that $f$ has linear growth in $(\o, \mu)$, we have
\bea
\label{Vreg2}
&&|V(t_1,\mu)-V(t_2,\nu)| \nonumber\\
&&\le \sup_{\a\in\cA_{t_1}} \Big[\big|V(t_2, \dbP^{t_1, \mu, \a}) - V(t_2, \nu)\big|+ \int_{t_1}^{t_2} \dbE^{\dbP^{t_1,\mu,\a}}\big[\big|f(s, X, \dbP^{t_1,\mu,\a}, \a_s)\big|\big] ds\Big]\\
&&\le \rho\big(\cW_2(\dbP^{t_1, \mu, \a}_{[0, t_2]}, \nu_{[0, t_2]})\big) + C_0 \int_{t_1}^{t_2} \big[1 + \dbE^{\dbP^{t_1,\mu,\a}}[\|X_{s\wedge \cd}\|] + \cW_2(\dbP^{t_1, \mu, \a}_{[0, s]}, \d_{\{0\}})\big] ds.\nonumber
\eea
Note that, since $b$ and $\si$ are bounded, for $s\in [t_1. t_2]$,
\beaa
&&\cW_2(\dbP^{t_1, \mu, \a}_{[0, t_2]}, \nu_{[0, t_2]})  \le \cW_2( \mu_{[0, t_1]}, \nu_{[0, t_2]}) + \cW_2(\dbP^{t_1, \mu, \a}_{[0, t_2]}, \mu_{[0, t_1]}) \\
&&\le \cW_2( \mu_{[0, t_1]}, \nu_{[0, t_2]}) + \Big(\dbE^{\dbP^{t_1, \mu, \a}}\big[\sup_{t_1 \le s\le t_2} |X_s - X_{t_1}|^2\big]\Big)^{1\over 2}\\
&&\le \cW_2( \mu_{[0, t_1]}, \nu_{[0, t_2]})   + C(t_2-t_1)^{1\over 2} \le C\cW_2((t_1, \mu), (t_2, \nu);\\
&&\Big(\dbE^{\dbP^{t_1,\mu,\a}}\big[\|X_{s\wedge \cd}\| \big] + \cW_2(\dbP^{t_1, \mu, \a}_{[0, s]}, \d_{\{0\}})\Big)^2 \le   C\dbE^{\dbP^{t_1, \mu, \a}}[\|X_{s\wedge \cd}\|^2] \\
&& \le C\dbE^{\dbP^{t_1,\mu,\a}}\big[ \|X_{t_1\wedge \cd}\|^2  + \sup_{t_1\le r\le s} |X_{t_1, r}|^2\big]  = C\dbE^\mu [ \|X_{t_1\wedge \cd}\|^2] +  C \dbE^{\dbP^{t_1,\mu,\a}}\big[  \sup_{t_1\le r\le s} |X_{t_1, r}|^2\big] \\
&& \le C\dbE^\mu[\|X_{t_1\wedge \cd}\|^2] + C[t_2-t_1] \le C\dbE^\mu[\|X_{t_1\wedge \cd}\|^2] + C.
\eeaa
Then
\beaa
&&|V(t_1,\mu) - V(t_2, \nu)| \le  \rho\big(C\cW_2((t_1, \mu), (t_2, \nu)\big)  + C\Big(1+\dbE^\mu[\|X_{t_1\wedge \cd}\|^2] \Big)^{1\over 2} [t_2-t_1].
\eeaa
This proves \reff{Vreg}. 

Moreover, if $f$ is bounded, then \reff{Vreg2} implies that 
\beaa
|V(t_1,\mu)-V(t_2,\nu)| &\le& \rho\big(\cW_2(\dbP^{t_1, \mu, \a}_{[0, t_2]}, \nu_{[0, t_2]})\big)  + C[t_2-t_1]\\
&\le& \rho\big(C\cW_2((t_1, \mu), (t_2, \nu)\big)  + C[t_2-t_1].
\eeaa
This implies that $V$ is uniformly continuous in $(t,\mu)$. 
}
\qed

\subsection{A state dependent property}
We conclude this section with the following state dependent property.
\begin{thm}
\label{thm-Markov}
Let Assumption \ref{assum-standing} hold. Assume further that $b, \si, f, g$ are state dependent, namely $(b,\si, f)(t, \o, \mu,a) = (b,\si, f)(t,\o_t, \mu_t, a)$ and $g(\o, \mu) = g(\o_T, \mu_T)$, then $V(t, \mu) = V(t,\mu_t)$ is also state dependent.
\end{thm}
\proof By Lemma \ref{lem-VhatV}, it suffices to show that $V_0(t, \mu) = V_0(t, \nu)$ for all $t, \mu, \nu$ such that $\mu_t = \nu_t$. We proceed in three steps.

{\it Step 1.} First, one may construct $\overline \dbP\in \cP(\mu, \nu)$ such that $\ol \dbP(X_t = X'_t) = 1$. Indeed, one may construct it such that the conditional distributions are independent: for any $\xi,\xi' \in C^0_b(\O)$,
\beaa
\dbE^{\ol \dbP}\Big[\xi(X_{t\wedge \cd}) \xi'(X'_{t\wedge \cd})\Big] := \dbE^{\mu_t}\Big[ \dbE^\mu\big[\xi(X_{t\wedge \cd}) \big| X_t\big] ~ \dbE^\nu\big[\xi'(X'_{t\wedge \cd}) \big| X'_t = X_t\big]\Big].
\eeaa

{\it Step 2.}  For any $\pi: 0=s_0<\cds<s_m= t$ and $\e>0, \d>0$,  we may mimic  the arguments in Lemma \ref{lem-construct}  and construct $(\tilde \O, \tilde \cF, \tilde \dbP, \tilde B, \xi, \eta)$ such that

$\bullet$ $\cL_\xi=\mu, \cL_\eta = \nu$, and $\eta$ is independent of $\tilde B$;

$\bullet$ $\xi_\pi$ is measurable to the $\si$-algebra $\sigma(\eta_\pi,\tilde{B}_{[0,\d]})$.

$\bullet$ $\dbE^{\tilde \dbP}[|\xi_t-\eta_t|^2] \le \e^2$.

\no Indeed, since $\ol \dbP(X_t = X'_t) = 1$, in Cases 1 and 2 in Lemma \ref{lem-VhatV}, it is obvious that $\xi_t = \eta_t$. In Case 3, we can show that $\dbE^{\tilde \dbP}[|\xi_t-\eta_t|^2] \le \e^2$.

{\it Step 3.}  We now mimic the arguments in Lemma \ref{lem-hatVreg} to prove $V_0(t,\mu)=V_0(t,\nu)$.  Fix an arbitrary $\a = \sum_{i=0}^{n-1} h_i(X_\pi, X_{[t, t_i]}) 1_{[t_i, t_{i+1})}\in \cA^0_t$ with the corresponding partition $\pi: 0 \le s_1 <\cds<s_m =t$. Consider the notations in Steps 1 and 2 in this proof, and introduce $\tilde \dbP, B', B^\d, X^\a, \tilde X, \d$ as in Lemma \ref{lem-hatVreg}. Similar to \reff{tildeX2}  we can prove
\beaa
\dbE^{\tilde \dbP}\Big[\sup_{t\le s\le T} |\tilde X_s - X^\a_s|^2\Big] \le C\dbE^{\tilde \dbP}[|\xi_t - \eta_t|^2] + C_\mu \rho_0^2(\d) \le C\e^2 + C_\mu \rho_0^2(\d).
\eeaa
Moreover, following the arguments in \reff{tildeX5} and \reff{DJ}, we can show that 
\beaa
J(t, \mu, \a)  - V_0(t, \nu) &\le& C\dbE^{\tilde \dbP}\Big[\rho_0\big(\sup_{t\le s\le T} [|\tilde X_s - X^\a_s|  + \cW_2(\cL_{\tilde X_s}, \cL_{X^\a_s})]\big)\Big] + C_{\mu,\nu} \rho_0(\d) \\
&\le& C\rho\big(C\e + C_\mu \rho_0(\d)\big) + C_{\mu,\nu} \rho_0(\d),
\eeaa
for some modulus of continuity function $\rho$. Send $\e, \d\to 0$, we obtain: $J(t, \mu, \a)  - V_0(t, \nu)\le 0$. Since $\a$ is arbitrary, this implies that $V_0(t,\mu) \le V_0(t, \nu)$. The opposite inequality can be proved similarly, and thus $V_0(t,\mu) = V_0(t, \nu)$.
\qed


\begin{thebibliography}{99}

\bibitem{BCFP1}
Bandini, E.,  Cosso, A.,  Fuhrman, M., and  Pham H.  {\it Randomization Method and Backward SDEs for Optimal Control of Partially Observed Path-Dependent Stochastic Systems}, {\sl Annals of Applied Probability} , 28 (2018), 1634-1678.

\bibitem{BCFP2}
Bandini, E.,  Cosso, A.,  Fuhrman, M., and  Pham H.  {\it Randomized Filtering and Bellman Equation in Wasserstein Space for Partial Observation Control Problem.}  {\sl Stochastic Processes and their Applications}, 129 (2019), 674-711.

\bibitem{BCP}
Bayraktar, E., Cosso, A., and  Pham, H.  {\it Randomized dynamic programming principle and Feynman-Kac representation for optimal control of McKean-Vlasov dynamics}. {\sl Transactions of the American Mathematical Society}, 370 (2018), 2115-2160.


\bibitem{BFY} 
Bensoussan, A., Frehse, J., and Yam,  S.~C.~P.  {\sl Mean Field Games and Mean Field Type Control Theory}.  Springer Briefs inMathematics, N.Y., Heidelberg, Dordrecht, London, 2013.

\bibitem{BGY}
Bensoussan, A., Graber, P.,  and Yam, S.C.P. {\it Stochastic Control on Space of Random Variables}, preprint, arXiv:1903.12602.



\bibitem{BY}
Bensoussan, A. and Yam, S.C.P. {\it  Control problem on space of random variables and master equation}, {\sl ESAIM: Control, Optimisation and Calculus of Variations},  accepted, arXiv:1508.00713.

\bibitem{BLPR}
Buckdahn, R., Li, J., Peng, S.;, and Rainer, C. {\it Mean-field Stochastic Differential Equations and Associated PDEs}. {\sl Ann. Probab.}, 45 (2017), 824-878.

\bibitem{CHM}
Caines, P.E., Huang, M., and Malhame, R.P.  {\it Large population stochastic dynamic games: closed-loop McKean-Vlasov systems and the Nash certainty equivalence principle}, {\sl Communications in Information and Systems}, 6 (2006), 221-252.

\bibitem{Cardaliaguet}
Cardaliaguet, P. {\it Notes on Mean Field Games (from P.-L. Lions lectures at Coll`ege de France)}, http://www.college-de-france.fr, 2013.

\bibitem{CDLL}
Cardaliaguet, P., Delarue, F., Lasry, J.M., and Lions, P.L. {\it The master equation and the convergence problem in mean field games}, Princeton University Press, 2019.

\bibitem{CD}
Carmona, R. and Delarue, F. {\it The Master Equation for Large Population Equilibriums}, {\sl  Stochastic Analysis and Applications 2014}. Crisan D., Hambly B., Zariphopoulou T. (eds) Springer Proceedings in Mathematics \& Statistics, vol 100. Springer, Cham.


\bibitem{CD1}
Carmona, R. and Delarue, F. {\sl Probabilistic Theory of Mean Field Games I - Mean Field FBSDEs, Control, and Games}. Springer, 2018.

\bibitem{CD2}
 Carmona R. and Delarue, F. {\sl Probabilistic Theory of Mean Field Games II - Mean Field Games with Common Noise and Master Equations}. Springer, 2018.


\bibitem{CCD}
 Chassagneux, J.-F.,  Crisan, D.,  and  Delarue, F. {\it A Probabilistic Approach to Classical Solutions of the Master equation for Large Population Equilibria}. preprint, arXiv: 1411.3009.

\bibitem{CF} Cont, R. and  Fournie, D. {\it Functional It\^o calculus and stochastic integral representation of martingales}. {\sl Ann. Probab.} 41 (2013), 109-133.

\bibitem{CP}
Cosso, A.  and  Pham, H. {\it Zero-sum stochastic differential games of generalized McKean-Vlasov type},  {\sl Journal de Math\'{e}matiques Pures et Appliqu\'{e}es}, accepted, arXiv:1803.07329.

\bibitem{Dupire}
Dupire, B. {\it Functional It\^{o} Calculus}. http://ssrn.com/abstract=1435551, 2009.

\bibitem{EKTZ} 
Ekren, I., Keller, C., Touzi, N., and Zhang, J. {\it On viscosity solutions of path dependent PDEs}. {\sl Annals of Probability}, 42 (2013), 204-236.

\bibitem{ETZ1}
Ekren, I., Touzi, N., and Zhang, J. {\it Viscosity Solutions of Fully Nonlinear Parabolic Path Dependent PDEs: Part I}, {\sl Annals of Probability}, 44 (2016), 1212-1253.

\bibitem{ETZ2}
Ekren, I., Touzi, N., and Zhang, J. {\it Viscosity Solutions of Fully Nonlinear Parabolic Path Dependent PDEs: Part II}, {\sl Annals of Probability}, 44 (2016), 2507-2553.

\bibitem{FGS}
 Fabbri, G.,  Gozzi, F.,  and Swiech, A. {\sl Stochastic Optimal Control in Infinite Dimension -- Dynamic Programming and HJB Equations}, Springer, 2017.
 
\bibitem{GS1}
Gangbo, W.  and Swiech, A. {\it Metric viscosity solutions of Hamilton-Jacobi equations depending on local slopes}. {\sl Calculus of Variations}, 54 (2015), 1183-1218.

\bibitem{GS2}
Gangbo, W. and Swiech A. {\it Existence of a solution to an equation arising from the theory of mean field games}, {\sl Journal of Differential Equations}, 259 (2015), 6573-6643.

\bibitem{GT}
Gangbo, W.  and  Tudorascu, A. {\it On differentiability in the Wasserstein space and well-posedness for Hamilton-Jacobi equations}. {\sl Journal de Math\'{e}matiques Pures et Appliqu\'{e}es}, accepted.


\bibitem{Lacker}
Lacker, D. {\it Limit Theory for Controlled McKean-Vlasov Dynamics}, {\sl SIAM Journal on Control and Optimization}, 55 (2017),  1641-1672.


\bibitem{LL}
Lasry, J. and Lions, P.L. {\it Mean field games}, {\sl Jpn. J. Math.}, 2 (2007), 229-260.

\bibitem{Lions1}
Lions , P.L. {\it Viscosity solutions of fully nonlinear second-order equations and optimal stochastic control in infinite dimensions. I. The case of bounded stochastic evolutions.} {\sl Acta Math.} 161 (1988), 243-278.

\bibitem{Lions2}
Lions , P.L. {\it Viscosity solutions of fully nonlinear second order equations and optimal stochastic control in infinite dimensions. II. Optimal control of Zakai's equation}. In {\sl Stochastic Partial Differential Equations and Applications, II} (Trento, 1988).  Lecture Notes in Math. 1390, 147-170. Springer, Berlin, 1989.

\bibitem{Lions3}
Lions , P.L. {\it Viscosity solutions of fully nonlinear second-order equations and optimal stochastic control in infinite dimensions. III. Uniqueness of viscosity solutions for general second-order equations.} {\sl J. Funct. Anal.} 86 (1989), 1-18.

\bibitem{Lions4} Pierre L. Lions. \emph{Cours au College de France}, www.college-de-france.fr.

\bibitem{MWZ}
Ma, J.,  Wong, L., and Zhang, J. {\it Time Consistent Conditional Expectation under Probability Distortion}, preprint, arXiv:1809.08262.

\bibitem{MZ}
Mou, C. and Zhang, J. {\it Weak Solutions of Mean Field Game Master Equations}, preprint, arXiv:1903.09907.


 \bibitem{PengW}
 Peng, S., Wang, F. {\it BSDE, path-dependent PDE and nonlinear Feynman-Kac formula}. {\sl Sci. China Math.} 59(2016), 19-36.
 
\bibitem{PW}
Pham, H. and Wei, X. {\it Bellman equation and viscosity solutions for mean-field stochastic control problem}, {\sl ESAIM: COCV}, 24 (2018), 437-461.

\bibitem{PZ}
Pham, T., Zhang, J. {\it Two person zero-sum game in weak formulation and path dependent Bellman-Isaacs equation}. {\sl SIAM J. Control Optim.} 52(2014), 2090-2121.

\bibitem{PTZ}
Possamai, D., Touzi, N., and Zhang, J. {\it Zero-sum path-dependent stochastic differential games in weak formulation}, preprint, arXiv:1808.03756.


\bibitem{RR}
Ren, Z. and Rosestolato, M. {\it Viscosity solutions of path-dependent PDEs with randomized time}, preprint, arXiv:1806.07654.

\bibitem{RY}
Revuz, D. and  Yor, M.  {\sl Continuous Martingales and Brownian Motion}, Springer-Verlag, 1991.



\bibitem{SZ}
Saporito, Y. and Zhang, J.  {\it Stochastic Control with Delayed Information and Related Nonlinear Master Equation}, {\sl SIAM J. Control Optim.}, 57 (2019), 693-717.


\bibitem{Sirbu}
Sirbu, M. {\it Stochastic Perron's method and elementary strategies for zero-sum differential games.} {\sl SIAM J. Control Optim.}  52 (2014), no. 3, 1693-1711.

\bibitem{VZ}
  Viens, F.  and  Zhang,  J. {\it A Martingale Approach for Fractional Brownian Motions and Related Path Dependent PDEs}, preprint, arXiv:1712.03637.
  
\bibitem{WZ}
Wu, C.  and  Zhang, J. {\it An Elementary Proof for the Structure of Wasserstein Derivatives}. preprint, arXiv:1705.08046.

\bibitem{Zhang}
Zhang, J. {\sl Backward Stochastic Differential Equations -- from linear to fully nonlinear theory},
Springer, New York, 2017.

\bibitem{W.A.Zheng85} Zheng,  W. A.  {\it Tightness results for laws of diffusion processes application to stochastic mechanics}, {\sl Ann. Inst. Henri Poincar\'{e}}, 1985.

\bibitem{Zhou} 
Zhou, X. Y. (2010). {\it Mathematicalising behavioural finance}. {\sl Proceedings of the International Congress of Mathematicians}. Volume IV 3185-3209. Hindustan Book Agency, New Delhi. 2010.

\end{thebibliography}
\end{document}